\documentclass[11pt]{article}

\usepackage{amsthm,amsfonts}
\usepackage{amsmath,amssymb, enumerate, fullpage}
\usepackage{graphicx,color}
\usepackage[dvipsnames]{xcolor}
\usepackage{hyperref}
\usepackage{multirow}
\usepackage{caption,subcaption}
\usepackage{algorithm}
\usepackage{tabularx}
\usepackage[noend]{algpseudocode}
\usepackage{booktabs,longtable}
\usepackage{algorithm}
\usepackage{algpseudocode}
\usepackage{tikz}
\usetikzlibrary{arrows.meta}

\newtheorem{Proposition}{Proposition}
\newtheorem{corollary}{Corollary}
\newtheorem{definition}{Definition}
\newtheorem{lemma}{Lemma}[subsection]
\newtheorem{theorem}{Theorem}[subsection]

\begin{document}

\title{Identifying changing partial differential equations using \\
Sampled Local WeakIdent}
\author{Wenbo Hao \thanks{School of Mathematics, Georgia Institute of Technology, Atlanta, GA 30332-0160. Email: whao36@gatech.edu},  
Mengyi Tang \thanks{School of Mathematics, Georgia Institute of Technology, Atlanta, GA 30332-0160. Email:tangmengyi@gatech.edu}, and 
Sung Ha Kang\thanks{School of Mathematics, Georgia Institute of Technology, Atlanta, GA 30332-0160. Email: kang@math.gatech.edu. Kang's research is partially supported by Simons Foundation travel grant SFI-MPS-TSM-00014145.}}

\date{}
\maketitle

\begin{abstract} 

We propose Sampled Local WeakIdent (SLW-Ident), a framework for identifying changing governing equations from a single set of given data. Different from a typical approach of using finite element based approximation to represent varying coefficients, this paper explores a local approach in identification of differential equations.  First, we present the power of Local WeakIdent which gives good local identification and is also computationally efficient with a small patch size, yet it can be sensitive to local perturbations. 
We propose SLW-Ident which stabilizes the identification process and also incorporates global information: we first sample patches in the whole given domain, identify equations for each sampled patch, then use residual error of these  equations to find the transitions between different equations. We refer to a region where the support of the identified equation does not change to be a region of one equation.  Within each region of one equation, we pick the most frequently identified equation as the identified equation, and find constant as well as varying coefficient PDEs within each region of one equation. 
This is justified by an uncertainty quantification theory that gives the relation between statistical error of dominant support selection and the number of patches.  We provide various numerical experiments showing that SLW-Ident accurately recovers the regions of one equation and the governing equations for changing PDEs with varying coefficients even with noisy given data.  
\end{abstract}

\section{Introduction}

The problem of identifying differential equations from data is fundamental in science and engineering as witnessed by various early studies such as \cite{akaike1974new,bellman1969new}.  
Many approaches focus on parameter estimation of finding differential equations, e.g., \cite{baake1992fitting,bar1999fitting,bock1983recent,ljung1998system,muller2004parameter}, and there are various related works such as  \cite{bongard2007automated,bongini2017inferring,lu2018nonparametric,schmidt2009distilling}.  
With an assumption that the governing equation is a form of linear combination of linear and nonlinear terms, identification of differential equation can be expressed as a linear system. By using sparsity constraints, recent methods show successful results in identifying the governing equation, e.g.,  Sparse Identification of Nonlinear Dynamics \cite{brunton2016discovering} and IDENT approaches \cite{cheng2025weightedweak,he2023group, kang2021ident,tang2025wgident,tang2023weakident, tang2023fourierident}, see \cite{he2025identreview} for review on IDENT approaches. 
Recent progress using weak formulations has further improved the robustness of differential equation identification, especially in the presence of noisy data \cite{cheng2025weightedweak, Gurevich_Golden_Reinbold_Grigoriev_2024, gurevich2019robust, messenger2021weakpde,messenger2021weakode,reinbold2020using, tang2025wgident,tang2023weakident}.  
For linear and nonlinear feature terms $f_l$s, weak formulation becomes 
\begin{equation} \label{E:weakintegral}
\int_{\Omega} u_t(x, t) \phi(x, t)\,dx\,dt  = \sum^L_{l = 1}c^{(l)}\int_{\Omega} f_l (x, t)\phi(x, t)\,dx\,dt. 
\end{equation} 
with a test function $\phi(x,t)$.  Here $L$ denotes the total number of features in the dictionary, and the feature set $\{ f_l \}_{l=1}^L$ includes all monomials $u^{\beta_l}$, as well as the derivatives of monomials $u^{\beta_l}$, with non-negative integers $\beta_l = 0,1,2,...,\bar{\beta}$, and the derivative orders  $\alpha _l = 1,2,..., \bar{\alpha}$, i.e., $f_l (x,t) =\partial_x^{\alpha_l}\,u^{\beta_l}$.  This feature set is beneficial for Weak formulation, since all the derivatives can be moved to the test function via integration by parts. 
The identification of differential equation becomes finding the coefficient $\mathbf{c}(t) =\bigl(c^{(1)}(t), c^{(2)}(t), \dots, c^{(L)}(t) \bigr)^T$ values.  
In WeakSINDy \cite{messenger2021weakode,messenger2021weakpde} and WeakIdent \cite{tang2023weakident}, a particular test function is used where the parameters are automatically computed from the given data (see \cite{tang2023weakident} Section 2.2 for details). 
The test function has a compact support, but identification in \eqref{E:weakintegral} is carried out using the full spatiotemporal dataset.

In this paper, we explore a notion of locality, and identify differential equations in a local patch.  This local patch is much smaller than the full data region and of the magnitude of the support size of test function.  
This local identification is the  foundation of the present framework: by identifying the governing equation on each small patch independently, we gain the flexibility to identify  and characterize changing dynamics without assuming a fixed global structure. 
There are related works exploring changing or spatially dependent equations. In IDENT \cite{kang2021ident} and GP-IDENT \cite{he2023group}, equations with varying coefficients are modeled through finite element basis representations of the coefficients, that changes in equations are modeled through coefficients that transition between zero and nonzero values. These methods construct a single global system using the entire dataset and recover the coefficients of all candidate terms through a global optimization procedure.  Such approaches are well suited for smoothly varying dynamics, but may require a large number of basis functions when coefficients change abruptly or when there are many governing equations in different regions of the domain. 
Phase-IDENT \cite{yang2026phaseident} adopts a more localized strategy, and explicitly addresses identifying dynamics when it changes once and focus on locating the transition.  
This work complements the proposed work in this paper.  We consider patch-wise identification starting from WeakIdent \cite{tang2023weakident}, and focus on identifying multiple regions of one equation which are the areas where given dynamics can be represented by one differential equation.

\begin{figure}
\centering
\begin{tabular}{cc}
(a) & (b) \\
\parbox[c]{17em}{
\includegraphics[width=0.4\textwidth]{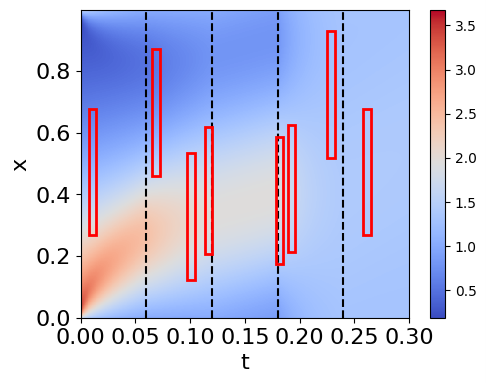}}&
\begin{tabular}{lll}
\toprule
& Interval & Equation \\ \midrule
True & $[0,\;0.060)$  & $u_t=-u_x+0.2\,u_{xx}$ \\
SLW & $[0.001,\;0.048)$ & $u_t=-1.000\,u_x+0.200\,u_{xx}$ \\\midrule
True & $[0.060,\;0.120)$ & $u_t=u_{xx}$\\
SLW & $[0.057,\;0.109)$ & $u_t=1.000\,u_{xx}$ \\ \midrule
True & $[0.120,\;0.180)$ & $u_t=-u_x$\\
SLW & $[0.114,\;0.168)$ & $u_t=-1.001\,u_x$ \\ \midrule
True & $[0.180,\;0.240)$ & $u_t=u_{xx}$\\
SLW & $[0.174,\;0.231)$ & $u_t=0.994\,u_{xx}$ \\ \midrule
True & $[0.240,\;0.300)$ & $u_t=0.1\,u_{xx}-0.75\,(u^2)_x$ \\
SLW & $[0.235,\;0.285)$ & $u_t=0.100\,u_{xx}-0.750\,(u^2)_x$ \\ \bottomrule
\end{tabular}
\end{tabular}
\caption{The proposed SLW-Ident result: (a) Given data $\mathcal{D}$; black dashed lines are the true transition locations, and red rectangles show representative local patches.  (b) SLW-Ident result (SLW) is presented in the second row, compared with the true equations and the boundaries in the first row. }\label{fig:first_image}
\end{figure}

We propose a Sampled Local WeakIdent (SLW-Ident) framework for identifying changing equations and detecting regions of one equation representing different dynamics. We refer to a region where the support of the identified equation does not change to be a region of one equation.   The proposed framework is based on three key ideas: First, we explore the power of  WeakIdent on a local patch.  Even with limited dataset on a local patch, Local WeakIdent results can be quite stable.  We generalize this process with sampling, and identify underlying equations independently on each patch. Second, we use the residual error associated with identified PDEs to locate the boundaries of the regions of one equation.  The behavior of each residual error provides key information on the location of the region of one equation, since it gives a small error when the correct PDE is used to compute the residual error.  Third, although only finitely many patches are sampled, the dominant support converges to the true dominant one with high probability. We provide an uncertainty quantification analysis for patch sampling. 
We present one of the results of the proposed SLW-Ident in Figure~\ref{fig:first_image}.  The proposed method identified the regions of one equation as well as equations in each region well. 
The contributions of this paper are summarized as follows:
\begin{enumerate}
\item{We explore the possibility and the power of WeakIdent for identifying equations in a local patch. Local identification enables us to explore complex dynamics, and gives a good understanding of local behavior. }
\item{Building on the localized identification, we propose Sampled Local WeakIdent (SLW-Ident) for identifying changing partial differential equations.  Local WeakIdent is performed on randomly sampled local patches, which makes identification versatile.  This approach scales naturally to large datasets independently of the global domain size.}
\item{We establish an uncertainty quantification framework for patch-based sampling that yields rigorous, computable confidence bounds for dominant-support selection, providing principled guidance on the number of patches required for statistically reliable identification.}
\end{enumerate}

This paper is organized as follows. In Section~\ref{s: local_identification}, we provide the problem setup and introduce the Local WeakIdent method for conducting identification locally. In Section~\ref{sec: SLW-Ident}, we present the details of the proposed SLW-Ident framework for identifying changing differential equations. In Section~\ref{sec: uq}, we provide the uncertainty quantification analysis for local patch sampling. Section~\ref{sec: algo} discusses implementation details of the algorithms. Numerical experiments are provided in Section \ref{sec:numerical} to validate the proposed framework for identifying changing PDE with varying coefficients.  We conclude the paper in Section~\ref{s: conclusion}.  
We present the proofs of uncertainty quantification results in Appendix \ref{Asec:proofs}, details of the test-function in Appendix \ref{appendix:compute_windows}, and the details of more numerical results in Appendix \ref{append: noise15_stats} and \ref{append:additional_experiments}.

\section{Problem set-up and the effectiveness of Local WeakIdent}
\label{s: local_identification}

Let the given data be 
\[
  \mathcal{D}=\{U_j^n\mid j=0,\dots,N_x;\;n=0,\dots,N_t\},\quad
  U_j^n = u(x_j,t^n)+\epsilon_j^n,
\]
defined on the rectangular domain
$\Omega=[X_1,X_2]\times[0,T]$
with uniform spatial and temporal discretizations
$\Delta x=(X_2-X_1)/N_x$ and
$\Delta t=T/N_t$, where
$x_j=X_1+j\Delta x$,
$t^n=n\Delta t$, and
$N_x$ and $N_t$ denote the numbers of spatial and temporal grid intervals, respectively.
We assume that the observed data $U_j^n$ are obtained from the true solution $u(x,t)$ with a possible independent Gaussian noise, \(\epsilon_j^n \sim \mathcal{N}(0,\sigma^2 I),\)
where $\sigma$ denotes the noise level.

We utilize WeakIdent \cite{tang2023weakident} which identifies a single PDE in the entire spatiotemporal given domain. Using the full dataset, a Fourier transform is applied, and the spectral and temporal cutoff frequencies $k_x^*$ and $k_t^*$ are  identified by fitting a two-piece linear model to the Fourier spectrum along the spatial and temporal frequency axes. The spatial and temporal half-widths of the test function support, denoted by $m_x$ and $m_t$, are then computed by approximating the test function using a Gaussian function and setting $k_x^*$ and $k_t^*$ to the $95\%$ Gaussian tail.  The approximately Gaussian shape of the test function suppresses high-frequency components while preserving the more stable and physically relevant features of the data.  See Appendix~\ref{appendix:compute_windows} for details.  

We apply a simplified WeakIdent only locally in this paper. 
Let $B_x = 3m_x\Delta x$ and $B_t = m_t \Delta t$, and we define the \textit{local patch} as   
\begin{equation}
\label{eq:patch_def}
\mathcal{B}_i = \{(x, t)\mid x_{j_i} \leq x\leq x_{j_i} + B_x, \, t^{n_i} \leq t \leq t^{n_i} + B_t \},
\end{equation}
where $j_i$ and $n_i$ are spatial and temporal indices corresponding to the starting point of patch $\mathcal{B}_i$.  We apply WeakIdent \cite{tang2023weakident} only on  local patch, and do not use high-dynamic-region for the coefficient recovery.  This is a simplified  WeakIdent on a local patch and we refer to this algorithm as the \textit{Local WeakIdent}:  
\begin{equation}\label{eq:localWeakIdent}
 (\mathbf{c}_i, \mathbf{s}_i) = \mathrm{LocalWeakIdent}(\mathcal{B}_i)
\end{equation}
which gives outputs $\mathbf{c}_i \in \mathbb{R}^L$ which is the identified sparse coefficient vector of dimension $L$ (the size of the dictionary), and $\mathbf{s}_i$ is an array of indices where the vector $\mathbf{c}_i$ has non-zero entry.  This support gives the equation form indicating which linear and nonlinear terms are identified.  This is computed from the linear system $ W_i \mathbf{c}_i = b_i $  where $W_i$ and $b_i$ are the weak feature matrix and weak form approximating $u_t$ as in the equation \eqref{E:weakintegral}, defined on the local patch $\mathcal{B}_i$ (not on $\Omega$).

\begin{figure}
\centering
\includegraphics[width=0.7\linewidth]{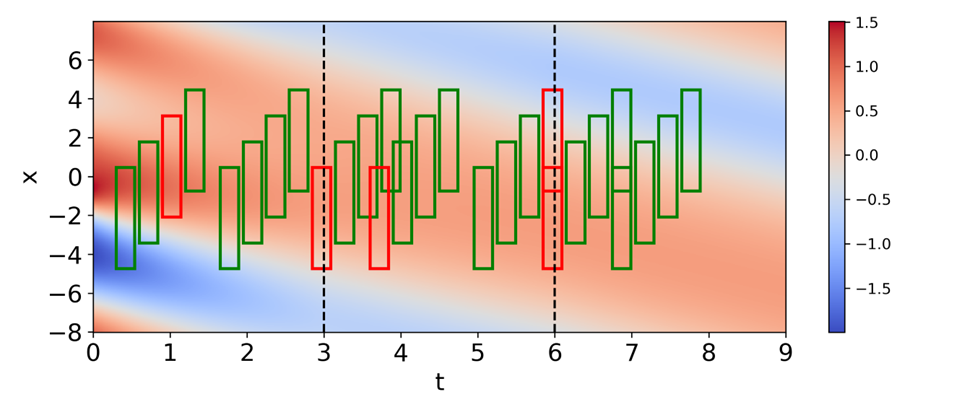}\\
\begin{tabular}{lll}
\hline
Local Patch & True Equation & Identified Equation \\
\hline
$\mathbf{\mathcal{B}_1}$ \textbf{to} $\mathbf{\mathcal{B}_2}$
&
$u_t = 1.0\,u_x + 1.0\,u_{xx}$
&  $\mathbf{u_t = 1.0\,u_x + 1.0\,u_{xx}}$ \\

$\mathcal{B}_3$
& $u_t = 1.0\,u_x + 1.0\,u_{xx}$
&  $u_t = 0.802\,u_{xx} - 2.003\,u_{xxx} + 0.393\,u^3$ \\

$\mathbf{\mathcal{B}_4}$ \textbf{to} $\mathbf{\mathcal{B}_8}$
& $u_t = 1.0\,u_x + 1.0\,u_{xx}$
&  $\mathbf{u_t = 1.0\,u_x + 1.0\,u_{xx}}$ \\
\hline

$\mathcal{B}_9$
& transition
&  $u_t = 1.309\,u_x - 2.144\,u_{xxxx}  - 0.265\,(u^2)_{xx}  - 0.3\,u^3$ \\
\hline

$\mathbf{\mathcal{B}_{10}}, \mathbf{\mathcal{B}_{11}}$
& $u_t = 0.2\,u + u_x + u_{xx} - 0.2\,u^2$
&  $\mathbf{u_t = 0.2\,u + 1.0\, u_x + 1.0\,u_{xx} - 0.2\,u^2}$ \\

$\mathcal{B}_{12}$
& $u_t = 0.2\,u + u_x + u_{xx} - 0.2\,u^2$
& $u_t = g_1$ \\

$\mathbf{\mathcal{B}_{13}}$ \textbf{to} $\mathbf{\mathcal{B}_{19}}$
&$u_t = 0.2\,u + u_x + u_{xx} - 0.2\,u^2$
&  $\mathbf{u_t = 0.2\,u + 1.0\,u_x + 1.0\,u_{xx} - 0.2\,u^2}$ \\
\hline
$\mathcal{B}_{20}$
& transition
& $u_t = g_2$ \\%

$\mathcal{B}_{21}$
& transition
& $u_t = g_3$ \\
\hline
 $\mathbf{\mathcal{B}_{22}}$ \textbf{to} $\mathbf{\mathcal{B}_{28}}$
& $u_t = 1.0\,u_x$
& $\mathbf{u_t = 1.0\,u_x}$ \\
\hline
\end{tabular}
\caption{[Effectiveness of Local WeakIdent] Given data \eqref{e:const feature 3 piece} is shown, with black dashed lines indicating the true  transitions between different equations.  Local patches $\mathcal{B}_i$ (rectangular boxes) are indexed from the left as $\mathcal{B}_1$,  $\mathcal{B}_2$, \dots,  $\mathcal{B}_{28}$. Green patches  represent correct identification, while red patches $\mathcal{B}_3$, $\mathcal{B}_9$, $\mathcal{B}_{12}$, $\mathcal{B}_{20}$ and $\mathcal{B}_{21}$ represent incorrect identification as illustrated in the table. 
Within one region of one equation most identifications are correct, even though the size of the regions used $\mathcal{B}_i$ is much smaller than the given domain. {\small Here $g_1= 0.0778  + 0.754\,u_x  + 1.038\,u_{xx} + 0.249\,(u^2)_x  - 0.548\,u^3  + 0.701\,u^4$, $g_2 = -0.357\,u  + 1.149\,u^2  + 2.856\,(u^2)_x - 0.925\,u^3  - 3.709\,(u^3)_x  + 1.764\,(u^4)_x $ and $g_3 =-2.89  + 2.876\,u  - 20.876\,u_x  - 36.033\,u_{xx}  + 21.542\,(u^2)_x  - 22.711\,(u^2)_{xx}  + 14.755\,u_{xxx}  - 5.555\,u_{xxxx}
$.}}
\label{fig:localW_patch}
\end{figure}

This Local WeakIdent gives good identification of differential equations.  We demonstrate the accuracy and feasibility of Local WeakIdent with an example.  We consider a dynamics with three piecewise constant-coefficient PDE: 
\begin{equation}
    u_t(x, t) = \begin{cases}
        u_x+ u_{xx}, & 0 \le t < 3, \\
        0.2u + u_x + u_{xx} - 0.2u^2,      & 3 \le t < 6, \\
        u_x,                    & 6 \le t \le 9,
    \end{cases}
    \label{e:const feature 3 piece}
\end{equation}
with the initial condition $u_0 = \exp(-(x+1)^2) + \sin {\frac{\pi x}{8}} + \cos{\frac{\pi x}{4}}$, and periodic boundary condition for $t \in [0, 9]$ and $x \in [-8, 8]$.  We set $N_x = N_t = 600$, and the abrupt changes in the governing equation occur at \(t = t^{200}\) and \(t = t^{400}\).  We simulate the given data with \(\Delta x = 0.0267\), \(\Delta t = 0.015\) using a Fourier pseudospectral method for spatial discretization combined with a method-of-lines formulation and adaptive time integration via the LSODA solver. 
In Figure \ref{fig:localW_patch}, we present the given data with various patches $\mathcal{B}_i$ for local identification.  We chose 28 patches at different spatiotemporal locations, labeled $\mathcal{B}_1$ to $\mathcal{B}_{28}$ from left to right.
The detailed identification result for each patch is presented in Table \ref{fig:localW_patch}. Among the 28 patches, there are three patches that cross two regions of one equation, i.e. $\mathcal{B}_9$, $\mathcal{B}_{20}$, and $\mathcal{B}_{21}$, and cannot give meaningful identification results. Among the 25 patches that lie entirely in a single region of one equation, only 2 patches, $\mathcal{B}_3$ and $\mathcal{B}_{12}$, give inaccurate identification results. 
The task of local identification is inherently more challenging than global identification, since there are limited data within a  local patch, e.g., may show insufficient dynamic variations, and it can be more sensitive to measurement noise, yet Local WeakIdent \eqref{eq:localWeakIdent} shows good identification using $B_x = 3m_x\Delta x$ and $B_t = m_t \Delta t$.

\section{Identification of changing partial differential equations using Sampled Local WeakIdent (SLW-Ident)}
\label{sec: SLW-Ident}

Building on the effectiveness of the Local WeakIdent \eqref{eq:localWeakIdent}, we  propose Sampled Local WeakIdent (SLW-Ident) for identifying changing partial differential equations.   This method focuses on finding the \textit{regions of one equation}, and corresponding coefficients vector $\mathbf{c}_i$ and $\mathbf{s}_i$ its support within each region of one equation.   The region of one equation is a subdomain  of $\Omega$ where we identify one equation: this equation can have varying coefficients, but the support of the vector does not change within this region.  We denote each region of one equation as $\Omega_m \subset \Omega$. 

To make the method efficient while incorporating global information, we randomly sample a finite number of local patches to apply Local WeakIdent, instead of considering all possible local patches in $\Omega$.  We explore uncertainty quantification to analyze the effect of such sampling in Section \ref{sec: uq}.   From each sampled local patch,  we record the coefficient and its support $(\mathbf{c}_i, \mathbf{s}_i)$ identified by the Local WeakIdent \eqref{eq:localWeakIdent}.  We use this collected information to find transitions between different PDEs, i.e., regions of one equation, as well as the coefficient $\mathbf{c}_i$ in each region of one equation. 

\begin{figure}
  \centering
  \begin{minipage}[t]{0.3\textwidth}
    \centering
    \includegraphics[width=\linewidth]{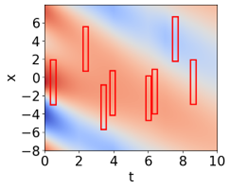}\\[4pt]
    \textbf{Step 1}\\Random patch sampling and Local WeakIdent \eqref{eq:localWeakIdent}.
  \end{minipage}
  \hfill
  \raisebox{2cm}{\begin{tikzpicture}
    \draw[-{Stealth[length=8pt,width=8pt]},line width=3pt] (0,0) -- (0.7,0);
  \end{tikzpicture}}
  \hfill
  \begin{minipage}[t]{0.28\textwidth}
    \centering
    \includegraphics[width=\linewidth]{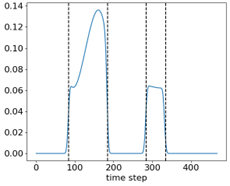}\\[4pt]
    \textbf{Step 2}\\Regions of one equation identification.
  \end{minipage}
  \hfill
  \raisebox{2cm}{\begin{tikzpicture}
    \draw[-{Stealth[length=8pt,width=8pt]},line width=3pt] (0,0) -- (0.7,0);
  \end{tikzpicture}}
  \hfill
  \begin{minipage}[t]{0.26\textwidth}
    \centering
    \includegraphics[width=\linewidth]{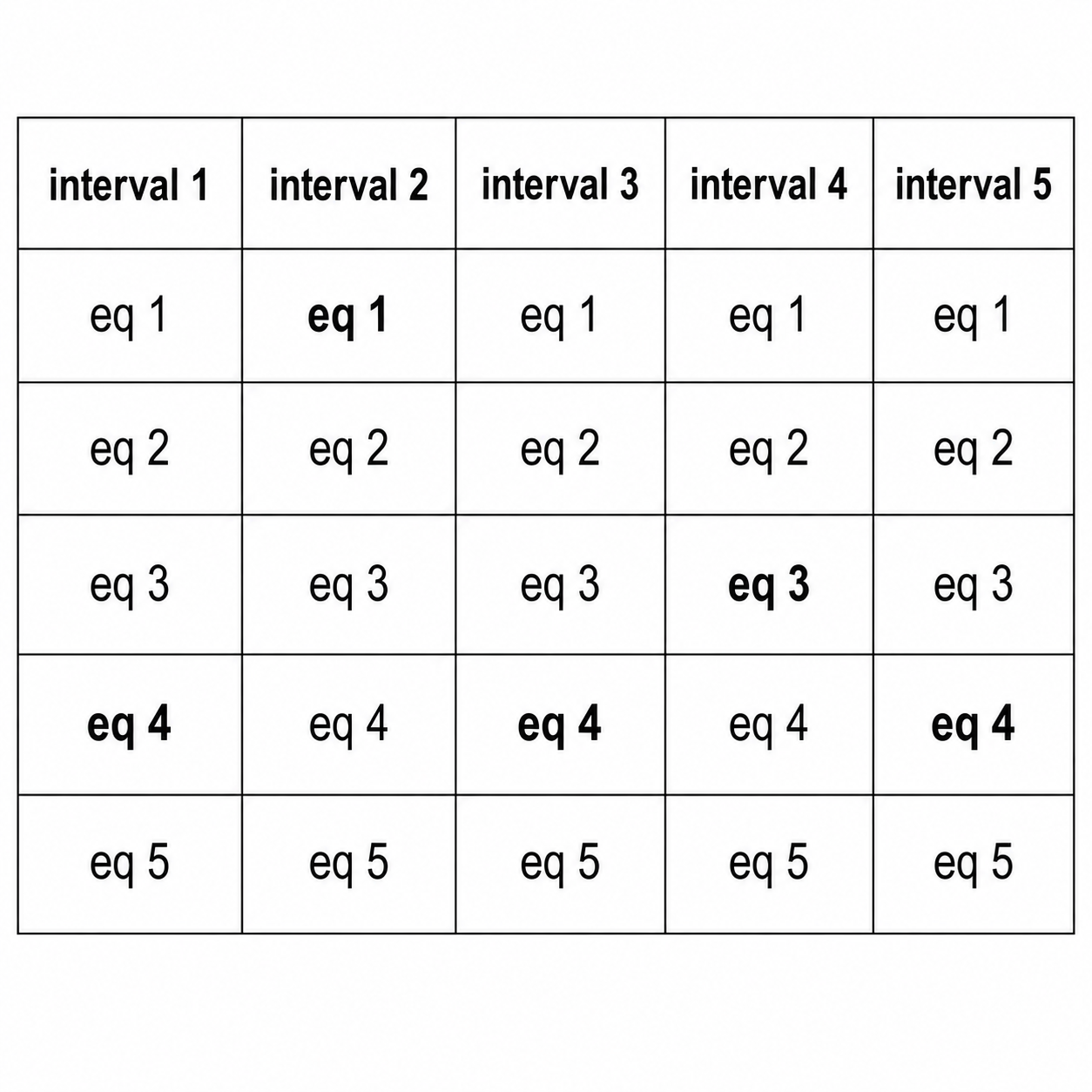}\\[4pt]
    \textbf{Step 3}\\Support selection and coefficient recovery
  \end{minipage}
  \caption{Overview of the SLW-Ident for identifying time-varying changing PDEs. First, Local WeakIdent is applied to randomly sampled local patches. Second, the regions of one equation are  located, and third, coefficient values are computed for each region of one equation.}
  \label{fig:flowchart}
\end{figure}
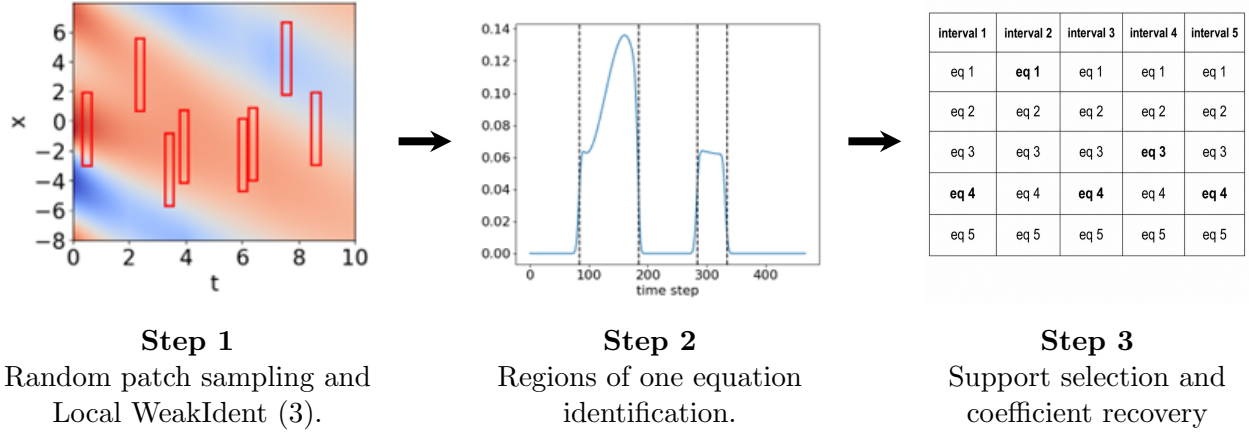

The SLW-Ident framework consists of three main steps as illustrated in Figure~\ref{fig:flowchart}.   First, Local WeakIdent is applied to randomly sampled local patches. Second, the regions of one equation are  located, and third, coefficient values are computed for each region of one equation.  

\vspace{0.3cm}
\noindent\textbf{[Step 1] Patch sampling and Local WeakIdent.} 
We sample local patches $\mathcal{B}_i$ throughout the domain for $i=1,\dots, N_\mathrm{sample}$. Each patch $\mathcal{B}_i$ is defined by its starting indices $(j_i, n_i)$ as in~\eqref{eq:patch_def}, and the sampling is done by choosing the starting indices $(j_i, n_i)$ from the set of admissible starting indices: 
\begin{equation}\label{eq:admissibleSet}
j_i \in \{0, 1, \dots, N_x - 3m_x - 1\} \quad \text{ and }\quad n_i \in \{0, 1, \dots, N_t - m_t - 1\}.
\end{equation}
The upper bounds ensure that no patch extends beyond the given domain. Instead of independent Monte Carlo sampling, we use a structured two-level stratified sampling strategy: first sample in spatial range, then on the selected spatial index, we sample temporal index over the admissible set \eqref{eq:admissibleSet}. 
The sampling procedure is designed to achieve good coverage of the spatiotemporal domain while avoiding duplicate patches. Details of the sampling strategy are provided in Section~\ref{subsec: patch_sampling}. 
Local WeakIdent~\eqref{eq:localWeakIdent} is computed on each local patch $\mathcal{B}_i$, and all identified supports are collected into the \textit{candidate equation set}:
\begin{equation}\label{eq:Scand}
\mathcal{S}_\mathrm{cand}
= \{ \mathbf{s}_i \;:\; \mathbf{s}_i \text{ identified from } \mathcal{B}_i,\; i = 1, 2, \dots, N_\mathrm{sample} \}.
\end{equation}
We denote the size of this set as   $|\mathcal{S}_\mathrm{cand}|= N_\mathrm{cand} \leq N_\mathrm{sample}$, and each $\mathbf{s}_k \in \mathcal{S}_\mathrm{cand}$ represents a unique support (unique equation). 
Since the success rate for local identification is high (as illustrated in Section~\ref{s: local_identification}),  and the structured sampling provides broad domain coverage, the union of identified supports contains the true equation set with high probability. Further numerical details are provided in Section~\ref{sec: algo}.

\vspace{0.3cm}
\noindent\textbf{[Step 2] Regions of one equation: transition detection and interval refining.}
To determine regions of one equation, we consider the collection of residual error of each candidate equation given by $\mathbf{s}_k$.  
For each candidate support $\mathbf{s}_k \in \mathcal{S}_\mathrm{cand}$, we define the \textit{residual error} at $t^n$  as 
\begin{equation}
\label{eqn:residual_def}
r_k(t^n) = \sum_{j=1}^{N_x} \left| u_t(x_j, t^n) - \sum_{l=1}^L c_k^{(l)}(t^n)\,f_l(x_j, t^n) \right|,
\end{equation}
which measures the PDE mismatch at each time step $t^n$. 
Here $\mathbf{c}_k(t^n) = (c_k^{(1)}(t^n), \dots, c_k^{(L)}(t^n))$ is the coefficient vector using the support $\mathbf{s}_k$, recomputed only at each time step $t^n$.  That is, we let $W_k^n$ and $b^n$ be the weak feature matrix and approximation to $u_t$, restricted to the spatial rows at time $t^n$ and feature matrix only containing columns corresponding to the support $\mathbf{s}_k$. Their column-normalized versions are $\widetilde{W}_k^n$ and  $\widetilde{b}^n$.  Then, $\mathbf{c}_k(t^n)$ is computed by the \textit{varying-coefficient regression}:   
\begin{equation}
\label{eqn: time-varying-regres}
\mathbf{c}_k(t^n) = \arg\min_{\mathbf{c}}\|\widetilde{W}_k^n\,\mathbf{c} - \widetilde{b}^n\|_2^2\; .
\end{equation}
This local fitting at $t^n$ allows for more accurate coefficient identification and identification of varying coefficient possible. 
The key idea is that when $\mathbf{s}_k$ is the true governing equation of a region of one equation, the residual $r_k$ remains uniformly low throughout the interval and rises sharply wherever the governing equation differs from $\mathbf{s}_k$.
This is illustrated in Figure \ref{fig:flowchart} [Step 2] the residual error \eqref{eqn:residual_def} gives clear transitions between each region of one equation.   This graph is  $r_k(t^n)$ for one of the supports $\mathbf{s}_k \in \mathcal{S}_\mathrm{cand}$, which shows there are three regions with the equation given by $\mathbf{s}_k$ and there are two other regions which represent different equations. 

To find the transition locations, e.g., the black dotted lines in Figure \ref{fig:flowchart} [Step 2], we propose \textit{grouping transition points}  for more stable identification.
We consider the gradient of $r_k$ to represent the transitions as spikes, and only consider large spikes as the transition points. Then, we collect all transition points from each support $\mathbf{s}_k \in \mathcal{S}_\mathrm{cand}$, and consolidate them using a clustering method to simplify and stabilize locating the transition points. Details are presented in Section~\ref{Subsec:regionofOne}.  
 
This process finds a region of one equation $\Omega_m$, which may be smaller than the true region of one equation as illustrated in Figure \ref{fig:first_image}(b).  This is due to the width of the local patch where the information is mixed.  We refer to such region as a transitional region.  Studies such as Phase-IDENT \cite{yang2026phaseident}  can be applied to further refine the transition region as necessary.  We focus on the global changing equations identification, and the power of local identification in this work. 

\vspace{0.3cm}
\noindent\textbf{[Step 3] Changing differential equation identification: dominant support selection and coefficient recovery.}
Given the regions of one equation $\Omega_1, \Omega_2, \dots, \Omega_M$ and the candidate equation set $\mathcal{S}_\mathrm{cand}$, now we find the governing equation for each region of one equation.  

There are many sampled local patches within $\Omega_m$ in [Step 1].  For each region $\Omega_m$, let $\mathcal{T}(\Omega_m) = \{\, i : \mathcal{B}_i \cap \Omega_m \neq \emptyset \,\}$ denote the set of sampled local patches intersecting $\Omega_m$.  We define the empirical frequency of each candidate support $\mathbf{s}_k$ within $\Omega_m$ as 
\[
\pi_m(\mathbf{s}_k) = \frac{|\{\, i \in \mathcal{T}(\Omega_m) : \mathbf{s}_i = \mathbf{s}_k \,\}|}{|\mathcal{T}(\Omega_m)|}.
\]
This is the ratio of how many times the support $\mathbf{s}_k$ is identified among all sampled local patches in $\Omega_m$. 
We choose the dominating support as the underlying equation within the region of one equation $\Omega_m$:
\begin{equation} \label{eq:dominantS}  
\mathbf{s}^{opt}_m = \arg\max_{\mathbf{s}_k \in \mathcal{S}_\mathrm{cand}} \pi_m(\mathbf{s}_k). 
\end{equation}
This frequency-based selection is theoretically justified by our sampling design and the choice of $N_\mathrm{sample}$ (see Section~\ref{sec: uq}), the most frequently identified support across patches intersecting $\Omega_m$ equals the true dominant support with high probability. 
We construct the coefficients by applying varying-coefficient regression~\eqref{eqn: time-varying-regres} restricted to $\mathbf{s}^{opt}_m$ for the data within $\Omega_m$. Using time-varying rather than fixed-coefficient regression allows the method to capture varying coefficients, while the support remains the same within the region. 

The \textit{identified governing equation} is given in each region of one equation $\Omega_m$ as 
\[
u_t(x,t) = \sum_{l \in \mathbf{s}^{opt}_m} c^{(l)}(t)\, f_l(x,t), \qquad (x, t) \in \Omega_m, \quad m = 1, \dots, M.
\]
We present the details in Section~\ref{sec: algo}. In the following section, we develop an uncertainty quantification theory for local patch sampling, which provides probabilistic guarantees for the dominant-support selection. 

\section{Uncertainty quantification for local patch sampling}
\label{sec: uq}

In the proposed method, we sample local patches \eqref{eq:patch_def} and consider the dominant support \eqref{eq:dominantS} in [Step 1] and [Step 3], since it is neither computationally feasible nor necessary to perform local identification on all admissible patches across the spatiotemporal domain.  To ensure correct identification of the dominant support, we quantify the uncertainty arising from finitely many patches in this section. 

Let $p_k$ denote the empirical probability of identifying support $\mathbf{s}_k \in \mathcal{S}_\mathrm{cand}$ based on the sampled local patches.  Let $q_K$ denote the corresponding probability of identifying support $s_\mathrm{K} \in \mathrm{S}_\mathrm{true}$ under exhaustive sampling over all admissible patches.  The lowercase $k$ and the uppercase $\mathrm{K}$ refer to support from the empirical and true support collection, respectively. Here $\mathrm{S}_\mathrm{true}$ denotes the set of supports with nonzero true probabilities of being identified. Since the empirically observed supports must have nonzero true probability of being identified, $\mathcal{S}_\mathrm{cand} \subset \mathrm{S}_\mathrm{true}.$
For the convenience of analysis, we use the same index to indicate the same support in $\mathcal{S}_\mathrm{cand}$ and $\mathrm{S}_\mathrm{true}$, so that $k = K$ whenever $\mathbf{s}_k = s_\mathrm{K}$, i.e., the supports in $\mathcal{S}_\mathrm{cand}$ are the first $N_\mathrm{cand}$ supports of $\mathrm{S}_\mathrm{true}$.
The collection of probabilities $p_k$ and $q_K$ form an empirical distribution $P = (p_k)_{\mathbf{s}_k \in \mathcal{S}_\mathrm{cand}}$ and true distribution $Q = (q_K)_{s_\mathrm{K} \in \mathrm{S}_\mathrm{true}}$ which are distributions of the set of $\mathcal{S}_\mathrm{cand}$ and $\mathrm{S}_\mathrm{true}$, respectively.  Let $\mathcal{I}$ denote the sample set of local patches.  We consider one region of one equation without any transition, i.e., one of $\Omega_m$, and $\mathcal{S}_\mathrm{cand}$ is only considered for local patches in $\mathcal{T}(\Omega_m)$ as in [Step 3].

Our goal is to quantify the probability that the empirically identified dominant support, $\arg\max_{k} p_k$, coincides with the true dominant support, $\arg\max_{K} q_K$, given the number of patches $|\mathcal{I}| = N_\mathrm{sample}$. 
We analyze the deviation between $p_k$ and $q_K$ when $\mathbf{s}_k= s_\mathrm{K}$, since controlling this deviation preserves the ordering of supports and ensures correct identification of the dominant support. 
We define the \textit{dominance ratio} as
\begin{equation}
\label{dominance_ratio}
\mathrm{R} = \pi_m(\mathbf{s}^{opt}_m).
\end{equation}
using the frequency of the dominating support \eqref{eq:dominantS}.
We consider two different sampling regimes that lead to distinct types of probabilistic guarantees. Under general independent sampling in Subsection~\ref{subsec: max_dev}, we establish bounds based on the maximum deviation $\max_k |p_k - q_k|$, which provide robust worst-case control. Under a Monte-Carlo sampling scheme in Subsection~\ref{subsec: exp_dev}, we further obtain probabilistic guarantees by analyzing the squared deviation $(p_k - q_k)^2$ for each support $\mathbf{s}_k \in \mathcal{S}_\mathrm{cand}$. Together, these results yield explicit guarantees for correctly identifying the dominant support, and in Subsection \ref{s:uq-validation}, we present numerical validation of confidence bounds. 

\subsection{Confidence guarantees under general sampling schemes}
\label{subsec: max_dev}

We consider the general setting in which patches are sampled independently, in which case the maximum deviation $\max_{k} |p_k - q_k|$ provides a natural criterion for quantifying the sampling error.

To get concentration bounds on the maximum deviation, we start from Hoeffding's inequality and
derive a concentration bound on the entropy deviation $|H(P) - H(Q)|$, where $H(P) = -\sum_{k:\ \mathbf{s}_k \in \mathcal{S}_\mathrm{cand}} p_k \log p_k$ and $H(Q) = -\sum_{K: s_\mathrm{K} \in \mathrm{S}_\mathrm{true}} q_K\log q_K$ denote the entropies of the empirical and true distributions, respectively. 

\begin{theorem}\label{thm:max-dev}
Suppose that the $|\mathcal{I}| = N_\mathrm{sample}$ patches are sampled independently from a fixed distribution. Then, the $N_\mathrm{sample}$ identification outcomes are i.i.d random variables, and for any $\eta >0$, there exists a constant $N_{\eta}>0$ such that whenever $|\mathcal{I}|\ge N_{\eta}$, for any $\epsilon>0$,
\[
  \mathbb{P}\!\bigl(|H(P)-H(Q)|\ge\epsilon + \eta\bigr)
  \;\le\;
  2\exp\!\Bigl(-\tfrac{2\,N_\mathrm{sample}\,\epsilon^2}{(\log q_{\min})^2}\Bigr) + \eta.
\]
Here $q_{\min} := \min_{K:\ s_\mathrm{K} \in \mathrm{S}_\mathrm{true}} q_K > 0$ is the smallest nonzero probability of the distribution $Q$.
\end{theorem}
\begin{proof}
See Appendix~\ref{app:proof-thm-entropy-concentration}.
\end{proof}

With a concentration bound for the entropy deviation, we bound the distance between the two probability distributions \(P\) and \(Q\) directly.  We establish a connection between the entropy deviation and the KL divergence in Lemma~\ref{lem:kl-entropy-bound}. 

\begin{lemma}\label{lem:kl-entropy-bound}
Let \(P\) and \(Q\) be two discrete probability distributions such that 
\(\mathrm{supp}(P) \subseteq \mathrm{supp}(Q)\). Then the Kullback--Leibler divergence satisfies
\[
  D_{\mathrm{KL}}(P \parallel Q)
  \;\le\;
  \frac{1}{q_{\min}}\bigl(H(Q) - H(P)\bigr)
  + \frac{1 - q_{\min}}{q_{\min}}\, H(P),
\]
\end{lemma}
\begin{proof}
See Appendix~\ref{app:proof-lem-kl-entropy-bound}.
\end{proof}

Using Lemma~\ref{lem:kl-entropy-bound}, we relate the entropy difference to the total variation distance via Pinsker's inequality, which yields an upper bound on the total variation. This bound, in turn, leads to a concentration inequality for \(\max_k |p_k - q_k|\).

\begin{theorem}\label{thm:max-dev-kl}
In addition to the conditions of Theorem~\ref{thm:max-dev}, assume that $\text{supp}(P) = \text{supp}(Q)$, i.e., the sampled local patches are sufficient to identify all possible outcomes. Then for any $\eta >0$, there exists a constant $N_{\eta}>0$ such that whenever $N_\mathrm{sample}\ge N_{\eta}$, for any $\epsilon>0$, 
\[
  \mathbb{P}\!\Bigl(\max_{k: \ \mathbf{s}_k \in \mathcal{S}_\mathrm{cand}}|p_k - q_k|\le 
    \sqrt{\frac{\epsilon + \eta}{2\,q_{\min}}
      +\frac{1-q_{\min}}{2\,q_{\min}}\,H(P)}\Bigr)
  \;\ge\;
  1-2\exp\!\Bigl(-\tfrac{2\,N_\mathrm{sample}\,\epsilon^2}{(\log q_{min})^2}\Bigr)-\eta.
\]
\end{theorem}
\begin{proof}
See Appendix~\ref{app:proof-thm-max-dev-kl}.
\end{proof}

In practice, the condition \(\mathrm{supp}(P) = \mathrm{supp}(Q)\) is satisfied with high probability when the number of sampled local patches is sufficiently large. Under this assumption, Theorem~\ref{thm:max-dev-kl} provides a high-probability bound on the maximum deviation \(\max_k |p_k - q_k|\).
Let \(k' = \arg\max_k p_k\) and \(\mathrm{R} = p_{k'}\) denote the empirically identified dominant support and its dominance ratio. If the deviation satisfies
\[
\max_k |p_k - q_k| < \mathrm{R} - 0.5,
\]
then it follows that \(q_{k'} > 0.5\), implying that \(k'\) is also the dominant support under the true distribution.
This observation allows us to translate the deviation bound in Theorem~\ref{thm:max-dev-kl} into a probabilistic guarantee for correct dominant support selection, leading to the following corollary.

\begin{corollary}\label{cor:argmax}
In addition to the condition in Theorem~\ref{thm:max-dev-kl}, suppose we observed a dominance ratio \(\mathrm{R} > 0.5\). Then for any $\eta > 0$, there exists $N_{\eta} >0$ such that whenever $N_\mathrm{sample} \geq N_{\eta}$, 
\[
  \mathbb{P}\bigl(\arg\max_k p_k = \arg\max_K q_K \mid P\bigr)
  \;\ge\;
   1 - 2\exp\!\Bigl(-\tfrac{2N_\mathrm{sample}}{(\log q_{\min})^2}\bigl(2(\mathrm{R}-0.5)^2q_{\min}-(1-q_{\min})H(P)\bigr)^2 - \eta\Bigr) - \eta.
\]
\end{corollary}
\begin{proof}
See Appendix~\ref{app:proof-cor-argmax}.
\end{proof}

The corollary shows that, with probability bounded below by the right-hand side, the empirically identified dominant support coincides with the true dominant support. In this case, additional patch sampling is unlikely to change the identification result.

To obtain a computable bound, we approximate \(q_{\min}\) by the empirical lower bound \(\frac{1}{N_\mathrm{sample}} \leq p_{\min}\), and neglect the vanishing parameter \(\eta\), which becomes negligible for sufficiently large sample size \( N_\mathrm{sample} \). This leads to a practical confidence measure that can be evaluated during the sampling process.

\begin{definition}[Hoeffding Confidence]
Assume that the dominance ratio \eqref{dominance_ratio} satisfies \(\mathrm{R} > 0.5\). The \emph{Hoeffding confidence} \(C_H\) is defined as
\begin{equation}\label{def:hoeffding}
 C_H \;=\;
  1 - 2\exp\!\Bigl(-\tfrac{2N_\mathrm{sample}}{(\log N_\mathrm{sample})^2}
  \bigl(2(\mathrm{R}-0.5)^2 p_{\min} - (1-p_{\min})\,H(P)\bigr)^2\Bigr).
\end{equation}
\end{definition}

The quantity \(C_H\) provides a lower bound on the probability that the empirically identified dominant support matches the true dominant support, and thus giving us useful information on the minimum number of patches we need to take to ensure statistically significant support selection. 

\subsection{Confidence guarantees under Monte-Carlo sampling}
\label{subsec: exp_dev}
\newtheorem{theoremII}{Theorem}[subsection]
In this subsection, we consider an idealized Monte Carlo sampling model in which the patch starting indices are sampled i.i.d., according to the uniform distribution on the admissible set of starting indices \eqref{eq:admissibleSet}. Under this stronger assumption, we characterize the distribution of the deviation $p_j - q_j$ using classical Monte Carlo convergence results, e.g., law of large numbers and central limit theorem. This yields the following result. 

\begin{theoremII}\label{thm:variance-clt}
  For each identified support \(\mathbf{s}_k\), the expected square deviation satisfies 
\[
  \mathbb{E}[(p_k - q_k)^2]
  = \frac{q_k(1-q_k)}{N_\mathrm{sample}}.
\]
Moreover, as \(N_\mathrm{sample}\to\infty\), 
\[
  p_k - q_k \;\xrightarrow{d}\;
  \mathcal{N}\!\Bigl(0,\tfrac{q_k(1-q_k)}{N_\mathrm{sample}}\Bigr).
\]
\end{theoremII}
\begin{proof}
See Appendix~\ref{app:proof-thm-variance-clt}. 
\end{proof}

Since we characterized the asymptotic distribution of the deviation $p_k - q_k$, we  use the Gaussian density to derive an empirical lower bound on this deviation for each support, and consequently obtain a confidence measure for the selected dominant support. 

However, the true variance 
$\displaystyle{
\mathrm{Var}(p_k - q_k) = \frac{q_k(1-q_k)}{N_\mathrm{sample}}}
$
is unknown since $q_k$ is not observable. To obtain a computable approximation, we replace it with the sample-based estimate 
$\displaystyle{
\frac{p_k(1-p_k)}{N_\mathrm{sample}-1},
}$
which converges to the true variance as $N_\mathrm{sample}\to\infty$ by the law of large numbers. 
Based on this approximation, we introduce an auxiliary distribution $P'$ that captures the asymptotic behavior of $P$, and derive the following result.

\begin{Proposition}[Monte-Carlo confidence approximation]\label{prop:monte-carlo}
Let $\{p'_k\}_{\mathbf{s}_k\in \mathcal{S}_\mathrm{cand}}$ be random variables approximating the empirical frequencies such that for each $\mathbf{s}_k \in \mathcal{S}_\mathrm{cand}$,
\[
p'_k - q_k \sim \mathcal{N}\!\left(0, \tfrac{p_k(1-p_k)}{N_\mathrm{sample}-1}\right).
\]
Let $\mathrm{R} := \max_{k: \ \mathbf{s}_k \in \mathcal{S}_\mathrm{cand}} p_k$. If $0.5 < \mathrm{R} < 1$, then under the approximation that $\max{p'_k} \approx \max{p_k} = \mathrm{R}$, we have the approximate probabilistic lower bound
\begin{align*}
\mathbb{P}\bigl(\arg\max_k p'_k = \arg\max_K q_K \mid \arg\max_k{p_k} = \arg \max_k{p'_k}\bigr) \\
 \;\gtrsim\; \sqrt{\frac{N_\mathrm{sample}-1}{2\pi\,\mathrm{R}(1-\mathrm{R})}}\,
  \int_{-(\mathrm{R}-0.5)}^{\mathrm{R}-0.5}
    \exp\!\Bigl(-\tfrac{(N_\mathrm{sample}-1)x^2}{2\,\mathrm{R}(1-\mathrm{R})}\Bigr)\,dx.
\end{align*}
\end{Proposition}
\begin{proof}
See Appendix~\ref{app:proof-cor-monte-carlo}.
\end{proof}

By the central limit behavior established above, the distribution of $p_k' - q_k$ is an approximation of $p_k - q_k$, with accuracy improving as the sample size increases. Under sufficiently large $N_\mathrm{sample}$, the quantity on the right-hand side of Proposition~\ref{prop:monte-carlo} provides an accurate and computable approximation of the probabilistic lower bound for the event that the empirically selected dominant support coincides with the true dominant support. Motivated by this interpretation, we introduce the following quantity as a practical measure of confidence associated with the observed dominant support.

\begin{definition}[Monte-Carlo Confidence]
Assume that the dominance ratio \(0.5<\mathrm{R}<1\).  The \emph{Monte–Carlo confidence} \(C_M\) is defined by
\begin{equation}\label{def:monteCarlo}     
  C_M \;=\;
  \sqrt{\frac{N_\mathrm{sample}-1}{2\pi\,\mathrm{R}(1-\mathrm{R})}}\,
  \int_{-(\mathrm{R}-0.5)}^{\mathrm{R}-0.5}
    \exp\!\Bigl(-\tfrac{(N_\mathrm{sample}-1)x^2}{2\,\mathrm{R}(1-\mathrm{R})}\Bigr)\,dx.
\end{equation}
\end{definition}

The Monte--Carlo confidence relies on the assumption of uniform patch sampling, whereas the Hoeffding confidence remains applicable under more general independent sampling schemes.  
These bounds \eqref{def:hoeffding} and  \eqref{def:monteCarlo} are nontrivial and are often substantially above zero. 

\subsection{Numerical validation of confidence bounds}
\label{s:uq-validation}

We present a numerical example to compare the empirical frequency that the identified dominant support
\(
\arg\max_k p_k
\)
coincides with the true dominant support
\(\arg\max_K q_K,
\)
against the theoretical lower bound of $\mathbb{P}\big (\arg\max_k p_k = \arg\max_K q_K|P\big)$ provided by the
Hoeffding confidence \(C_H\) in \eqref{def:hoeffding} and the Monte-Carlo confidence \(C_M\) in \eqref{def:monteCarlo}.

We consider a dataset generated from the following PDE:
\[
  u_t = -u_x + 0.05\,u_{xx},
  \quad x\in[0,1],\ t\in[0,0.3],
  \quad u(x,0)=\exp\bigl(-(x+1)^2\bigr),
\]
discretized on a uniform rectangular grid with
\(N_x=N_t=100\).
We choose the patch sizes
$B_x = m_x\Delta x$ and $B_t = m_t\Delta t$, and sample the starting locations of local patches uniformly from all grid points in
\([0,1-m_x\Delta x)
\times
[0,0.3-m_t\Delta t).
\)
For the true dominant support, we perform exhaustive local
identification over all valid patch locations.
The resulting true dominant support is
\(
\arg\max_K q_K=\{u_x,u_{xx}\},
\)
with the corresponding probability
\(q_{\max}:=\max_K q_K\approx0.717.\)
We conduct \(1000\) independent experiments, and in each experiment, we randomly sample \(8\) patches.  There are \(198\) cases with \(\mathrm{R} \le 0.5\) which are excluded from the analysis, since both confidence measures are defined only when \(\mathrm{R} >0.5\).  
\begin{figure}
\centering
\begin{tabular}{ccc}
(a) & (b) & (c)\\
\includegraphics[width=0.31\textwidth,height=2in]{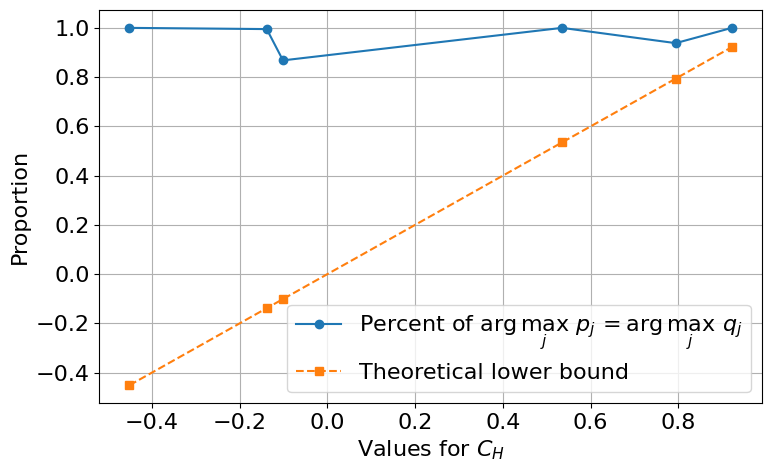} & 
\includegraphics[width=0.31\textwidth,height=2in]{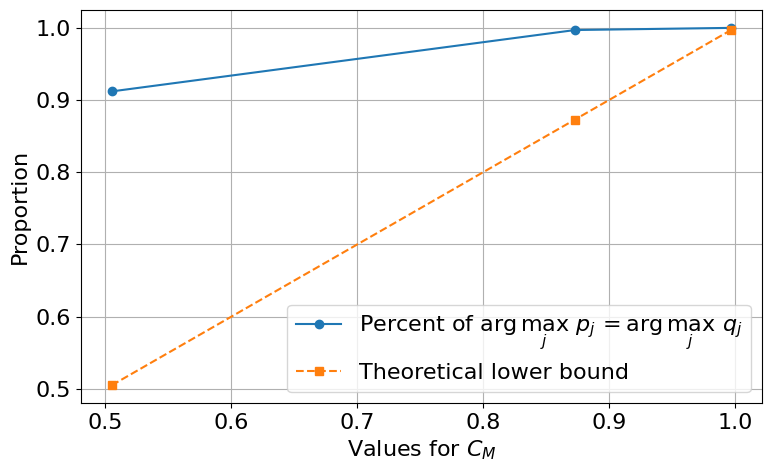} &
\includegraphics[width=0.31\textwidth,height=2in]{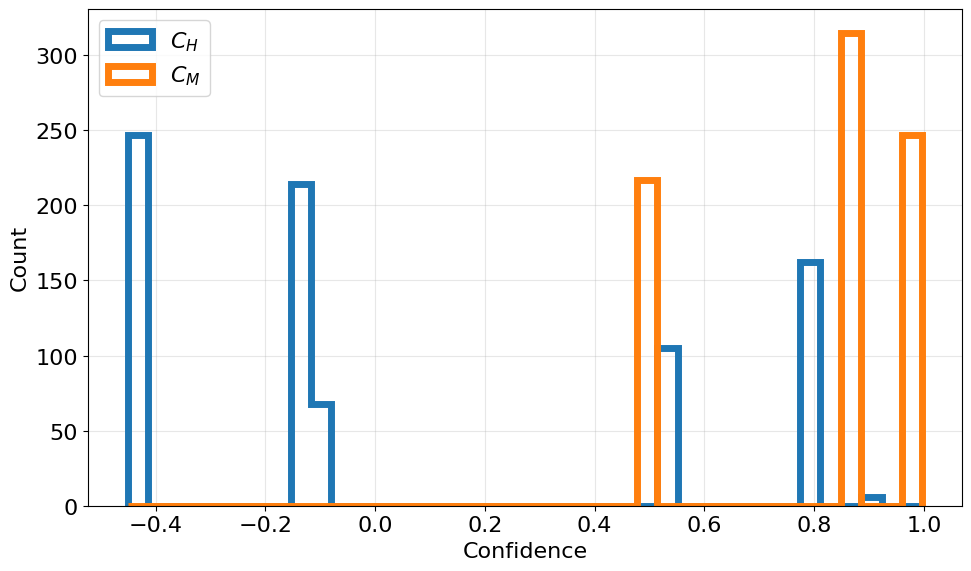}
\end{tabular}
\caption{[Empirical vs.\ theoretical validation of UQ bounds] 
Comparison of the observed probability of
$\arg\max_k p_k = \arg\max_K q_K$
with the theoretical confidence measures:
(a) Hoeffding confidence $C_H$ in \eqref{def:hoeffding};
(b) Monte--Carlo confidence $C_M$ in \eqref{def:monteCarlo};
(c) empirical distributions of $C_H$ and $C_M$. }
\label{fig:exp-validation}
\end{figure}

In Figure \ref{fig:exp-validation}, we show  the percentage of when the dominating support is identical to the true dominating support, i.e., \(\arg\max_k p_k=\arg\max_K q_K\) in blue.  We observed six values of Hoeffding confidence $C_H$ and three values of Monte-Carlo confidence $C_M$. The theoretical lower bounds are presented in orange.  In (a), for the six values of Hoeffding confidence $C_H$, and in (b), for the three values of Monte-Carlo confidence $C_M$.  In (c), we show the empirical distribution of $C_H$ and $C_M$, i.e., the frequency counts for each observed values. 
It shows that the empirical proportions exceed the theoretical lower bounds of both \eqref{def:hoeffding} and \eqref{def:monteCarlo}. 
Under the present sampling setup, the Monte--Carlo confidence is more informative than the Hoeffding confidence. Figure~\ref{fig:exp-validation} (c) shows that the Monte--Carlo confidence values are concentrated at substantially higher levels across the \(1000\) experiments, with the majority exceeding \(80\%\). This is consistently observed in our numerical experiments, and in practice, we use the Monte--Carlo confidence as the primary guideline for determining the number of sampled local patches.  For example, for a single differential equation, i.e. for one region of one equation, we observe that about 100 local patches with dominance ratio \(\mathrm{R}>0.6\) are good. 
This corresponds to Monte--Carlo confidence approximately being 
\[
C_M \approx 0.958,
\]
which is a strong statistical evidence that additional patch sampling is unlikely to alter the selected dominant support.

\section{Algorithms and numerical analysis}
\label{sec: algo}

We present the algorithmic details of SLW-Ident: [Step 1] the local patch sampling strategy in Subsection~\ref{subsec: patch_sampling}, and [Step 2] the procedure for identifying regions of one equation in Subsection~\ref{Subsec:regionofOne}. The complete algorithm is presented in Section~\ref{sec: complete_algo}.

\subsection{[Step 1] Local patch sampling and size of the patch}
\label{subsec: patch_sampling}

\textbf{Patch sampling:} We adopt a two-level stratified sampling strategy: first $N_\mathrm{sample}^x$ number of indices are drawn uniformly without replacement from the spatial  admissible range \eqref{eq:admissibleSet}, and  for each selected spatial index, $N_\mathrm{sample}^t$ number of  temporal starting indices are drawn independently and uniformly without replacement from the temporal   admissible range  \eqref{eq:admissibleSet}, yielding $N_\mathrm{sample} = N_\mathrm{sample}^x \times N_\mathrm{sample}^t$ patches in total.  This structured strategy achieves better coverage of the whole domain than pure Monte-Carlo sampling and achieves at least as good coverage as pure Monte-Carlo sampling. The UQ guarantees of Section~\ref{sec: uq} applies. 

For a typical example, for $N_x = N_t = 600$, and time varying differential equations, we took $N_\mathrm{sample}^x = 20$ and  $N_\mathrm{sample}^t = 40$, and took $N_\mathrm{sample} = 800$ samples. 
These values are chosen to have at least one local patch that lies entirely within one region of one equation.  According to the UQ analysis, choosing $N_\mathrm{sample} = 800$ ensures Monte-Carlo confidence $C_M > 95\%$ when the dominance ratio is over $0.6$, so this number is sufficient to guarantee that the dominant-support selection in [Step~3] is correct with high probability.

\vspace{0.3cm}
\noindent\textbf{Patch size:} 
For the local patch, we choose $B_x = 3m_x\Delta x$ and $B_t = m_t\Delta t$ for time-varying coefficient experiments, e.g., for Figure~\ref{fig:localW_patch}, we used $B_x = 195\Delta x$ and $B_t = 16\Delta t$.  Data with richer spatial or temporal variation encodes more information per grid point, allowing a smaller patch to suffice, while spatially or temporally flat data requires a larger patch. The data-driven parameters $m_x$ and $m_t$ capture this directly, since a larger spectral cutoff $k_x^*$ and $k_t^*$ indicates stronger spatial and temporal variation, respectively. 
The larger spatial patch is used when the governing equation varies only in time, incorporating more spatial data without crossing equation boundaries. 

The local identification should be correct enough times for the dominance ratio \eqref{dominance_ratio} to be $\mathrm{R} > 0.6$  for each region of one equation.  
We present the relation between the patch size and $\mathrm{R}$ with an example. We consider a given data generated from a single-equation 
\begin{equation}\label{eq:oneEquation} 
 u_t = 0.1\,u + 1.0\,u_x + 1.0\,u_{xx} - 0.1\,u^2,
\end{equation}
with an initial condition $u_0(x) = \exp(-(x+1)^2) + \sin\!\bigl(\tfrac{\pi}{8}x\bigr) + \cos\!\bigl(\tfrac{\pi}{4}x\bigr)$ on $x\in[-8,8]$, $t\in[0,5]$, $(N_x, N_t) = (300, 300)$. The data-driven procedure yields $m_x = 35$ and $m_t = 15$. 
In Table~\ref{tab:patch-dom}, we present the dominance ratio  for different patch sizes: time direction ranging from $5 \Delta t$ to $40 \Delta t$ where $m_t\Delta t = 15 \Delta t$, and spatial direction ranging from $20 \Delta x$ to $160 \Delta x$ where $m_x\Delta x = 35 \Delta x$.  Without any transition, larger patches improve the dominance ratio further.  However, if the region of one equation is very small, larger local patch may not lie completely inside the region often which will reduce the accuracy and also increase computational cost. 
For one equation case as in \eqref{eq:oneEquation}, using values around $m_t\Delta t$ and $m_x\Delta x$, already shows $\mathrm{R}$ to be above 60 \%.  Thus, we typically used  $3 m_x\Delta x$ and $m_t \Delta t$ for time-varying PDEs to make space interval a bit larger. 
\begin{table}
\centering
\begin{tabular}{c|cccccccc}
\toprule
$B_x \backslash B_t$ & 5$\Delta t$ & 10$\Delta t$ & 15$\Delta t$ & 20$\Delta t$ & 25$\Delta t$ & 30$\Delta t$ & 35$\Delta t$ & 40$\Delta t$ \\
\midrule
20 $\Delta x$ & 31.8 & 37.8 & 23.1 & \textbf{26.0} & \textbf{56.5} & \textbf{51.6} & \textbf{60.3} & \textbf{55.5} \\
40 $\Delta x$& 34.0 & \textbf{26.1} & \textbf{63.1} & \textbf{60.6} & \textbf{60.4} & \textbf{59.4} & \textbf{58.5} & \textbf{63.6} \\
60 $\Delta x$& 20.5 & \textbf{63.4} & \textbf{65.9} & \textbf{62.3} & \textbf{60.8} & \textbf{65.1} & \textbf{67.9} & \textbf{62.4} \\
80 $\Delta x$& \textbf{27.6} & \textbf{75.9} & \textbf{69.3} & \textbf{69.4} & \textbf{71.9} & \textbf{73.0} & \textbf{70.9} & \textbf{73.6} \\
100$\Delta x$ & \textbf{78.1} & \textbf{74.5} & \textbf{77.4} & \textbf{80.9} & \textbf{76.8} & \textbf{76.5} & \textbf{81.6} & \textbf{81.0} \\
120$\Delta x$& \textbf{85.3} & \textbf{85.4} & \textbf{87.8} & \textbf{86.0} & \textbf{86.9} & \textbf{89.3} & \textbf{89.9} & \textbf{92.8} \\
140$\Delta x$& \textbf{89.9} & \textbf{90.9} & \textbf{92.1} & \textbf{90.9} & \textbf{90.5} & \textbf{93.4} & \textbf{93.8} & \textbf{92.5} \\
160$\Delta x$& \textbf{94.6} & \textbf{91.9} & \textbf{94.3} & \textbf{95.3} & \textbf{94.3} & \textbf{93.1} & \textbf{94.3} & \textbf{96.6} \\
\bottomrule
\end{tabular}
\caption{Dominance ratio $R (\%)$  for different patch sizes. Entries shown in \textbf{bold} indicate correct identification.  Larger region is better for identification, while it may include transition.   We typically used  $3 m_x\Delta x$ and $m_t \Delta t$ for time-varying PDEs. }
\label{tab:patch-dom}
\end{table}

\subsection{[Step 2] Residual-based transition detection and interval refining}
\label{Subsec:regionofOne}

In [Step 1], we sample local patches $\mathcal{B}_i$, find coefficients and their supports $(\mathbf{c}_i, \mathbf{s}_i)$ from  Local WeakIdent \eqref{eq:localWeakIdent} on each $\mathcal{B}_i$, and the candidate equation set $\mathcal{S}_\mathrm{cand} = \{\mathbf{s}_i | i = 1, \dots, N_\mathrm{sample}\}$  \eqref{eq:Scand} is defined.  In  [Step 2], for each candidate support $\mathbf{s}_k \in \mathcal{S}_\mathrm{cand}$, varying-coefficient regression \eqref{eqn: time-varying-regres} is used to compute $\mathbf{c}_k(t^n)$ for each time step $t^n$, which is used to compute the residual error $r_k(t^n)$ in  \eqref{eqn:residual_def}.  This residual error helps to find the transition locations to define regions of one equation, and we present the details here.  We explain for time-varying case, but this approach can be extended to space, or space and time varying cases. 

\begin{figure}
\centering
\begin{tabular}{ccc}
  (a) & (b) & (c) \\
\includegraphics[width=0.3\textwidth]{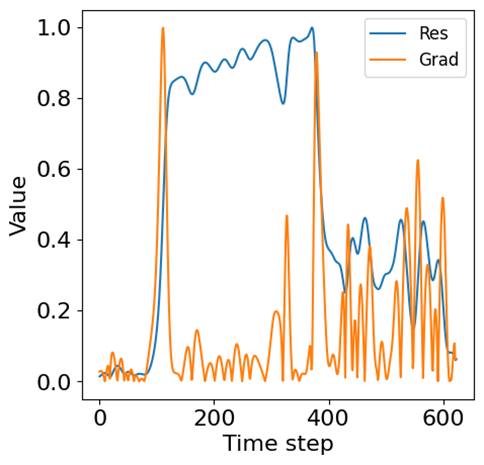} &
\includegraphics[width=0.3\textwidth]{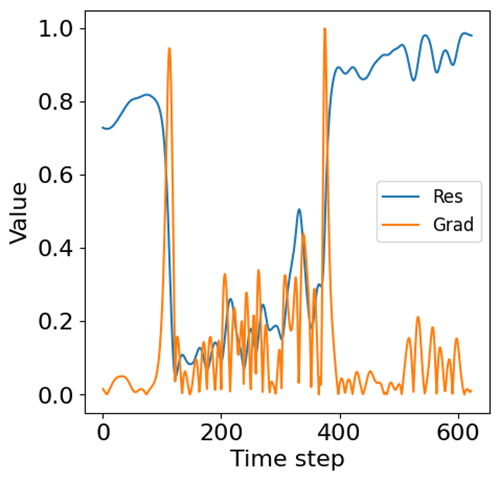}
&
    \includegraphics[width=0.3\textwidth]{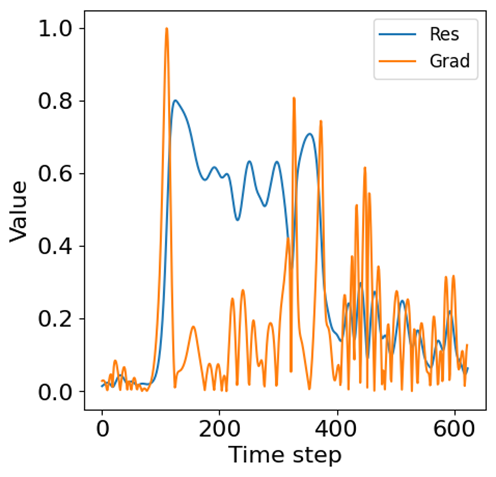}
\end{tabular}
\caption{[An example of $r_k$ and $g_k$] In the graph, $g_k$ (orange) amplifies sudden changes, while $r_k$ (blue) changes gradually.  
(a) $r_1$ and $g_1$ for a support $\mathbf{s}_1$ corresponding to $\{ u_{x}\}$. (b) $r_2$ and $g_2$ for a support $\mathbf{s}_2$ for $\{ u_{xx}, ( u^2)_x\}$.  (c) $r_3$ and  $g_3$ for a support $\mathbf{s}_3$ for $\{ u_{x}, u_{xx}\}$. 
Since the interval end point values $(a_m^k, b_m^k)$ are noisy, we apply clustering \eqref{eq:RegKmeans} to group the regions to find a simpler $\mathcal{P}$. }
\label{fig:transport_burger_noise15_residuals}
\end{figure}

\textbf{Transition point detection:}
From the residual error $r_k(t^n)$ in \eqref{eqn:residual_def}, we compute the discrete temporal gradient by the forward differences
\begin{equation}
\label{eq:gradient_def}
g_k(t^n) = r_k(t^{n+1}) - r_k(t^n).
\end{equation}
While $r_k$ can change gradually near an interval boundary,  $g_k$ amplifies sudden changes and concentrates the transitions into sharp, localized spikes.  Figure~\ref{fig:transport_burger_noise15_residuals} shows an example.

We sort the values of $g_k$ in ascending order and choose a threshold to separate between the meaningful transition points and noise.  
To find this threshold, we use a simple two-piece linear model fit method.  Let $y_0 \le y_1 \le \cdots \le y_{N-1}$ be the input sequence $g_k$ values sorted in ascending order. For each candidate breakpoint $\hat n_b \in \{1,\ldots,N-1\}$, define two line segments that share the point $(\hat n_b-1,\,y_{\hat n_b-1})$:
\begin{align*}
L_1(n) &= y_0 + \frac{y_{\hat n_{b}-1}-y_0}{\hat n_b-1}\,n, \quad n = 0,\ldots,\hat n_b-1,\\
L_2(n) &= y_{\hat n_{b-1}} + \frac{y_{N-1}-y_{\hat n_{b}-1}}{N-\hat n_b}\,\left(n-(\hat n_b-1)\right), \quad n = \hat n_b-1,\ldots,N-1.
\end{align*}
The fitting energy is defined as 
\begin{equation}\label{eq:twoline_fit}
\mathcal{E}(\hat n_b)
=
\left(
\sum_{n=0}^{\hat n_b-1}\bigl(L_1(n)-y_n\bigr)^2
+
\sum_{n=\hat n_b-1}^{N-1}\bigl(L_2(n)-y_n\bigr)^2
\right)^{1/2},
\end{equation}
and we let the optimal minimum breakpoint be  $ n_b = \arg\min_{\hat n_b}\,\mathcal{E}(\hat n_b)$. 
For computation, Python package Numpy's argmin operation is used, which performs a linear scan over the candidate break points.  Candidate breakpoints $\hat n_b$ for which $y_{\hat n_b-1} \le 0.1$ are discarded before the search to reduce unnecessary evaluations. 
From the breakpoint  $n_b$, the threshold is then chosen as the gradient value at the third entry below the breakpoint, i.e., the value corresponding to index $n_b-3$. In all experiments, we use an offset of three, which is empirically found to provide stable and accurate transition detection. 
Time indices $t^n$ for which $g_k(t^n)$ exceeds this threshold are recorded as transition points. 

\textbf{Region of one equation:} 
Regions between transition points $t^n$ define  candidate regions of one equation.  
For each candidate support $\mathbf{s}_k \in \mathcal{S}_\mathrm{cand}$, the residual error $r_k(t^n)$ in  \eqref{eqn:residual_def} is used to find transition points, and there are many regions of one equation $\Omega_m^k$ for $m=1, \dots, M_k$  
\[ 
\mathcal{P}_k = \{\Omega_1^k, \Omega_2^k, \dots, \Omega_{M_k}^k\}.
\] 
We represent each candidate region of one equation
\(
\Omega_m^k=[X_1,X_2]\times[a_m^k,b_m^k)
\)
by the vector
\(
\mathbf{v}_m^k=(a_m^k,b_m^k)\in\mathbb{R}^2.
\)
Typically there are multiple supports identified, and these partitions $\mathcal{P}_k$ are often partial and coarser than the full true partition. 

To stabilize and find meaningful regions of one equation, we propose \textit{grouping transition points}, and apply regularized $k$-means clustering~\cite{kang2011regularized} to the vectors
$\{\mathbf v_m^k:\Omega_m^k\in\mathcal U\}$ where $\mathcal{U} = \bigcup_{k=1}^{N_\mathrm{cand}}\mathcal{P}_k$ :
\begin{equation}
\label{eq:RegKmeans}
\min_{ M,\{{\mathcal C}_m\}_{m=1}^{ M},\{{\boldsymbol\mu}_m\}_{m=1}^{ M}}
\left\{
\lambda
\sum_{m=1}^{ M}
\frac{1}{|{\mathcal C}_m|}
+
\sum_{m=1}^{ M}
\sum_{\mathbf v_m^k\in{\mathcal C}_m}
\|
\mathbf v_m^k-{\boldsymbol\mu}_m
\|_2^2
\right\}.
\end{equation}
Here, $ M$ is the number of clusters to be optimized, ${\mathcal C}_m$ denotes the $m$th cluster, $|{\mathcal C}_m|$ is the number of vectors assigned to ${\mathcal C}_m$, and ${\boldsymbol\mu}_m$ is the corresponding cluster center. Each  representative vector of the minimizer 
\(
\boldsymbol\mu_m=(a_m,b_m)
\)
defines a region
\(
\tilde \Omega_m=[X_1,X_2]\times[a_m,b_m), 
\)
which serves as a candidate for the final computed region of one equation. 
The representative vectors provide principled boundary estimates: since different supports independently detect the same true region with slightly different boundaries due to noise, the representative vector yields more stable boundaries than any individual support prediction.

We define the \textit{region of one equation} 
\[
\mathcal{P} = \{\Omega_1, \Omega_2, \dots, \Omega_M\},
\] 
if it satisfies two conditions: (i) Its temporal length is $b_m - a_m > 3\Delta t$.  (ii) We sort the representative regions $\tilde{\Omega}_m$ in decreasing order of cluster size $|\mathcal{C}_m|$, and starting from the largest cluster, we greedily accept a region if it does not overlap with any previously accepted region.
Large cluster collects interval predictions from many different supports, providing strong multi-source evidence that the corresponding region is a true region of one equation, while a small cluster is more likely to represent a spurious detection from an incorrect candidate support.

Figure~\ref{fig:transport_burger_noise15_residuals} shows an example generated on \(x\in[-2\pi,2\pi]\), \(t\in[0,1]\) with \(N_x=N_t=701\), 
\begin{align*}
&u_t = \bigl(-2 + 0.025\cos t\bigr)\,u_x 
& \text{ in }& \quad [t^0,t^{140}) \\
&u_t = 0.01\,u_{xx} + 0.5\bigl(1 + 0.025\sin t\bigr)\,(u^2)_x 
& \text{ in }& \quad [t^{140},t^{420})\\
& u_t = 0.5\bigl(1 + 0.025\cos t\bigr)\,u_x + 0.02\,u_{xx} 
& \text{ in }& \quad [t^{420},t^{700}),
\end{align*}
with 15\% Gaussian noise (NSR = 0.15), and initial condition 
\[
  u(x,0)
  = \sin(2\pi x)
    +0.5\,\sin\bigl(4\pi x + \tfrac\pi4\bigr)
    +0.25\,\cos\bigl(8\pi x - \tfrac\pi3\bigr)
    +e^{-5(x-0.3)^2}.
\]
Figure~\ref{fig:transport_burger_noise15_residuals} (a) shows $r_1$ (blue) and $g_1$ (orange) for a support $\mathbf{s}_1$ with $\{ u_{x}\}$, (b) $r_2$ (blue) and $g_2$ (orange) for a support $\mathbf{s}_2$ with $\{ u_{xx}, ( u^2)_x\}$ and (c) $r_3$ (blue) and $g_3$ (orange) for a support $\mathbf{s}_3$ with $\{ u_{x}, u_{xx}\}$.  
The residual error $r_k$ is consistently flat and low throughout the interval where $\mathbf{s}_k$ is the correct governing equation, even under 15\% noise and time-varying coefficients, and rises sharply outside it. The gradient $g_k$ makes transitions more visible: spikes in $g_k$ are sharper and better separated from the noise floor than the corresponding variations in $r_k$.
We observe that noise introduces a small systematic shift in the detected transition location.  Thus, clustering step helps to simplify the representation.  
Requiring non-overlapping accepted regions prevents the effect of super long intervals identified by universally wrong supports from polluting the final partition, and the minimum-length condition $b_m - a_m > 3\Delta t$ removes clusters whose centroids correspond to short noise-induced oscillations. 

\subsection{SLW-Ident algorithm}
\label{sec: complete_algo}

\begin{algorithm}
\caption{SLW-Ident: Identification of changing PDEs by Sampled Local WeakIdent}
\label{alg:SLW-Ident}
\begin{algorithmic}[1]

\Require Noisy spatiotemporal data $\mathcal{D}$, dictionary $\mathcal L=\{f_l\}_{l=1}^{L}$, and the number of local patches $N_{\mathrm{sample}}$.
\Ensure Regions of one equation $\{\Omega_m\}_{m=1}^{M}$, supports $\{s_m^\star\}_{m=1}^{M}$, and coefficients $\{\mathbf c_m(t)\}_{m=1}^{M}$.

\Statex \textbf{Step 1: Patch sampling and local identification}
\State Compute data-driven patch parameters $(m_x,m_t)$ (Appendix~\ref{appendix:compute_windows}).
\State Sample $N_{\mathrm{sample}}$ patches $\mathcal{B}_i$ in \eqref{eq:patch_def} with $B_x=3m_x\Delta x$, $B_t=m_t\Delta t$.
\For{$i=1,\dots,N_{\mathrm{sample}}$}
    \State Construct weak features $(b_i,W_i)$ from $\mathcal B_i$, and perform Local WeakIdent~\eqref{eq:localWeakIdent} to obtain $(\mathbf{c}_i,\mathbf{s}_i)$.
    \State Add $\mathbf{s}_i$ to the candidate pool $\mathcal S_{\mathrm{cand}}$.
\EndFor

\Statex \textbf{Step 2: Region of one equation} 
\For{each $\mathbf{s}_k\in\mathcal S_{\mathrm{cand}}$}
    \State Find transition points by \eqref{eq:twoline_fit} and form a candidate partition $\mathcal P_k$ represented by $\{\mathbf v_m^k\}$.
\EndFor

\State Apply regularized $k$-means clustering~\eqref{eq:RegKmeans} to obtain representative vectors $\{\boldsymbol\mu_m\}_{m=1}^{\tilde M}$ from $\mathcal{U} = \bigcup_{k=1}^{N_\mathrm{cand}}\mathcal{P}_k$,  and select non-overlapping regions as the regions of one equation $\mathcal{P} = \{\Omega_1, \dots, \Omega_M\}$.

\Statex \textbf{Step 3: Dominant support selection and coefficient recovery}
\For{$m=1,\dots,M$}
    \State Compute empirical support frequencies $\pi_m(\mathbf{s}_k) = \frac{\#\{i:\mathcal B_i\cap\Omega_m\neq\emptyset,\,\tilde s_i=\mathbf{s}_k\}}{\#\{i:\mathcal B_i\cap\Omega_m\neq\emptyset\}}$.
    \State Select dominant support $s_m^\star = \arg\max_{\mathbf{s}_k}\,\pi_m(\mathbf{s}_k)$.
    \State Recover coefficients on $\Omega_m$ by solving~\eqref{eqn: time-varying-regres} restricted to $s_m^\star$.
\EndFor
\State \Return $\{(\Omega_m,\,s_m^\star,\,\mathbf c_m(t))\}_{m=1}^{M}$.

\end{algorithmic}
\end{algorithm}

We summarize the SLW-Ident in Algorithm~\ref{alg:SLW-Ident}. 
The patch size $(B_x, B_t)$ and sampling count $N_\mathrm{sample}$ are determined by the  procedures in Sections~\ref{subsec: patch_sampling}, the transition threshold and  the residual-gradient approach is given in  Section~\ref{Subsec:regionofOne}, and the sufficiency of $N_\mathrm{sample} = 800$ patches for reliable dominant-support selection is guaranteed by the confidence analysis of Section~\ref{sec: uq}. Together, these choices ensure that each step of SLW-Ident is both theoretically grounded and practically reliable across a range of governing equation structures and noise levels.

\begin{table}
\centering
\begin{tabular}{lll}
\toprule
Step & Cost & Scales with $N_x \cdot N_t$ \\
\midrule
1: Patch sampling \& local WeakIdent 
& $\mathcal{O}(N_\mathrm{sample} \cdot m_x \cdot m_t \cdot L^2)$ & No (constant in $N_x \cdot N_t$) \\
2: Time-varying regression\eqref{eqn: time-varying-regres}  & $\mathcal{O}(N_\mathrm{cand} \cdot N_x \cdot N_t \cdot L^2)$ & Linear \\
2-1: Regularized k-means \eqref{eq:RegKmeans} 
& $\mathcal{O}(N_\mathrm{cand} \cdot M_\mathrm{max})$ & No (negligible) \\
3: 
Coefficient recovery by regression \eqref{eqn: time-varying-regres} & $\mathcal{O}(N_x \cdot N_t \cdot L^2)$ & Linear \\
\midrule
\textbf{Overall} & $\mathcal{O}(N_\mathrm{cand} \cdot N_x \cdot N_t \cdot L^2)$ & \textbf{Linear} \\
\bottomrule
\end{tabular}
\caption{[Computational complexity of SLW-Ident] In [Step1], the main computation is by Local WeakIdent which is proportional to $m_x \cdot m_t$, not the data size $N_x \cdot N_t$. Time varying regression in [Step 2] and [Step 3] is the most computationally heavy, since it is computed for each time step $t^n$.  Here $M_\mathrm{max} = \max\{M_1, \dots, M_{N_\mathrm{cand}}\}$ is the maximum number of intervals per candidate partition, and if the number of clusters is small, computation is relatively fast. }
\label{tab:complexity}
\end{table}
In Table~\ref{tab:complexity}, we present per-step and overall computational complexity of SLW-Ident.  In [Step 1], computational cost is proportional to the patch size $m_x \cdot m_t$ grid points, and is independent of the overall dataset size $N_x \cdot N_t$. The cost of [Step 1] is constant with respect to $N_x \cdot N_t$ and does not grow with the dataset. 
The major cost comes from time-varying regression in [Step 2] and [Step 3], which constructs and fits the feature system for each candidate support via least-squares regression \eqref{eqn: time-varying-regres} for every point $t^n$, and scales linearly in $N_x \cdot N_t$.   
In [Step 2] clustering step, $M_\mathrm{max} = \max\{M_1, \dots, M_{N_\mathrm{cand}}\}$ is the maximum number of clusters per candidate partition. If the number of clusters is small, computation is relatively fast.  The overall complexity of SLW-Ident is therefore $\mathcal{O}(N_\mathrm{cand} \cdot N_x \cdot N_t \cdot L^2)$, linear in the total data size, in contrast to global sparse-regression methods whose cost grows with $N_x \cdot N_t$ through the construction and solution of a large global system.  Furthermore, the $N_\mathrm{sample}$ local identifications in [Step 1] and $N_\mathrm{cand}$ global regressions in [Step 2] are mutually independent and can be parallelized across computing units, providing substantial acceleration in practice.

\section{Numerical experiments}
\label{sec:numerical}

We present detailed experimental results of the SLW-Ident for identifying changing equations. 
All datasets are simulated on uniform grid ranging from \(301\times301\) to \(701\times701\) points. We construct a dictionary consisting of $f_l (x,t) =\partial_x^{\alpha_l}\,u^{\beta_l}$ derivative order \( \alpha_l \le 4\) with polynomial degree \( \beta_l \le 4\), which includes many classical PDEs such as the Burgers equation, Korteweg--de Vries (KdV) equation, Kuramoto--Sivashinsky (KS) equation, Fisher--KPP equation, Allen--Cahn equation, Cahn--Hilliard equation, Swift--Hohenberg equation, and the Schr\"odinger equation.
When we consider only time varying cases, we use notation $I_1,\dots,I_M,$ intervals of one equation, instead of $\Omega_1,\dots,\Omega_M$ for clarity. 

We consider a number of measures to evaluate interval accuracy, support accuracy and coefficient value accuracy, as summarized in Table \ref{tab:metrics}.    
To evaluate region of one equation accuracy, we first define  $I_{\rm SLW}$ is \textit{paired with} a ground-truth interval $I_{\rm gt}$, if the computed region of one equation $I_{\rm SLW}$ overlaps with a ground-truth interval $I_{\rm gt}$ at least 70 \% of its own length:
\[
  |I_{\rm SLW}\cap I_{\rm gt}|
  \;\ge\; 0.7\,|I_{\rm SLW}|.
\]
Region of one equation $I_m$ accuracy is measured by interval True Positive Rate (interval TPR) and interval Positive False Rate (interval PPV). \textit{Interval TPR} is the length of prediction--truth overlap divided by the length of the ground-truth interval, and \textit{interval PPV} is divided by the length of the identified interval.  The interval boundary is quantified by the \textit{Inclusion score}, defined as 
\begin{equation}
\label{eq:match_score}
  \mathrm{Inclusion} =
  \begin{cases}
    1.0, & I_{\rm SLW}\subseteq I_{\rm gt},\\
    0.5, & \text{exactly one endpoint of $I_{\rm SLW}$ lies within $I_{\rm gt}$},\\
    0.0, & \text{otherwise}.
  \end{cases}
\end{equation}

For the support $\mathbf{s}_i$ accuracy, we use the support True Positive Rate (\textit{support TPR}), which is the fraction of terms in the true support that is correctly identified, and the support Positive Predictive Value (\textit{support PPV}), which is the fraction of identified terms that is in the true support.  
To quantify the error in choosing the dominant support introduced by finite sampling, we compute the dominance ratio $\mathrm{R}$ in \eqref{dominance_ratio}, which measures the empirical dominance of the most frequent support within $I_m$. In practice, if we observe \(\mathrm{R}\ge0.6\), this implies over $95\%$ probability that the true dominant support is the one we select, indicating that no additional patch sampling is needed. 

For coefficient $\mathbf{c}_i$ accuracy, we compute three averaged relative errors, $E_2$, $E_\infty$, and $E_{\mathrm{res}}$.   Each error is evaluated at each time $t^n$, then averaged over the overlapping domain $I_{\rm SLW}\cap I_{\rm gt}$ between the identified  interval and its paired ground-truth interval. In Table~\ref{tab:metrics}, $\hat{s}$ denotes the identified support and $s^*$ the true support of the paired ground-truth interval.

\begin{table}
\centering
\begin{tabular}{ll|ll|ll}
\toprule
Interval TPR:  & $\displaystyle\frac{|I_{\rm SLW}\cap I_{\rm gt}|}{|I_{\rm gt}|}$ 
& Interval PPV: & $\displaystyle\frac{|I_{\rm SLW}\cap I_{\rm gt}|}{|I_{\rm SLW}|}$ &
Inclusion:  & Eq. \eqref{eq:match_score}\\
\midrule
Support TPR:   &       $\displaystyle\frac{|\hat{\mathbf{s}}\cap \mathbf{s}^*|}{|\mathbf{s}^*|}$  &
Support PPV: & $\displaystyle\frac{|\hat{\mathbf{s}}\cap \mathbf{s}^*|}{|\hat{\mathbf{s}}|}$    &         
Dom.~ratio $\mathrm{R}$ : & Eq. \eqref{dominance_ratio}  \\ \midrule
$E_2$: averaged & $\displaystyle{\frac{\|\mathbf c-\mathbf c^*\|_2}{\|\mathbf c^*\|_2}}$ &
$E_\infty$:  averaged &  $\displaystyle\max_l\frac{|c^{(l)}-(c^*)^{(l)}|}{|(c^*)^{(l)}|}$  &    
$E_{\mathrm{res}}$: averaged &  $\displaystyle\frac{\|W\mathbf c-b\|_2}{\|b\|_2}$ \\ \bottomrule  
\end{tabular}
\caption{ [Accuracy Measures]  The top row for  interval accuracy, the middle row for support accuracy, and the bottom row for coefficient value accuracy measures.  For coefficient error, relative $\ell_2$, $\ell_\infty$ and weak-form residual error are computed at each time $t^n$, then averaged over the paired overlapping interval $I_{\rm SLW}\cap I_{\rm gt}$. }
\label{tab:metrics}
\end{table}
 
We examine the robustness of the algorithm under noise by testing a varying-coefficient PDE with noise-to-signal ratio (NSR) up to $15\%$, where
\[
  \mathrm{NSR}
  = \frac{\|\epsilon\|_2}{\max_{j,n}U_j^n - \min_{j,n}U_j^n},
  \quad \epsilon_j^n\sim\mathcal{N}(0,\sigma^2).
\]

In the following, we consider various cases to present the effectiveness of the proposed model. We first consider when PDE is changing between constant coefficient equations in Subsection \ref{subsec:5pc_toggle}.  In Subsection \ref{subsec:3pc_mixture}, we consider identifying  PDE which is changing among different PDEs with time varying coefficients, and in Subsection \ref{subsec:3pv_noise15}, we present results when observed data is noisy.  In Subsection \ref{subsec:global_vs_local}, we present comparison with GP-IDENT which is one of the best methods for identifying PDE with varying coefficient using global finite element set-up.  In Subsection \ref{subsec:localization_space_time}, we present results for PDE with space and time varying coefficients. Additional details and results are presented in Appendix~\ref{append: noise15_stats} and \ref{append:additional_experiments}.

\subsection{PDE changing between constant coefficient equations}
\label{subsec:5pc_toggle}

We first consider a challenging scenario in which PDE changes between two equations across five intervals.  The given data  $\mathcal{D}$ is generated by solving an advective Fisher--KPP equation in which the reaction term is intermittently activated:
\begin{equation}\label{exp:Fisher--KPP}
  u_t = c(t)\,u + u_x + u_{xx} - c(t)\,u^2
  \quad\text{on}\quad x\in[-8,8],\;t\in[0,10),
\end{equation}
with $c(t)=0$ when $ t \in [t^0, t^{100})$, activated to 0.2 for $t \in [t^{100},t^{200})$ and 0.3 for $ t \in [t^{300},t^{350})$, otherwise, $c(t)=0$.
We use periodic boundary conditions with initial condition
\[
  u(x,0)
  = \exp \bigl(-(x+1)^2\bigr)
    +\sin \bigl(\tfrac{\pi x}{8}\bigr)
    +\cos \bigl(\tfrac{\pi x}{4}\bigr).
\]
The numerical resolution is set to \(N_x \times N_t = 501 \times 501\).

\begin{figure}
\centering
\begin{tabular}{cc}
(a) & (b) \\
\includegraphics[height=4.5cm]{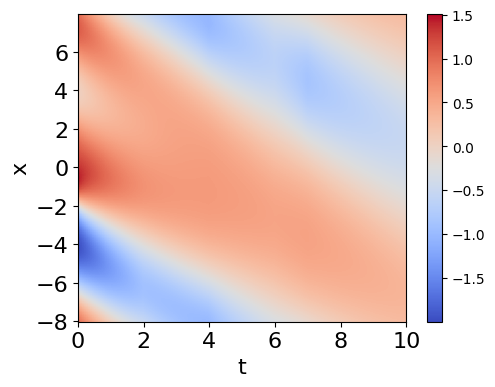} &  \includegraphics[height=4.5cm]{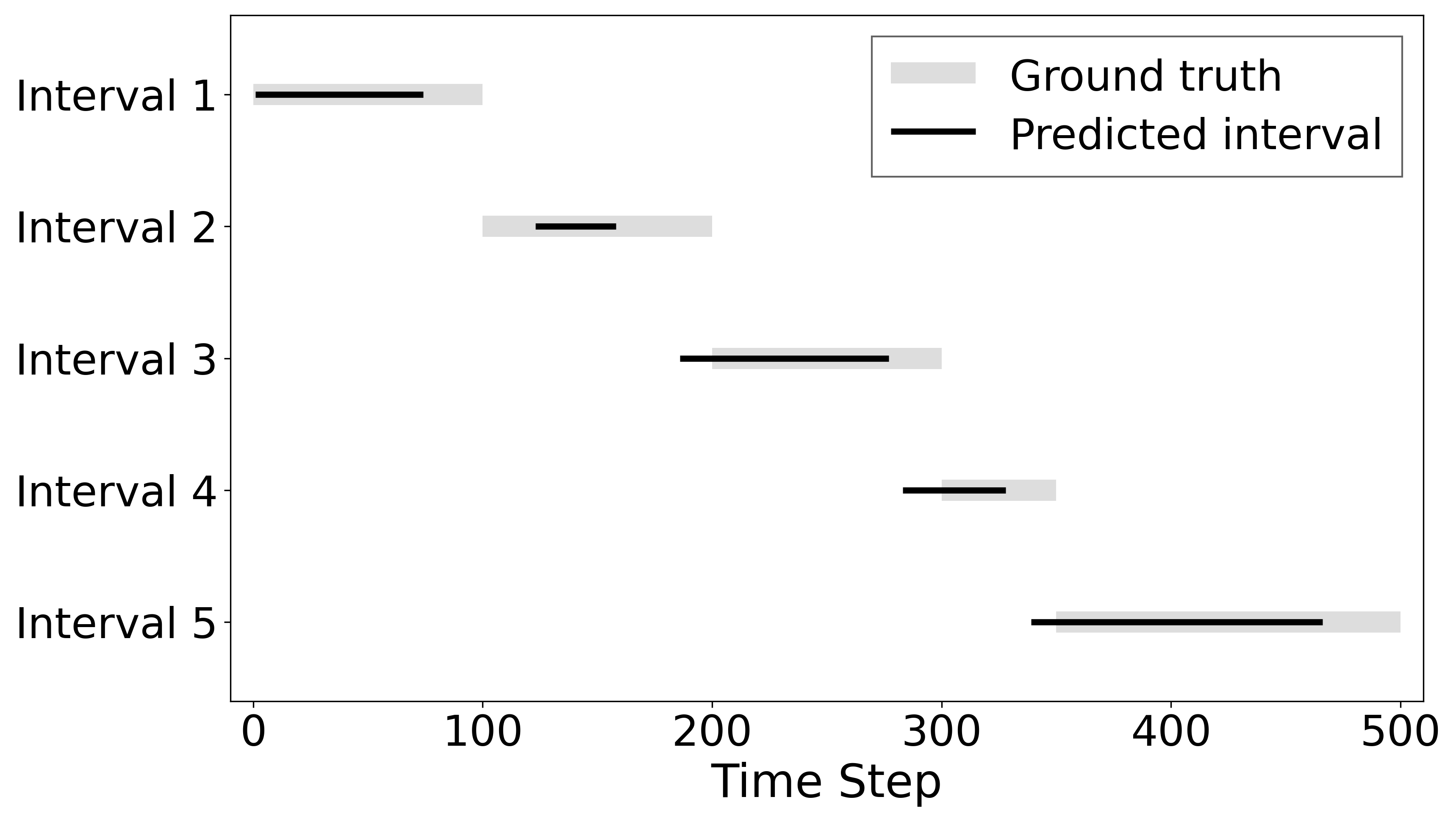}
\end{tabular}
{\small
\begin{tabular}{ll|ll}
\toprule
True Int. & SLW Int. & True Equation & SLW result \\ \midrule
$[t^{0},t^{100})$ & $[t^{1},t^{74})$ & $ u_t = u_x + u_{xx}$ & $u_t = 0.99\,u_x + 0.99\,u_{xx}$ \\
$[t^{100},t^{200})$ & $[t^{123},t^{158})$ & $ u_t = 0.2 u + u_x + u_{xx} - 0.2 u^2$ & $u_t = 0.20 u + 1.00\,u_x + 0.99\,u_{xx} - 0.20 u^2 $ \\
$[t^{200},t^{300})$ & $[t^{186},t^{277})$ & $ u_t = u_x + u_{xx}$ & $u_t = 1.00\,u_x + 0.99\,u_{xx}$ \\
$[t^{300},t^{350})$ & $[t^{283},t^{328})$ & $ u_t = 0.3 u + u_x + u_{xx} -0.3 u^2$ & $u_t = 0.29 u + 1.00\,u_x + 1.00\,u_{xx} -0.29 u^2$ \\
$[t^{350},t^{500})$ & $[t^{339},t^{466})$ & $ u_t = u_x + u_{xx}$ & $u_t = 1.00\,u_x + 0.99\,u_{xx}$ \\\bottomrule
\end{tabular}\\
\begin{tabular}{ccccccccc}
\toprule
Interval & Interval & Inclusion & Support & Support & $\mathrm{R (\%)}$ & $E_2$ & $E_\infty$ & $E_{\mathrm{res}}$ \\
TPR & PPV &  & TPR & PPV & & & & \\ \midrule
0.74 & 1.00 & 1.00 & 1.00 & 1.00 & 83.11 & $1.76 \times10^{-5}$ &  $2.44 \times10^{-5}$ &  $1.07 \times10^{-5}$ \\
0.35 & 1.00 & 1.0 & 1.00 & 1.00 & 96.61 & $8.47\times10^{-8}$ & $4.34\times10^{-7}$ & $8.52\times10^{-8}$ \\
0.77 & 0.85 & 0.5 & 1.00 & 1.00 & 80.35 & $2.97\times10^{-3}$ & $4.19\times10^{-3}$  & $1.31\times10^{-3}$ \\
0.56 & 0.62 & 0.5 & 1.00 & 1.00 & 44.30 & $9.26\times10^{-3}$  & $3.16\times10^{-2}$  & $3.78\times10^{-4}$  \\
0.77 & 0.91 & 0.5 & 1.00 & 1.00 & 99.47 & $4.58\times10^{-4}$  & $6.46\times10^{-4}$  & $1.05\times10^{-4}$ \\ \bottomrule
\end{tabular}
}
\caption{PDE changing between two constant coefficient equations \eqref{exp:Fisher--KPP} (a) Given data.  (b) Five identified intervals compared to the true interval. The middle table shows the SLW-Ident result: three digits for the coefficients are presented, but the computation is done for $10^{-5}$.  The bottom table shows the accuracy measures for Interval 1 to 5 from top to bottom. }
\label{fig:afkpp}
\end{figure}
 
Between successive intervals, the supports differ by one reaction term, making the structural variation minimal and difficult to detect. The proposed SLW-Ident can identify such subtle changes in the governing equation.  
Figure~\ref{fig:afkpp} (a) displays the given data $\mathcal{D}$, which shows how difficult it is to visually see the underlying structural changes in the governing equation.  
Figure~\ref{fig:afkpp} (b) shows the identified regions of one equation align closely with the ground truth.
The interval TPR in the first column is greater than 0.7 for all but Interval 2 and 4, and the interval PPV in the second column exceeds 0.8 for all but Interval 4.  For Interval 2, the magnitude of the reaction term is small and the dynamics is less distinguishable, so that the identified segment became usually short, which contributed to a low interval TPR value. Interval 4 shows a slight shift in the identified region.  Even a small boundary offset has a noticeable impact on both interval TPR and PPV. 
Despite these boundary discrepancies, the support identification is accurate: all supports are recovered perfectly, with TPR = PPV = 1.00 across all detected intervals. The dominance ratios exceed 60\% in all but Interval 4, since the true Interval 4 is half the size of other regions, making it difficult to sample enough patches. 
The coefficient errors mostly stay at the \(10^{-3}\) level or lower.  

\begin{figure}
\centering
\begin{tabular}{cc}
(a) & (b) \\
\includegraphics[height=4.5cm]{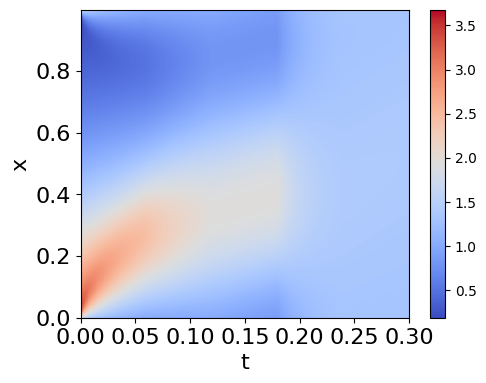} &  \includegraphics[height=4.5cm]{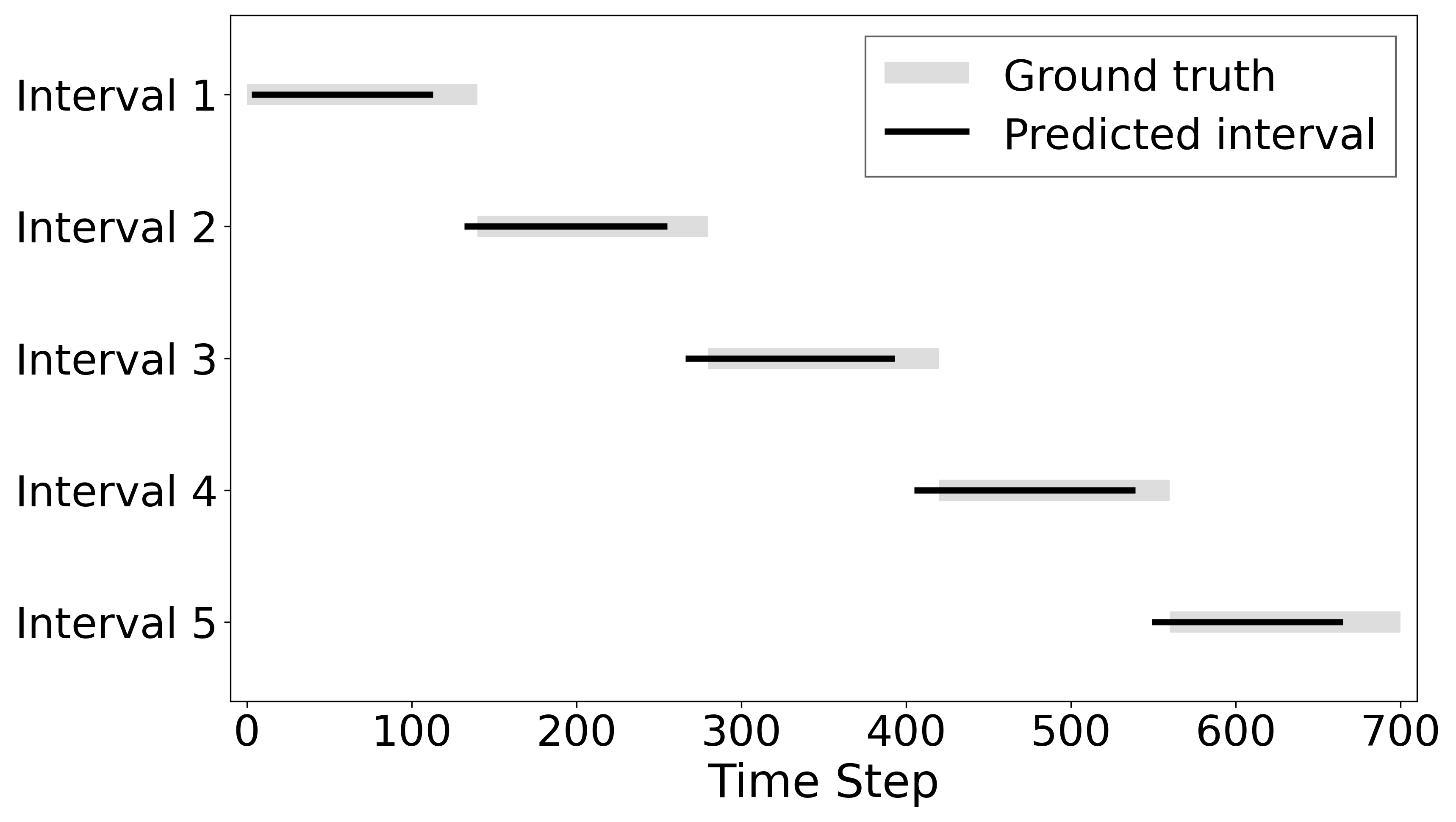}
\end{tabular}
{\small
\begin{tabular}{ll|ll}
\toprule
True Int. & SLW Int. & True Equation & SLW result \\ \midrule
$[t^0,t^{140})$  & $[t^{3},t^{113})$ & $u_t = -u_x + 0.2\,u_{xx}$ & $u_t = -1.00001\,u_x + 0.19998\,u_{xx}$ \\
$[t^{140},t^{280})$ & $[t^{132},t^{255})$& $u_t = u_{xx}$ & $u_t = 1.00000\,u_{xx}$ \\
$[t^{280},t^{420})$ & $[t^{266},t^{393})$ & $u_t = -u_x$ &  $u_t = -1.00093\,u_x$ \\
$[t^{420},t^{560})$ & $[t^{405},t^{539})$ & $u_t = \,u_{xx}$ & $u_t = 0.99360\,u_{xx}$\\
$[t^{560},t^{700})$ & $[t^{549},t^{665})$ & $u_t = 0.1\,u_{xx}-0.75\,(u^2)_x$ & $u_t =  0.10025\,u_{xx} -0.74982\,(u^2)_x$ \\ \bottomrule
\end{tabular}\\
\begin{tabular}{ccccccccc}
\toprule
Interval & Interval & Inclusion & Support & Support & $\mathrm{R (\%)}$ & $E_2$ & $E_\infty$ & $E_{\mathrm{res}}$ \\
TPR & PPV &  & TPR & PPV & & & & \\ \midrule
0.79 & 1.00 & 1.0 & 1.00 & 1.00 & 90.54 & $4.37\times10^{-6}$ & $8.18\times10^{-6}$ & $1.40\times10^{-6}$ \\
0.80 & 0.91 & 0.5 & 1.00 & 1.00 & 85.14 & $3.41\times10^{-5}$ & $1.26\times10^{-4}$ & $4.95\times10^{-6}$  \\
0.81 & 0.89 & 0.5 & 1.00 & 1.00 & 95.24 & $9.30\times10^{-4}$ & $9.30\times10^{-4}$ & $3.56\times10^{-3}$  \\
0.85 & 0.89 & 0.5 & 1.00 & 1.00 & 79.04 & $6.40\times10^{-3}$ & $6.40\times10^{-3}$ & $1.87\times10^{-3}$ \\
0.91 & 0.89 & 0.5 & 1.00 & 1.00  & 58.68  & $4.06\times10^{-4}$ & $2.48\times10^{-3}$ & $3.30\times10^{-6}$\\
\bottomrule
\end{tabular}
}
\caption{PDE changing between five constant coefficient equations. (a) Given data.  (b) Identified five intervals  compared against the true interval. The middle table shows the SLW-Ident result, and the bottom table shows the accuracy measures for Interval 1 to 5 from top to bottom.}
\label{fig:mix1}
\end{figure}

Next example in Figure \ref{fig:mix1}, we consider a dataset $\mathcal{D}$ when the governing equation changes between different constant coefficient equations over five  intervals.  The computational domain is $x\in[0,1]$ and $t\in[0,0.3)$ with $N_x=N_t=701$, and we use initial condition $u(x,0)=10e^{-(x+1)^2}$.  Interval TPR exceeds 0.7 in four of the five Intervals and Interval PPV above 0.7 in all cases. 
The method successfully identifies the governing PDE in each segment with perfect support recovery, Support TPR = Support PPV = 1.00 across all intervals.  The dominance ratios $\mathrm{R}$ are generally high, confirming the stability of support selection within each detected interval. Coefficient estimation errors remain small, with average $E_2$ on the order of $10^{-3}$ or below, indicating accurate recovery.  
The proposed method can localize and identify multiple structurally distinct PDEs within a single dataset, highlighting its robustness to abrupt changes in governing dynamics.

For the following example in Figure \ref{fig:linear}, we consider a dataset $\mathcal{D}$ when a coefficient changes linearly during short transitions between constant values.  We consider the dataset $\mathcal{D}$ given by 
\[
  u_t = 0.3\,u_{xx} - c(t)\,(u^2)_x
  \quad\text{on}\quad x\in[-5,5],\; t\in[0,10),
\]
with $c(t) = 1$ for $t^{0} \le t < t^{60}$, a short transition to zero by $c(t) = -5\,(t-t^{60})+1$ when  $t^{60} \le t < t^{66}$, and $c(t)=0$ for $t^{66} \le t < t^{150}$.  Another short transition $c(t) = 4\,(t-t^{150})$, when $ t^{150} \le t < t^{165}$ and $ c(t) = 2.0$ for $t^{165} \le t < t^{300}$. We use  periodic boundary conditions and initial condition to be 
\[
  u(x,0) = \sin \bigl(\tfrac{\pi x}{4}\bigr),
\]
and \(N_x \times N_t = 301 \times 301\).
Since the transition intervals are short, we take the ground-truth intervals of one equation to be
$  [t^0,t^{60}) : \{u_{xx}, (u^2)_x\}, \quad
  [t^{66},t^{150}) : \{u_{xx}\}, \quad
  [t^{165},t^{300}) : \{u_{xx}, (u^2)_x\}$.
\begin{figure}
\centering
{\small
\begin{tabular}{ll|ll}
\toprule
True Int. & SLW Int. & True Equation & SLW result \\ \midrule
$[t^{0},t^{60})$ & $[t^{0},t^{40})$  & $u_t =  0.3\,u_{xx} -1.0\,(u^2)_x$ & $u_t = 0.30000\,u_{xx}  -0.99976\,(u^2)_x$\\
$[t^{66},t^{150})$ & $[t^{57},t^{129})$ & $u_t = 0.3\,u_{xx}$ & $u_t = 0.29999\,u_{xx}$ \\
$[t^{165},t^{300})$ & $[t^{164},t^{270})$ & $u_t =  0.3\,u_{xx} -2.0(u^2)_x$ & $u_t =  0.30000\,u_{xx} -2.00000(u^2)_x$\\ \bottomrule
\end{tabular}\\
\begin{tabular}{ccccccccc}
\toprule
Interval & Interval & Inclusion & Support & Support & $\mathrm{R (\%)}$ & $E_2$ & $E_\infty$ & $E_{\mathrm{res}}$ \\
TPR & PPV &  & TPR & PPV & & & & \\ \midrule
0.67 & 1.00 & 1.0 & 1.00 & 1.00 & 63.08 & $1.99\times10^{-4}$  & $2.33\times10^{-4}$  & $1.40\times10^{-6}$  \\
0.75 & 0.88 & 0.5 & 1.00 & 1.00 & 86.99 & $1.40\times10^{-5}$  & $2.71\times10^{-5}$  & $5.13\times10^{-5}$ \\
0.78 & 0.99 & 0.5 & 1.00 & 1.00 & 99.32 & $2.66\times10^{-7}$  & $3.31\times10^{-7}$  & $2.10\times10^{-7}$  \\ \bottomrule
\end{tabular}
}
\caption{PDE changing linearly during short transitions between constant coefficient equations. Notice the gap between the true intervals, e.g., $t^{60}$ and $t^{66}$. Recovery of equations is accurate. }
\label{fig:linear}
\end{figure}
SLW-Ident localizes all three intervals of one equation with interval TPR greater than 0.6 and interval PPV greater than 0.8 even with the presence of short linear transition regions. 
The reduction in interval TPR values is due to identified intervals being shorter and shifted. Since the coefficient varies continuously, transition point detection becomes more difficult. The algorithm may either detect the transition slightly earlier or later than the nominal ground-truth boundary. However, these boundary offsets do not affect support identification or coefficient recovery within the intervals of one equation. In all detected intervals, the true support is recovered exactly (TPR = PPV = 1.00), and the estimated coefficients remain highly accurate, with relative errors at or below the $10^{-4}$. 
Additional experiments are provided in Appendix~\ref{append:additional_experiments}. 

\subsection{PDE which is changing among different PDEs with time varying coefficients}
\label{subsec:3pc_mixture}

\begin{figure}
\centering
\begin{tabular}{cc}
(a) & (b) \\
\includegraphics[height=4.5cm]{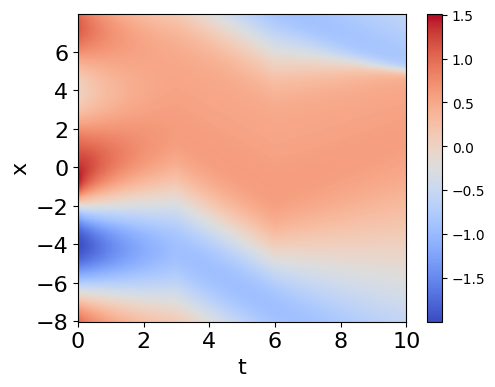} &  \includegraphics[height=4.5cm]{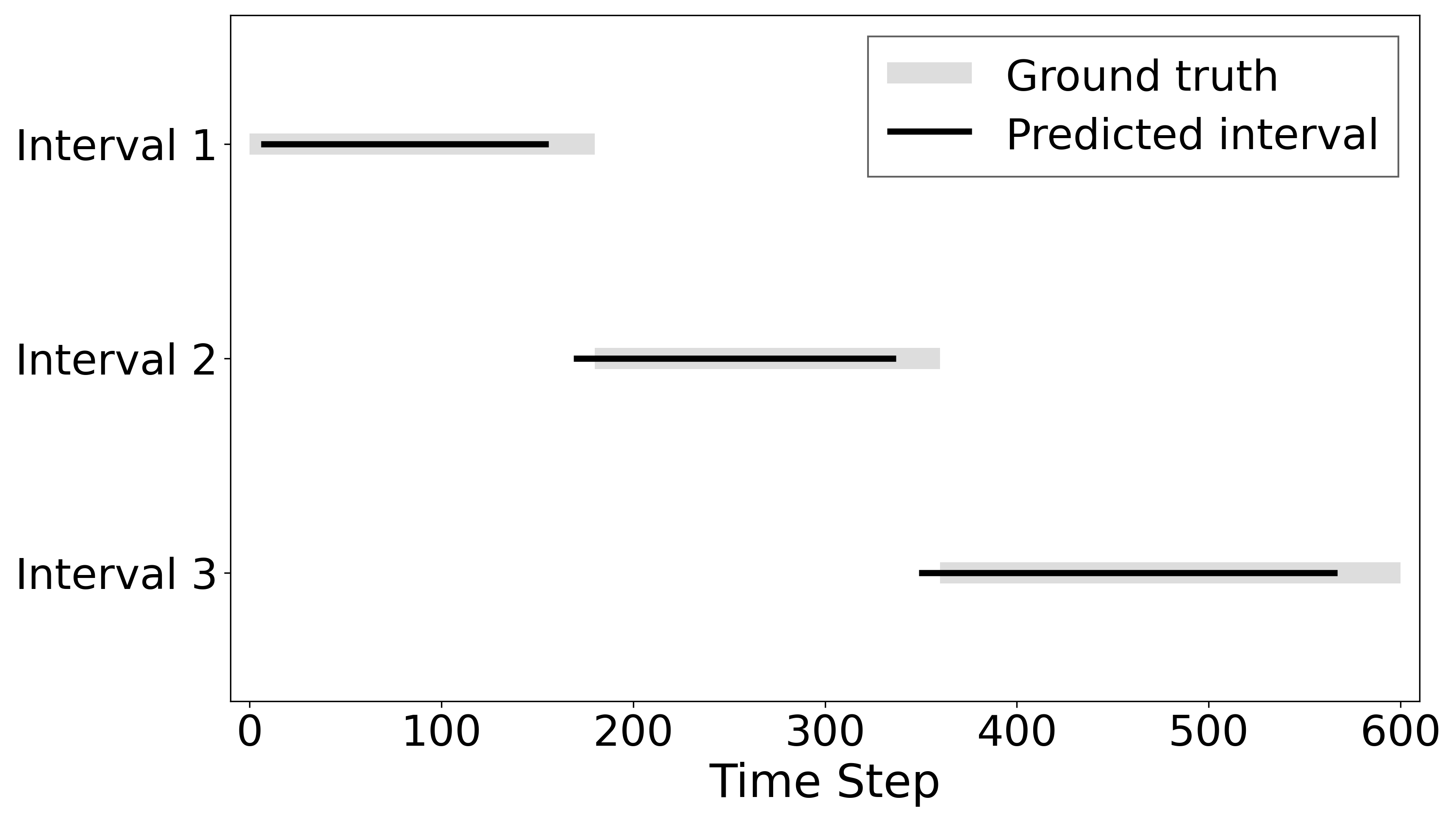}
\end{tabular}
{\small
\begin{tabular}{lll}
\toprule
True Int. & SLW Int.
& SLW result, to be compared to the true equation \eqref{eq:3changing}\\ \midrule
$[t^{0},t^{180})$ & $[t^{6},t^{156})$
& $u_t = 0.049996\,u + a_1(t)\,u_{xx} - 0.049995\,u^2$ \\
$[t^{180},t^{360})$ &  $[t^{169},t^{337})$
& $u_t = a_2(t)\,u_x$ \\
$[t^{360},t^{600})$ &  $[t^{349},t^{567})$
& $u_t = 0.09997\,u_{xx} + a_3(t)\,(u^2)_x$ \\ \bottomrule
\end{tabular}
}

\vspace{0.3cm}
\begin{tabular}{ccc}
(c) $a_1(t)$ vs. $(1+0.01 \cos t)$ & (d) $a_2(t)$ vs. $(1+0.01 \cos t) $ & (e) $a_3(t)$ vs. $0.5(-1+0.1 \cos t)$\\
\includegraphics[height=5.2cm]{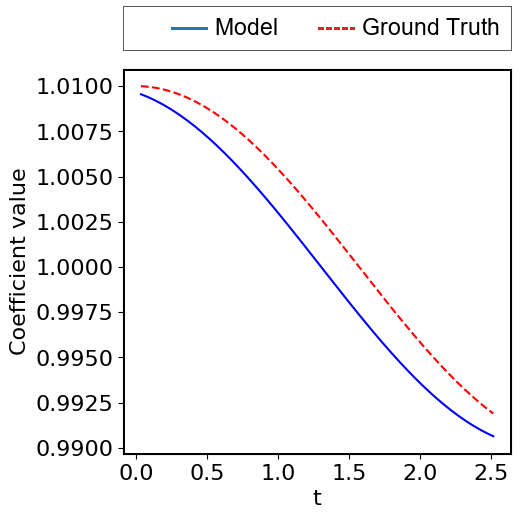} &
\includegraphics[height=5.2cm]{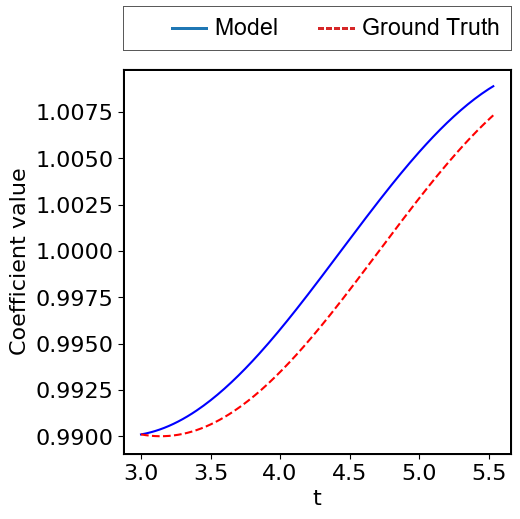} & 
\includegraphics[height=5.2cm]{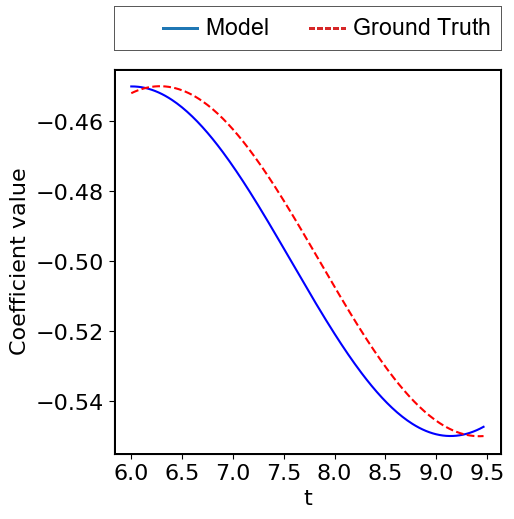}
\end{tabular}

{\small
\begin{tabular}{ccccccccc}
\toprule
Interval & Interval & Inclusion & Support & Support & $\mathrm{R (\%)}$ & $E_2$ & $E_\infty$ & $E_{\mathrm{res}}$ \\
TPR & PPV &  & TPR & PPV & & & & \\ \midrule
0.87 & 1.00 & 1.0 & 1.00 & 1.00 & 23.79 & $2.06\times10^{-3}$ & $2.10\times10^{-3}$ & $7.78\times10^{-5}$ \\
0.94 & 0.97 & 0.5 & 1.00 & 1.00 & 97.01 & $2.29\times10^{-3}$ & $2.29\times10^{-3}$ & $2.32\times10^{-4}$ \\
0.92 & 1.00 & 0.5 & 1.00 & 1.00 & 23.71 & $1.58\times10^{-2}$ & $1.64\times10^{-2}$ & $6.40\times10^{-4}$ \\
\bottomrule
\end{tabular}
}
\caption{PDE changing between three varying coefficient equations: Fisher--KPP, Transport and Burgers equations. (a) Given data \eqref{eq:3changing}. (b) Identified intervals  compared against the true interval. 
(c), (d) and (e) showing $a_1(t)$, $a_2(t)$ and $a_3(t)$  vs. ground truth \eqref{eq:3changing}.
The middle table shows the SLW-Ident result, and the bottom table shows the accuracy measures for Interval 1 to 3 from top to bottom.}
\label{fig:3changing}
\end{figure}

We generate the dataset $\mathcal{D}$ which is changing among three different PDEs with coefficients smoothly varying over three consecutive time intervals: a Fisher--KPP reaction--diffusion equation with cosine-modulated reaction strength, a transport equation with cosine-modulated advection speed, and a Burgers-type diffusion equation with cosine-modulated nonlinear advection:
\begin{align}
u_t &= 0.05\,u + \bigl(1 + 0.01\cos t\bigr)\,u_{xx} - 0.05\,u^2, &  t \in [t^0,t^{180}), \nonumber\\
u_t &= \bigl(1 + 0.01\cos t\bigr)\,u_x,  & t \in [t^{180},t^{360})\label{eq:3changing}\\
u_t &= 0.1\,u_{xx} + 0.5\bigl(-1 + 0.1\cos t\bigr)\,(u^2)_x &   t \in [t^{360},t^{600}) \nonumber
\end{align}
on \(x\in[-8,8]\), \(t\in[0,10)\) with \(N_x=N_t=601\) and
\[
  u(x,0)
  = \exp\!\bigl(-(x+1)^2\bigr)
    +\sin\!\bigl(\tfrac{\pi x}{8}\bigr)
    +\cos\!\bigl(\tfrac{\pi x}{4}\bigr).
\]
We present the result in Figure \ref{fig:3changing}. (a) shows the given data, (b) shows identified interval of one equation with the ground truth interval.  Table shows the SLW identified interval compared to the true interval.  
Since the underlying equations have varying coefficients, the SLW result also shows varying coefficient values, which we present as $a_1(t)$, $a_2(t)$, and $a_3(t)$ computed at each discrete value $t^n$ for each interval $I_1$, $I_2$ and $I_3$. These are plotted in (c) to (e) for each interval $I_{SLW} \cap I_{gt}$, i.e.,  $I_i$ and its paired true interval for $i=1,2$ and 3. The blue curves are $a_i(t)$ which are  compared to the ground truth in red.  The recovered coefficient accurately captures the trend of the true varying coefficient. The coefficient accuracy shows $E_2 \sim \mathcal{O}(10^{-2})$ to $\mathcal{O}(10^{-3})$. 
All three intervals of one equation are identified with perfect support recovery (TPR = PPV = 1.00) even if coefficient values are varying, and interval TPR is $\geq 0.8$ and interval PPV $\geq 0.9$.

\subsection{Noisy data of changing PDEs with time varying coefficients}
\label{subsec:3pv_noise15}

We consider the case with noisy observation with 15\% Gaussian noise (NSR = 0.15).  
We generate the given data $\mathcal{D}$ by PDE changing among three equations with smoothly varying coefficients: a transport equation with cosine-modulated advection speed, a Burgers-type diffusion with sine-modulated nonlinear advection, and a transport equation with cosine-modulated advection speed:
\begin{align}
u_t &= \bigl(-2 + 0.025\cos t\bigr)\,u_x, & t \in [t^0,t^{140}) \nonumber\\
u_t &= 0.01\,u_{xx} + 0.5\bigl(1 + 0.025\sin t\bigr)\,(u^2)_x,  & t \in [t^{140},t^{420}) \label{eq:varmix3_pdes}\\
u_t &= 0.5\bigl(1 + 0.025\cos t\bigr)\,u_x + 0.02\,u_{xx}, & t \in [t^{420},t^{700}) \nonumber
\end{align}
on \(x\in[-2\pi,2\pi]\), \(t\in[0,1)\) with \(N_x=N_t=701\) and initial condition 
\[
  u(x,0)
  = \sin(2\pi x)
    +0.5\,\sin\bigl(4\pi x + \tfrac\pi4\bigr)
    +0.25\,\cos\bigl(8\pi x - \tfrac\pi3\bigr)
    +e^{-5(x-0.3)^2}.
\]

\begin{figure}
\centering
\begin{tabular}{cc}
(a) & (b) \\
\includegraphics[height=4.5cm]{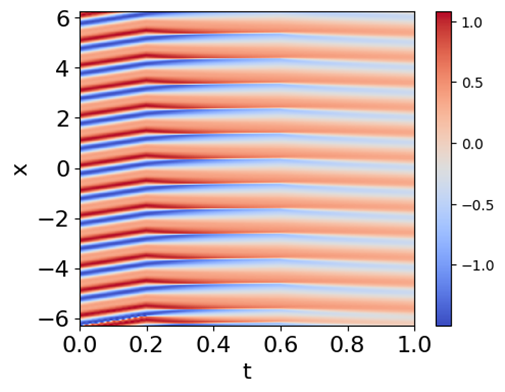}
&  \includegraphics[height=4.5cm]{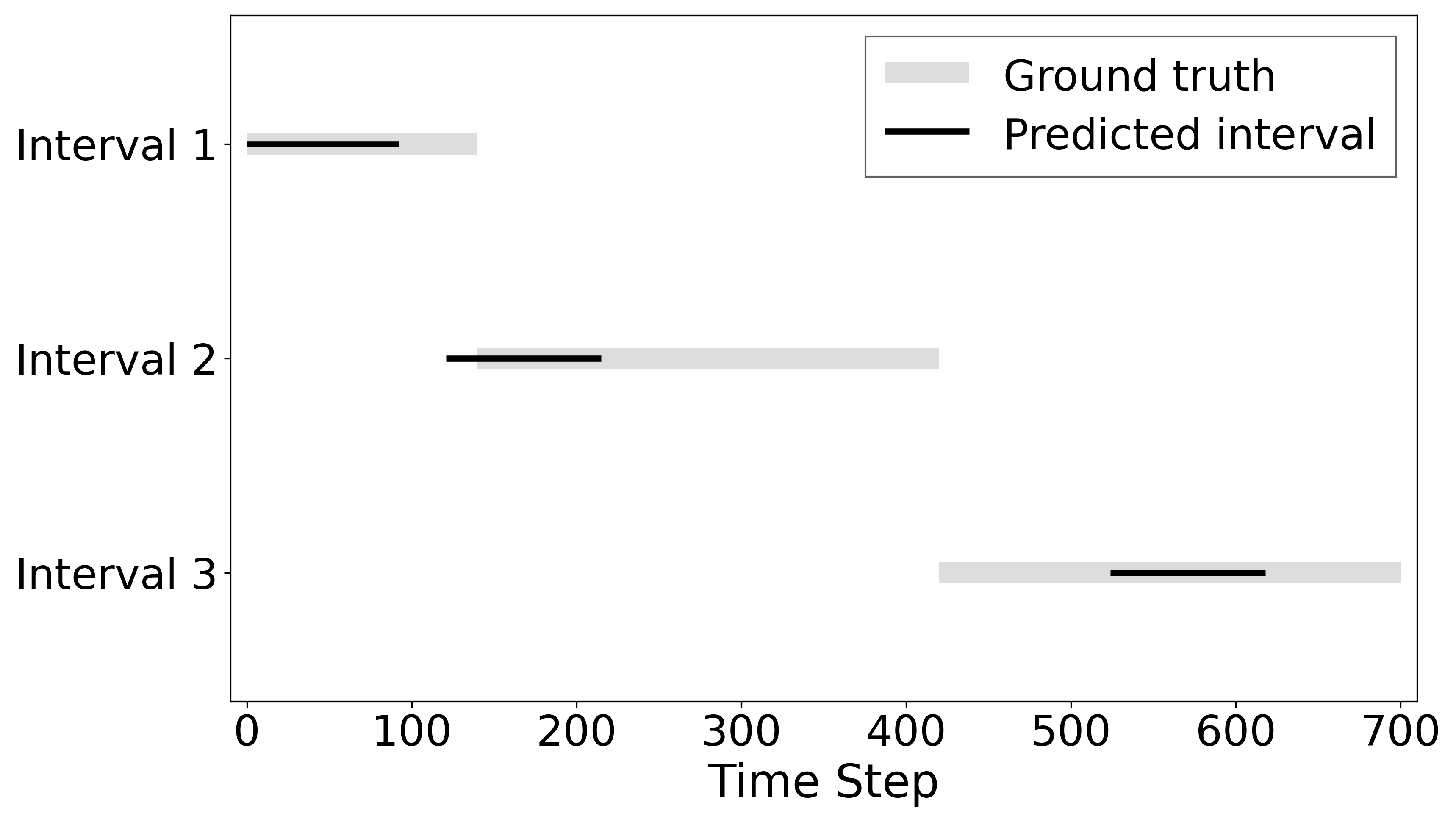}
\end{tabular}
{\small
\begin{tabular}{lll}
\toprule
True Int. & SLW Int.
& SLW result, to be compared to the true equation \eqref{eq:varmix3_pdes}\\ \midrule
$[t^{0},t^{140})$ & $[t^{0},t^{92})$
& $u_t = a_1(t)\,u_x$ \\
$[t^{140},t^{420})$ &  $[t^{121},t^{215})$
& $u_t = 0.01000\,u_{xx} + a_2(t)\,(u^2)_x$ \\
$[t^{420},t^{700})$ &  $[t^{524},t^{618})$
& $u_t = a_3(t)\,u_x + 0.01956\,u_{xx}$ \\ \bottomrule
\end{tabular}
}

\vspace{0.3cm}
\begin{tabular}{ccc}
(c) $a_1(t)$ vs. $(-2 +0.025\cos t)$ & (d) $a_2(t)$ vs. $0.5(1+0.025\sin t) $ & (e) $a_3(t)$ vs. $0.5(1+0.025\cos t)$\\
\includegraphics[height=5.5cm]{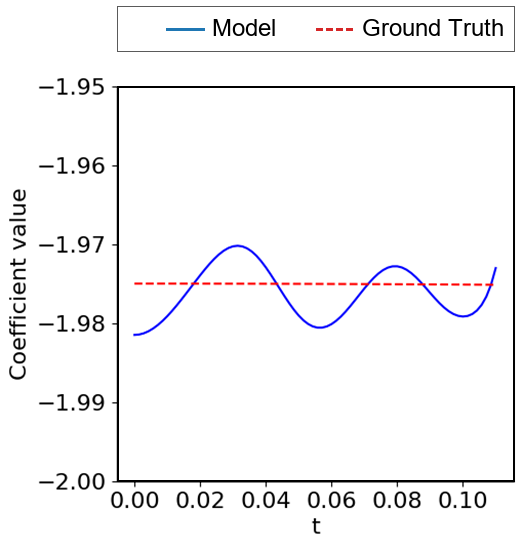} &
\includegraphics[height=5.5cm]{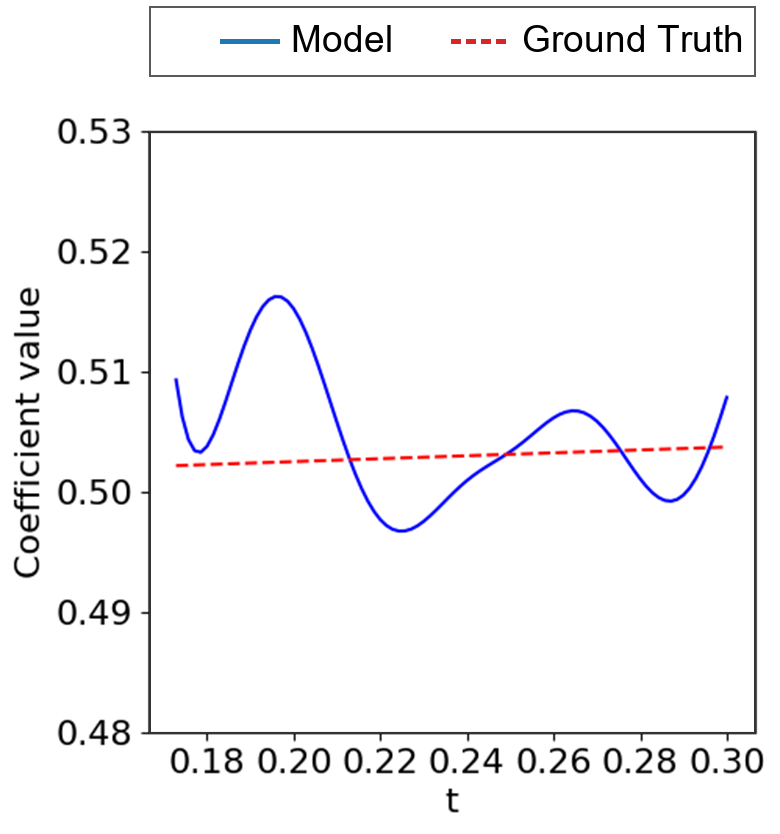}& 
\includegraphics[height=5.5cm]{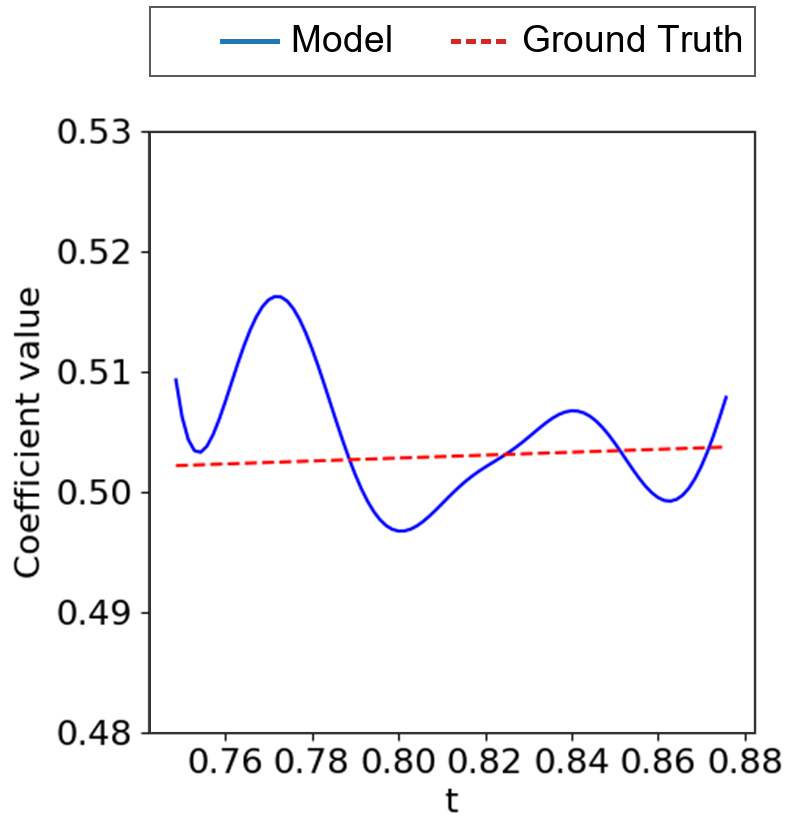}
\end{tabular}

{\small
\begin{tabular}{ccccccccc}
\toprule
Interval & Interval & Inclusion & Support & Support & $\mathrm{R (\%)}$ & $E_2$ & $E_\infty$ & $E_{\mathrm{res}}$ \\
TPR & PPV &  & TPR & PPV & & & & \\ \midrule
0.66 & 1.00 & 1.0 & 1.00 & 1.00 & 84.73 & $9.94\times10^{-3}$ & $9.94\times10^{-3}$ & $2.57\times10^{-2}$   \\
0.27 & 0.80 & 1.0 & 1.00 & 1.00 &42.36 & $9.07\times10^{-3}$ & $3.52\times10^{-2}$ & $1.34\times10^{-1}$   \\
0.34 & 1.00 & 1.0 & 1.00 & 1.00 & 59.49 & $1.86\times10^{-2}$ & $8.33\times10^{-2}$ & $2.27\times10^{-1}$  \\
 \bottomrule
\end{tabular}
}
\caption{PDE changing between three varying coefficient equations with 15\% Gaussian noise (NSR = 0.15).  (a) Given data \eqref{eq:varmix3_pdes}.  (b) Identified intervals  compared against the true interval.   (c), (d) and (e) showing $a_1(t)$, $a_2(t)$ and $a_3(t)$  vs. ground truth \eqref{eq:varmix3_pdes}.  The middle table shows the SLW-Ident result, and the bottom table shows the accuracy measures for Interval 1 to 3 from top to bottom.}
\label{fig:3varying_Noise1}
\end{figure}

Figure~\ref{fig:3varying_Noise1} (a) shows the given data with 15 \% Gaussian noise (NSR = 0.15), and Figure~\ref{fig:3varying_Noise1} (b) overlays the identified interval of one equation  compared against the true interval. 
We present the varying coefficients as  $a_1(t)$, $a_2(t)$, and $a_3(t)$ which are computed at each discrete value $t^n$ for each interval $I_1$, $I_2$ and $I_3$. These are plotted in (c) to (e) for each interval $I_{SLW} \cap I_{gt}$, i.e.,  $I_i$ and its paired true interval for $i=1,2$ and 3. The blue curves are $a_i(t)$ which are  compared to the ground truth in red.  Differences are less than $0.02$ throughout each interval. 
SLW-Ident recovers each interval's support perfectly, although identified intervals are shorter as indicated by low Interval TPR values, the identified intervals lie in the true intervals, with large interval PPV. 

We repeated this experiment 50 times.  In some cases, the last interval $\{u_x, u_{xx}\}$ was found separated into two intervals due to noise: the model identifies the support combination  $\{u_x\}$, $\{u_{xx},(u^2)_x\}$, $\{u_x,u_{xx}\}$ in three-piece form or  $\{u_x\}$, $\{u_{xx},(u^2)_x\}$, $\{u_x,u_{xx}\}$, $\{u_x,u_{xx}\}$  in four-piece form.
Table~\ref{tab:varmix3_support_freq} summarizes the identified support patterns, showing that in 72\% of the trials the proposed model finds the correct support.
\begin{table}[h]
  \centering
  \begin{tabular}{lr}
\toprule
\textbf{Support Combination} & \textbf{Number of Identification (out of 50) } \\
\midrule
\(\{u_x\},\;\{u_{xx},(u^2)_x\},\;\{u_x,u_{xx}\}\) & 27 + 9\\
\(\{u_x\},\;\{(u^2)_x,(u^3)_x\},\;\{u_x,u_{xx}\}\)  & 10 + 1\\
\(\{u_x\},\;\{u_x,u_{xx}\}\)  &  2 + 1 \\
\bottomrule
\end{tabular}
\caption{PDE changing between three varying coefficient equations with 15\% Gaussian noise (NSR = 0.15), 50 independent experiments.
First column shows intervals of one equation and the support identified in each interval.  In the second column, we present the number of times three intervals are identified, and the second number indicated with an addition is when the last interval $\{u_x, u_{xx}\}$ is separated to two intervals with the same support.} 
\label{tab:varmix3_support_freq}
\end{table}
Detailed statistics on the identified interval boundaries and coefficient errors for these successful trials are reported in Appendix~\ref{append: noise15_stats}, which shows that the recovered intervals and coefficients remain accurate.  An additional noisy experiment with three constant-coefficient PDEs at 10\% NSR is presented in Appendix~\ref{append:additional_experiments}.

\subsection{Comparison with global varying-coefficient methods under noise}
\label{subsec:global_vs_local}

In this section, we present a comparison with another method identifying PDE with varying coefficients.  Methods such as IDENT \cite{kang2021ident} and GP-IDENT \cite{he2023group} present equations with varying coefficients using finite element basis for the entire domain.  Using Group subspace pursuit, GP-IDENT finds all the terms which were activated during any part of the domain.   For accurate identification, a large number of basis functions, i.e., fine grid, is beneficial, which induces a large linear system.   The proposed method is based on local identification which reduces problem dimensionality and only solves a small size linear system.  This local WeakIdent is further stabilized by sampling, considering residual error, and grouping of transition locations to find region of one equation. 

\begin{figure}
  \centering
\includegraphics[width=0.7\textwidth]{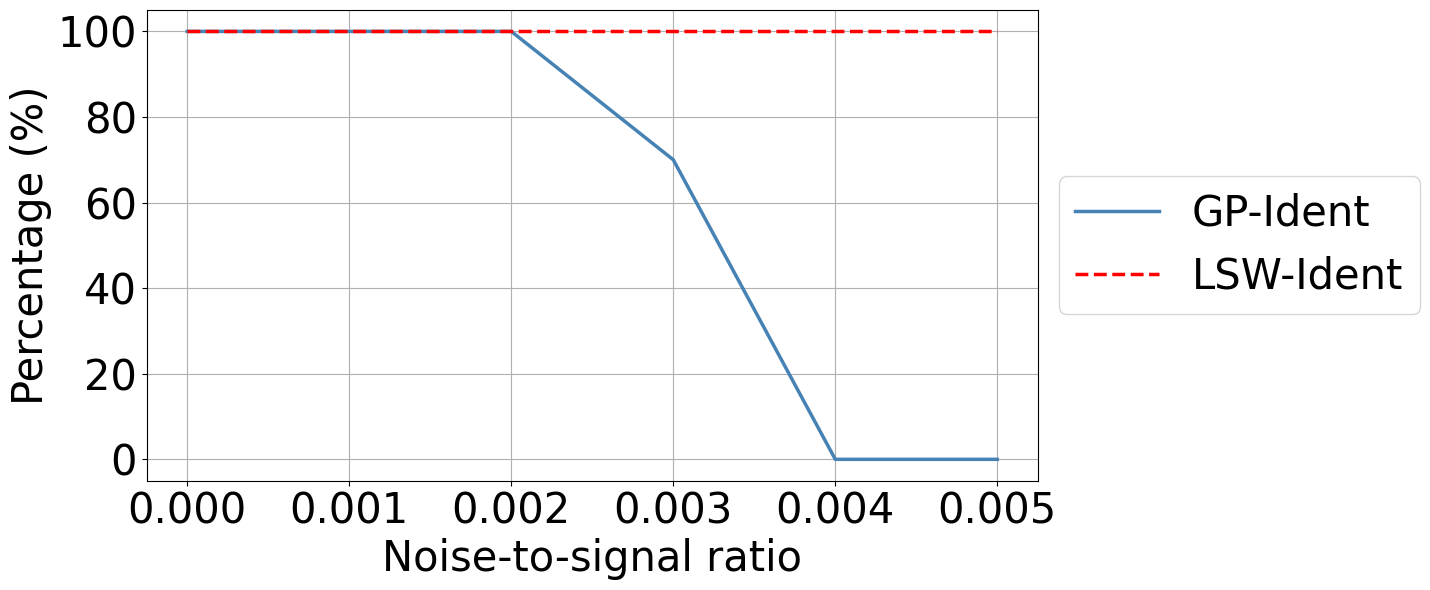}
  \caption{
  The two-piece viscous Burgers equation \eqref{eq:2burgers}.
  Support recovery success rate vs.\ noise-to-signal ratio (NSR = 0, 0.001, 0.002, 0.003, 0.004, and 0.005).  SLW-Ident result in blue and GP-IDENT \cite{he2023group} in red.}
  \label{fig:comparison_noise}
\end{figure}

Due to the difference in the set-up, comparison is not straightforward;
we restrict the comparison to where both intervals share the same support but differ in their coefficient values.  We show that the local approach can make the support recovery of SLW-Ident more robust to noise than the global design.
We compare SLW-Ident with a global finite-element expansion method, GP-IDENT \cite{he2023group} on a two-piece viscous Burgers problem with a sudden coefficient change:
\begin{align}
u_t & =  0.8\,u_{xx} + 2(u^2)_x, &  0<t\le0.01 \label{eq:2burgers}\\
u_t & = 1.6\,u_{xx}+ 4(u^2)_x, & 0.01<t\le0.02  \nonumber
\end{align}
on $x\in[-2,2]$, $t\in[0,0.02]$, and $N_x=N_t=512$. 
The initial condition is
\[
\begin{aligned}
u(x,0) &= \sin\!\bigl(\pi(2x-0.1)\bigr)
        +\cos\!\bigl(\pi(5x-0.2)\bigr)
        +\cos\!\bigl(\pi(3x-0.3)\bigr)\cos\!\bigl(\pi(x+0.1)\bigr)\\
       &\quad{}+\sin\!\bigl(\pi(4x+0.5)\bigr)
        +0.8\exp\!\Bigl(-\tfrac{(x-0.2)^2}{2(0.05)^2}\Bigr)
        +0.6\exp\!\Bigl(-\tfrac{(x-0.5)^2}{2(0.08)^2}\Bigr)
        +5,
        \quad x<1,\quad\text{and }0\text{ otherwise}.
\end{aligned}
\]
We add Gaussian perturbations of noise-to-signal ratios, NSR = 0, 0.001, 0.002, 0.003, 0.004, and 0.005, and perform 10 independent trials at each noise level.  To count for successful identification, for SLW-Ident, the correct support \(\{u_{xx},\,(u^2)_x\}\) should be identified in each of its identified intervals, and for GP-IDENT, the identification is considered successful as long as it identifies a single correct global support. 
Figure~\ref{fig:comparison_noise} presents success rates among 10 independent trials for each NSR level. SLW-Ident maintains a perfect support recovery rate across all six noise levels, while GP-IDENT's support accuracy reduces
for $\text{NSR}\geq 0.003$.  

\subsection{PDE with space--time varying coefficients}
\label{subsec:localization_space_time}

While the explanation and experiments are presented for PDEs with time-varying coefficients, the proposed method SLW-Ident can be extended to PDEs with space- and time-varying coefficients. 
We consider a domain partitioned into four quadrants, each governed by a different constant-coefficient PDE:
\begin{align} 
&u_t = 0.02\,u_{xx}, &&(x,t)\in[x^0,x^{300})\times[t^0,t^{280}), \nonumber\\
&u_t = 2.0u_x, &&(x,t)\in[x^{300},x^{700})\times[t^0,t^{280}),\label{eq:spaceTime}\\
&u_t = 1.0u_x, &&(x,t)\in[x^0,x^{300})\times[t^{280},t^{700}),\nonumber\\
&u_t = -1.0u_x + 0.03u_{xx}, &&(x,t)\in[x^{300},x^{700})\times[t^{280},t^{700}),\nonumber
\end{align}
on \(x\in[0,1]\), \(t\in[0,0.3]\) with \(N_x=N_t=701\) and initial condition 
\[
  u(x,0)
  = 10\left(\sin(2\pi x) + 5\sin(4\pi x + \frac{\pi}{4}) + 7\cos(8\pi x -\frac{\pi}{3}) + \exp\left(-5(x-0.3)^2\right)\right).
\]

\begin{figure}[t]
\centering
\begin{tabular}{ccc}
(a) & (b) & (c) \\
\includegraphics[width=0.32\textwidth]{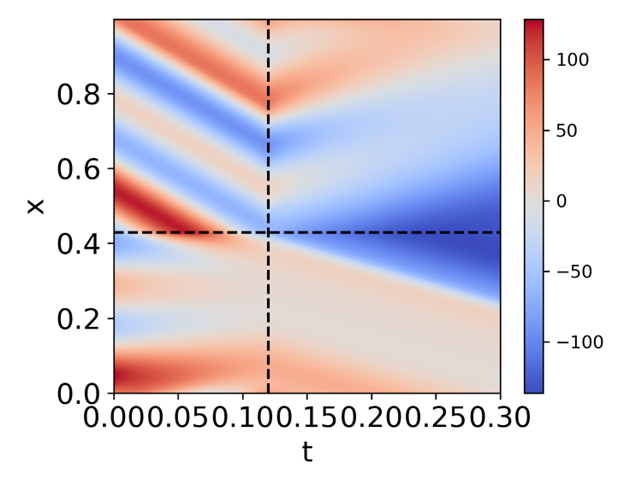} &
\includegraphics[width=0.32\textwidth]{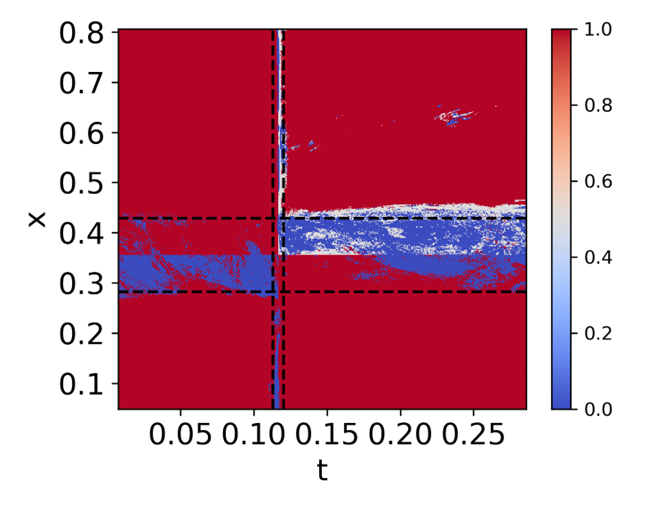} &
\includegraphics[width=0.32\textwidth]{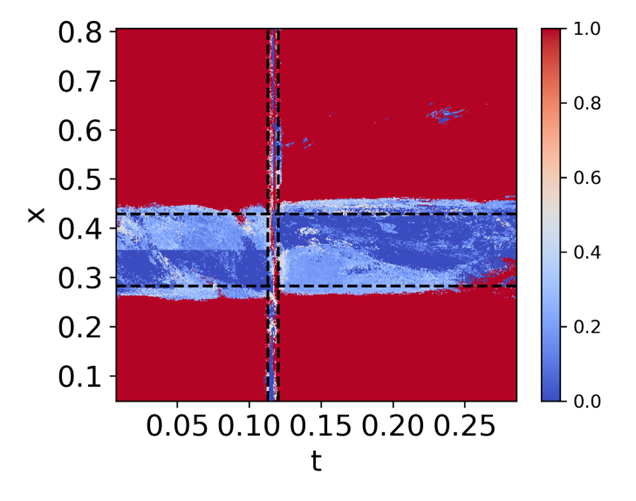}\\
(d) & (e) & (f) \\
\includegraphics[width=0.32\textwidth]{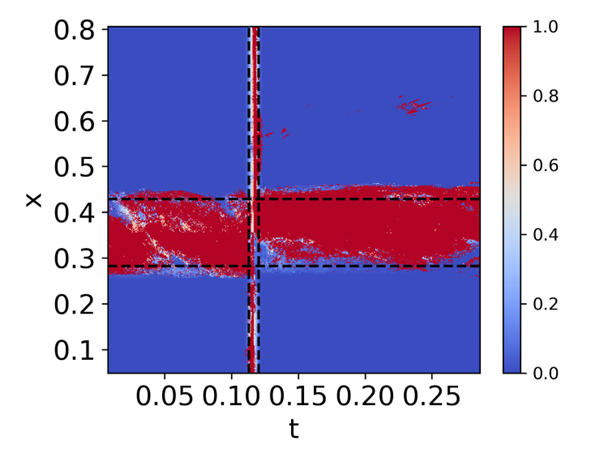} &
\includegraphics[width=0.32\textwidth]{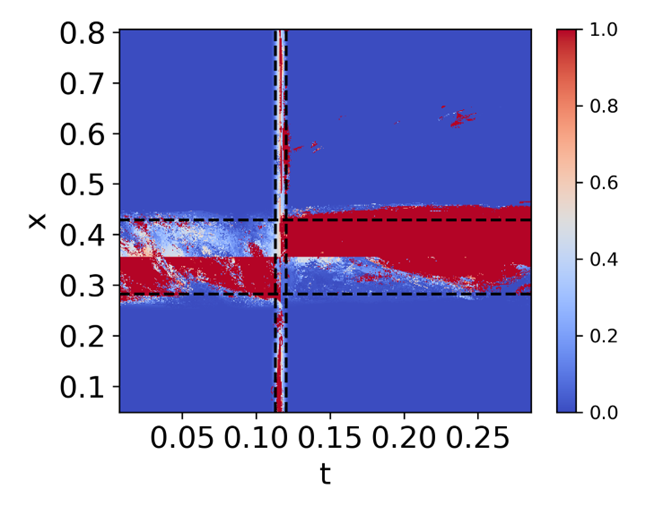} &
\includegraphics[width=0.32\textwidth]{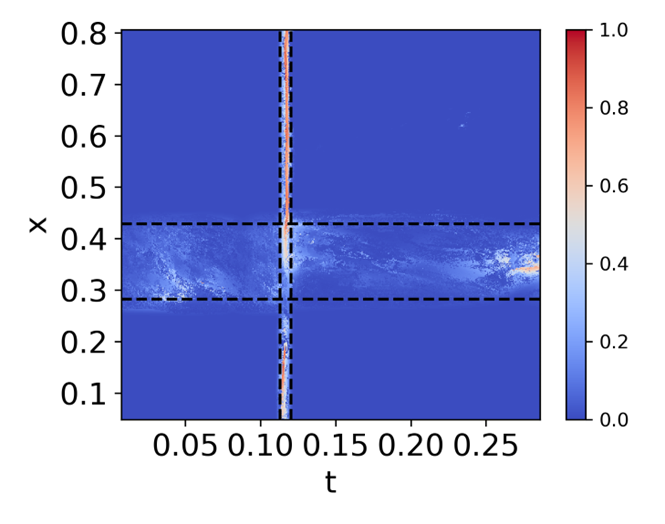}
\end{tabular}  
\caption{PDE with space and time varying coefficients \eqref{eq:spaceTime}.  (a) The given data $\mathcal{D}$; dashed lines present the true equation-change boundaries at \(x^{300}\) and \(t^{280}\). (b) Support TPR and (c) Support PPV show mostly values of 1, and coefficient errors (d) \( E_2\), (e) \( E_\infty\), and (f) \(E_{\mathrm{res}}\) show low values except at the transitions. Accuracy measures are evaluated at each patch-starting location. }
\label{fig:localization_space_time}
\end{figure}

Figure~\ref{fig:localization_space_time} shows (a) the given data $\mathcal{D}$ given by \eqref{eq:spaceTime} where dashed lines mark the true equation-change boundaries at \(x^{300}\) and \(t^{280}\). (b)--(f) present accuracy measures evaluated at each patch-start location.  (b) Support-TPR and (c) Support-PPV show mostly values of 1, and coefficient errors (d) \(E_2\), (e) \(E_\infty\), and (f) \( E_{\mathrm{res}}\) show low values except at the transition band given by the patch size.  

\section{Conclusion}
\label{s: conclusion}

We propose SLW-Ident to identify changing PDE from a single set of given data.  The main idea is to utilize the robustness of WeakIdent in a local domain.  The patch size for Local WeakIdent \eqref{eq:localWeakIdent} is about 3 times larger than the test function support, and it gives stable identification.  Local WeakIdent \eqref{eq:localWeakIdent} is a simplified version compared to the full WeakIdent \cite{tang2023weakident}.   Since the simplified Local WeakIdent only uses a small amount of local data, computation is efficient and we do not use high dynamic region regression. 
Due to noise and reduced domain (compared to using the full data), there may be some perturbations in the local identification, thus, we propose SLW-Ident in Section \ref{sec: SLW-Ident} to further stabilize the identification process.  We first sample patches in the domain, find equations for each sampled patch, then consider residual errors \eqref{eqn:residual_def} of many candidate equations to find the transitions. Within each region of one equation, we pick the most frequently identified equation as the identified equation and use regression \eqref{eqn: time-varying-regres} to find varying coefficients.  We provide confidence bounds for dominant-support selection to guide the number of patches.  
We demonstrate the performance of the proposed SLW-Ident by various experiments such as mixtures of constant and varying-coefficient PDEs, and datasets with up to 15\% Gaussian noise.  The local approach gives efficiency and flexibility to the identification of PDEs.  

\bibliographystyle{plain}
\bibliography{ref_localW}

\appendix

\section{Proofs of confidence guarantees under general and Monte-Carlo sampling schemes}\label{Asec:proofs}

In this section, we provide proofs of the results presented in Section~\ref{subsec: max_dev} and \ref{subsec: exp_dev}.

\subsection{Proof of Theorem \ref{thm:max-dev}}\label{app:proof-thm-entropy-concentration}

\begin{proof}

Let $\mathcal{B}_1, \dots, \mathcal{B}_{N_\mathrm{sample}}$ be random variables representing the sampled local patches from the domain. Since the patches are drawn independently from a fixed distribution, the sequence $\{\mathcal{B}_i\}_{i=1}^{N_\mathrm{sample}}$ is i.i.d. 
Let $\mathbf{s}_i$ denote the support identified from patch $\mathcal{B}_i$,
then $\{\mathbf{s}_i\}_{i=1}^{N_\mathrm{sample}}$ are also i.i.d.\ random variables, since local WeakIdent is deterministic, and
$\mathbf{s}_i$ takes values in $\mathcal{S}_\mathrm{cand}$ with distribution $\mathbb{P}(\mathbf{s}_i = \mathbf{s}_k) = q_k$. 
For each patch $i\in \mathcal{I}$, let $\tilde p_i$ denote the empirical probability associated with the support of patch $i$, i.e., $\tilde p_i = p_k$ if $\tilde {s}_i=\mathbf{s}_k$. Let $N_k$ denote the number of times support $\mathbf{s}_k$ is identified among the sampled local patches. Then $p_k = \frac{N_k}{N_\mathrm{sample}}$. 
We let $X_i := -\log \tilde p_i$ and $0 \le X_i \le -\log p_{\min}$ since $0<\tilde p_i \le 1$.  
Since $p_{\min} \ge \frac{1}{N_\mathrm{sample}}$, it follows that
$ X_i \le \log N_\mathrm{sample}$. 

\textbf{Step 1: Representation of empirical entropy:}
\begin{align*}
H(P)
&= -\sum_{k=1}^{N_\mathrm{cand}} p_k \log p_k = -\sum_{k=1}^{N_\mathrm{cand}} \frac{N_k}{N_\mathrm{sample}}\log \frac{N_k}{N_\mathrm{sample}} \\
&= -\frac{1}{N_\mathrm{sample}}\sum_{k=1}^{N_\mathrm{cand}} N_k \log p_k = -\frac{1}{N_\mathrm{sample}}\sum_{i=1}^{N_\mathrm{sample}} \log \tilde p_i = \frac{1}{N_\mathrm{sample}}\sum_{i=1}^{N_\mathrm{sample}} X_i.
\end{align*}

\textbf{Step 2: Limiting behavior:}
Let $\tilde q_i$ denote the true probability of identifying the dominating support of patch $\mathcal{B}_i$ among all sampled local patches. Then $\tilde q_i = q_K$ if $\tilde s_i = s_\mathrm{K}$, where $q_K$ is the true probability of identifying support $s_\mathrm{K}$.  
Define the ideal random variables
$
Y_i := -\log \tilde q_{i},
$
then $\{Y_i\}$ are i.i.d.\ and 
$
0 \le Y_i \le - \log q_{\min}, 
$
and
$
\mathbb{E}[Y_i] = -\sum_{K=1}^{N_\mathrm{cand}^{true}} q_j \log q_j = H(Q).
$
By the Law of Large Numbers,
\[
\frac{1}{N_\mathrm{sample}}\sum_{i=1}^{N_\mathrm{sample}} Y_i \to H(Q) \quad \text{as } N_\mathrm{sample}\to\infty.
\]

\textbf{Step 3: Approximation argument:} Since $p_k \to q_k$ almost surely as $N_\mathrm{sample}\to\infty$, we have
$
\tilde p_i \to \tilde q_i
$
that
$
X_i = -\log \tilde p_i \;\to\; Y_i = -\log \tilde q_i $ almost surely.
Moreover,
\[
\left|
\frac1{N_\mathrm{sample}}\sum_{i=1}^{N_\mathrm{sample}}X_i-\frac1{N_\mathrm{sample}}\sum_{i=1}^{N_\mathrm{sample}}Y_i
\right|
\le
\frac1{N_\mathrm{sample}}\sum_{i=1}^{N_\mathrm{sample}}|X_i-Y_i|
\le
\max_{k}\bigl|\log p_k-\log q_k\bigr|.
\]
Since $p_k\to q_k$ almost surely for each $k$ and there are finitely many $k\in N_\mathrm{cand}$, it follows that
\[
\max_{k}\bigl|\log p_k-\log q_k\bigr|\to 0
\qquad\text{almost surely}
\]
and
\[
\left|
\frac1{N_\mathrm{sample}}\sum_{i=1}^{N_\mathrm{sample}}X_i-\frac1{N_\mathrm{sample}}\sum_{i=1}^{N_\mathrm{sample}}Y_i
\right|\to 0
\qquad\text{almost surely}.
\]

Since almost sure convergence implies convergence in probability, it follows that
\[
\left|
\frac1{N_\mathrm{sample}}\sum_{i=1}^{N_\mathrm{sample}}X_i-\frac1{N_\mathrm{sample}}\sum_{i=1}^{N_\mathrm{sample}}Y_i
\right|
\xrightarrow{\mathbb P}0.
\]
Hence, for every $\eta>0$, there exists $N_{\mathrm{sample}}(\eta)$ such that for all $N_\mathrm{sample}\ge N_{\eta}$,
\[
\mathbb P\!\left(
\left|
\frac1{N_\mathrm{sample}}\sum_{i=1}^{N_\mathrm{sample}}X_i-\frac1{N_\mathrm{sample}}\sum_{i=1}^{N_\mathrm{sample}}Y_i
\right|>\eta
\right)
<\eta.
\]
\medskip
\noindent\textbf{Step 4: Application of Hoeffding's inequality.}

Since $\{Y_i\}$ are i.i.d.\ and bounded in $[0,\log |N_\mathrm{cand}|]$, Hoeffding's inequality yields
\[
\mathbb{P}\!\left(
\left|\frac{1}{N_\mathrm{sample}}\sum_{i=1}^{N_\mathrm{sample}} Y_i - H(Q)\right|
\ge \epsilon
\right)
\le
2\exp\!\left(-\frac{2N_\mathrm{sample}\epsilon^2}{(\log q_{\min})^2}\right).
\]

Combining this with the approximation in Step 3 above, by using triangle inequality, we obtain that for every $\eta >0$, there exists some $N_{\eta} >0$ such that for every large $|I| \ge N_{\eta}$,
\[
\mathbb{P}\!\bigl(|H(P)-H(Q)|\ge\epsilon + \eta \bigr)
\;\le\;
2\exp\!\Bigl(-\tfrac{2\,N_\mathrm{sample}\,\epsilon^2}{(\log q_{\min})^2}\Bigr) + \eta.
\]
This completes the proof.
\end{proof}

\subsection{Proof of Lemma \ref{lem:kl-entropy-bound}}\label{app:proof-lem-kl-entropy-bound}

\begin{proof}
By definition,
\[
D_{\mathrm{KL}}(P\parallel Q)
   =\sum_{k=1}^{N_\mathrm{cand}} p_k\log\frac{p_k}{q_k}
  =-H(P)-\sum_{k=1}^{N_\mathrm{cand}} p_k\log q_k
  \le -H(P)-\sum_{k=1}^{N_\mathrm{cand}} p_k\log q_{\min}
   =-H(P)-\log q_{\min}.
\]
Since \(-\log q_{\text{min}} =\log\frac1{q_{\min}}\le\sum_k\log\frac1{q_k}\le\sum_k\frac{q_k}{q_{\min}}\log\frac1{q_k}\), it follows
\[
  D_{\mathrm{KL}}(P\parallel Q)
  \le -H(P)+\sum_{k=1}^{N_\mathrm{cand}}\frac{q_k}{q_{\min}}\log\frac1{q_k}
  =-H(P)+\frac1{q_{\min}}H(Q)
  =\frac{H(Q)-H(P)}{q_{\min}}
   +\frac{1-q_{\min}}{q_{\min}}H(P).
\]
\end{proof}

\subsection{Proof of Theorem \ref{thm:max-dev-kl}}\label{app:proof-thm-max-dev-kl}
\begin{proof}
By Pinsker's inequality, 
\[
  d_{\mathrm{TV}}(P,Q)^2
  \;\le\;
  \tfrac12D_{\mathrm{KL}}(P\parallel Q)
  \;\le\;
  \frac{H(Q)-H(P)}{2\,q_{\min}}
  +\frac{1-q_{\min}}{2\,q_{\min}}\,H(P),
\]
using Lemma~\ref{thm:max-dev} where the total variation \(d_{\mathrm{TV}}(P,Q)=\tfrac12\sum_k|p_k-q_k|\).  On the other hand, for any index 
  \(k_0 \in \{1, 2, \dots, N_\mathrm{cand} \}\), 
\begin{align*}
  d_{\mathrm{TV}}(P,Q)
  & = \frac12\sum_k|p_k-q_k|
  \;=\;
  \frac12\Bigl(|p_{k_0}-q_{k_0}|+\sum_{k\neq k_0}|p_k-q_k|\Bigr)
  \\ & \ge\; \frac12\Bigl(|p_{k_0}-q_{k_0}|+|\sum_{k\neq k_0}p_k-\sum_{k\neq k_0}q_k|\Bigr) = \;|p_{k_0} - q_{k_0}|.
    \end{align*}
Hence
\[
  \max_{1 \le k \le N_\mathrm{cand} }|p_k-q_k|
  \le d_{\mathrm{TV}}(P,Q)
  \le \sqrt{\frac{H(Q)-H(P)}{2\,q_{\min}}
            +\frac{1-q_{\min}}{2\,q_{\min}}\,H(P)}.
\]
Substituting \(|H(Q)-H(P)|\le\epsilon + \eta\) and invoking Theorem~\ref{thm:max-dev} completes the proof. 
\end{proof}

\subsection{Proof of Corollary \ref{cor:argmax}}\label{app:proof-cor-argmax}
\begin{proof}
We want \(|p_k - q_k|\le R-0.5\) for all \(k\), which guarantees \(\arg\max p_k=\arg\max q_k\) because if \(p_k = p_{\text{max}}\) then \(|p_k-q_k|\;\leq\;R-0.5\;=\;p_k - 0.5\), which implies that $p_k - q_k \;\leq\; p_k-0.5$ so $q_k \;\ge\;0.5$ and $q_k = q_{\max}$. From Theorem~\ref{thm:max-dev-kl},
\[
  \max_k|p_k-q_k|
  \;\le\;
  \sqrt{\frac{\epsilon + \eta}{2q_{\min}}
    +\frac{1-q_{\min}}{2q_{\min}}\,H(P)}
\]
with probability at least \(1-2\exp\bigl(-\frac{2N_\mathrm{sample}\epsilon^2}{(\log q_{\min})^2}\bigr) - \eta\).  Setting
\[
  \sqrt{\frac{\epsilon + \eta}{2q_{\min}}
    +\frac{1-q_{\min}}{2q_{\min}}\,H(P)}
  = R-0.5
  \quad\text{by setting}\quad
  \epsilon
  = 2(R-0.5)^2\,q_{\min}
    - (1-q_{\min})\,H(P) - \eta,
\]
yields
\[
  \mathbb{P}\bigl(\max_k|p_k-q_k|\le R-0.5\bigr)
  \;\ge\;
  1 - 2\exp\!\Bigl(-\tfrac{2N_\mathrm{sample}}{(\log q_{\min})^2}\bigl(2(R-0.5)^2q_{\min}-(1-q_{\min})H(P)\bigr)^2 - \eta\Bigr) - \eta,
\]
which is exactly the stated bound.
\end{proof}

\subsection{Proof of Theorem~\ref{thm:variance-clt}}\label{app:proof-thm-variance-clt}
\begin{proof}
Let \(B_t\) and \(B_x\) be the patch sizes in the \(t\) and \(x\) directions, respectively. Define the map \(g_k: \mathcal{Z} \to\{0,1\}\) by
\[
  g_K(\vec z)
  = \begin{cases}
    1, & \text{if the patch starting at z identifies support }s_\mathrm{K},\\
    0, & \text{otherwise},
  \end{cases}
\]
where $\mathcal{Z}$ is the set of all valid starting location of patches, and each \(\vec z=(x,t)\) defines a patch \([x, x + B_x] \times [t, t + B_t) \).  Since patch start locations are sampled uniformly at random, the Monte Carlo estimate of \(\mathbb{E}[g_K]\) is
\[
  F_{N_\mathrm{sample}}[g_K]
  = \frac{1}{N_\mathrm{sample}}\sum_{i=1}^{N_\mathrm{sample}} g_K(\vec z_i),
\]
while the true expectation is
\[
  F[g_K]
  = \mathbb{E}_{\vec z}[g_K]
  = \int g_K(\vec z)\,\mathbb{P}(\vec z),
  \quad \vec z\sim\mathcal U\bigl([0,X-B_x]\times[0,T-B_t]\bigr).
\]
By the Central Limit Theorem, the Monte Carlo error 
\(\epsilon_{N_\mathrm{sample}}(g_j)=F_{N_\mathrm{sample}}[g_K]-F[g_K]\) satisfies 
\[
  \mathbb{E}[\epsilon_{N_\mathrm{sample}}^2]
  = \frac{\mathrm{Var}(g_K)}{N_\mathrm{sample}},
  \quad
  \epsilon_{N_\mathrm{sample}}(g_K)
  \xrightarrow{d}
  \mathcal N\!\Bigl(0,\tfrac{\mathrm{Var}(g_K)}{N_\mathrm{sample}}\Bigr).
\]
Noting that
\[
  F_{N_\mathrm{sample}}[g_K]
  = \frac{\#\{i: \text{patch }i\text{ yields } s_\mathrm{K}\}}{N_\mathrm{sample}}
  = p_K,
  \quad
  F[g_K]=q_K.
\]
On the other hand, 
\[
  \mathrm{Var}(g_K)
  = \int\bigl(g_K(\vec z)-q_K\bigr)^2\mathbb{P}(\vec z)\,dz
  = q_K(1-q_K).
\]
Thus, we obtain
\[
  \mathbb{E}[(p_K-q_K)^2]
  = \frac{q_K(1-q_K)}{N_\mathrm{sample}},
  \quad
  p_K - q_K
  \xrightarrow{d}
  \mathcal N\!\Bigl(0,\tfrac{q_K(1-q_K)}{N_\mathrm{sample}}\Bigr).
\]
\end{proof}

\subsection{Proof of Proposition~\ref{prop:monte-carlo}}\label{app:proof-cor-monte-carlo}
\begin{proof}[Justification]
Let \(k^\star=\arg\max_{k: \ \mathbf{s}_k\in \mathcal{S}_\mathrm{cand}}p_k\), so that \(p_{k^\star}=R>0.5\). 
Under the Gaussian approximation in Proposition~\ref{prop:monte-carlo}, the deviation of the dominant empirical support satisfies
\[
p'_{k^\star}-q_{k^\star}
\sim
\mathcal{N}\!\left(0,\frac{R(1-R)}{N_\mathrm{sample}-1}\right).
\]
The event
\[
|p'_{k^\star}-q_{k^\star}|<R-0.5
\]
ensures that the dominant support remains separated from the threshold \(0.5\) under the Gaussian perturbation. 

Now we use the approximation \(\max_k p'_k \approx R\). Under the condition 
\(\arg\max_k p_k=\arg\max_k p'_k\), this implies \(p'_{k^\star}\approx R\). 
Therefore, on the event above,
\[
q_{k^\star} \geq p'_{k^\star}-(R-0.5)\approx 0.5.
\]
Since a support whose probability exceeds \(0.5\) is necessarily the unique dominant support, this event provides a sufficient condition, under the stated approximation, for the empirical dominant support to agree with the population dominant support. Therefore,
\[
\mathbb{P}\bigl(\arg\max_k p'_k=\arg\max_k q_k \mid \arg \max{p_k} = \arg \max{p'_k}\bigr)
\;\gtrsim\;
\mathbb{P}\bigl(|p'_{k^\star}-q_{k^\star}|<R-0.5\bigr).
\]
Using the Gaussian density gives
\[
\mathbb{P}\bigl(|p'_{k^\star}-q_{k^\star}|<R-0.5\bigr)
=
\sqrt{\frac{|I|-1}{2\pi R(1-R)}}\,
\int_{-(R-0.5)}^{R-0.5}
\exp\!\Bigl(-\frac{(|I|-1)x^2}{2R(1-R)}\Bigr)\,dx.
\]
Combining the above approximation yields the desired estimate.
\end{proof}

\section{Details of the test function}
\label{appendix:compute_windows}

The test function $\phi_w(x,t)$ centered at $(x_j,t^n)$ is defined by
\begin{equation}
\phi_w(x,t)
=
C
\left(
1-
\left(
\frac{x-x_j}{m_x\Delta x}
\right)^2
\right)^{p_x}
\left(
1-
\left(
\frac{t-t^n}{m_t\Delta t}
\right)^2
\right)^{p_t},
\qquad
(x,t)\in\Omega_{w(x_j,t^n)},
\label{e:test_function}
\end{equation}
on the domain
\[
\Omega_{w(x_j,t^n)}
=
[x_j-m_x\Delta x,\;x_j+m_x\Delta x]
\times
[t^n-m_t\Delta t,\;t^n+m_t\Delta t].
\]
Here, $m_x$ and $m_t$ are the half-window sizes in space and time, respectively; $p_x$ and $p_t$ control the smoothness and boundary decay of the test function; and $C$ is a normalization constant. 

To determine $m_x$ and $p_x$, we follow \cite{tang2023weakident}: We apply the discrete Fourier transform to the data along the spatial direction, average the resulting magnitudes over all time levels, and form the cumulative spatial spectrum $E(k_x)$ by summing this time-averaged spectrum from frequency index $0$ to $k_x$, for $k_x = 0, \ldots, N_x$. We fit a two-piece linear model to $E(k_x)$ (as in Section~\ref{Subsec:regionofOne}) and denote the breakpoint as $k_x^*$.
We set
\[
m_x = \min\!\left\{\left\lfloor\frac{N_x-1}{2}\right\rfloor,\,\lceil m_x^* \rceil\right\}, \;\; \text{ and } \;\; p_x = \max\!\left\{\bar{\alpha} + 1,\;\left\lceil\frac{\log 10^{-10}}{\log\!\left(1-\left(1-\frac{1}{m_x}\right)^{\!2}\right)}\right\rceil\right\},
\]
where $m_x^*$ is the solution to
\[
\left(1-\left(1-\frac{1}{m_x}\right)^{\!2}\right)^{\!\frac{\pi^2(k_x^*)^2 m_x^2}{2\mathbb{N}_x^2} - \frac{3}{2}} = 10^{-10},
\]
and $\bar{\alpha} = \max_l\alpha_l$ is the maximum derivative order in the dictionary. With this choice, the spatial factor of $\phi_w$ is below $10^{-10}$ one grid step inside its boundary, and $p_x > \bar{\alpha}$, ensuring compatibility with integration by parts. The normalizing constant $C$ is chosen such that
\[
\int_{\Omega_{w(x_j,t^n)}} \phi_w(x,t)\,dx\,dt = 1,\quad \text{ and } \quad 
\phi_w(x,t) = 0 \text{ on } \partial\Omega_{w(x_j,t^n)}.
\]

\section{Details of experiment in Subsection~\ref{subsec:3pv_noise15}}
\label{append: noise15_stats}

In Subsection~\ref{subsec:3pv_noise15}, for the changing PDE with time varying coefficients~\eqref{eq:varmix3_pdes} under 15\% Gaussian noise, we presented an example in  Figure~\ref{fig:3varying_Noise1}, and reported the number of identification out of 50 independent trials in Table~\ref{tab:varmix3_support_freq}.  Here we provide more details of successful cases (27+9) in terms of  interval boundary and coefficient identification statistics. 

\begin{figure}
\centering  
(a) \\
\includegraphics[width=0.75\textwidth]{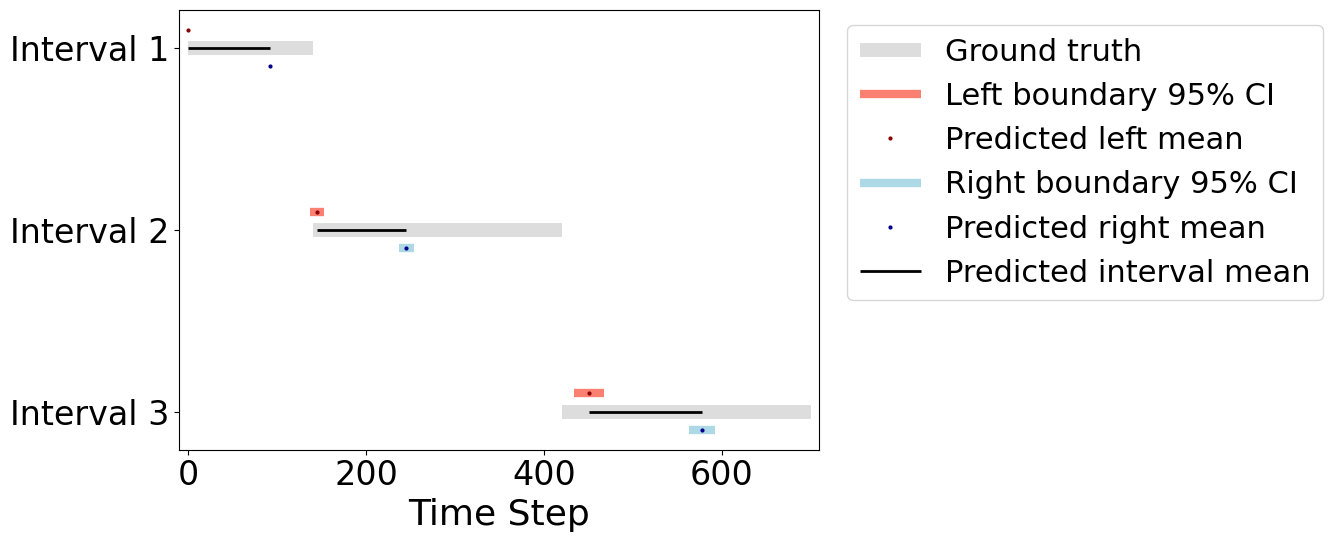}
\\ (b) \\
{\small
  \begin{tabular}{llcccccc}
    \toprule
    Interval & Stat & Left ($t$-idx) & Right ($t$-idx) & $ E_2$ & $ E_\infty$ & $ E_{\mathrm{res}}$ & $\mathrm{R (\%)}$ \\
    \midrule
    $[t^{0},t^{140})$
      & Mean & 0.000 & 92.185 & 0.0104 & 0.0104 & 0.0258 & 84.02 \\
      & SD & 0.000 & 1.210 & 0.002 & 0.002 & 0.0005 & 4.47 \\
      & SE & 0.000 & 0.233 & 0.0004 & 0.0004 & 0.0001 & 0.86 \\
      & 95\% CI lower & 0.00 & 91.73 & 0.0097 & 0.0097 & 0.0256 & 82.33 \\
      & 95\% CI upper & 0.00 & 92.64 & 0.0112 & 0.0112 & 0.0260 & 85.71 \\
    \midrule
    $[t^{140},t^{420})$
      & Mean & 145.074 & 245.482 & 0.0137 & 0.0356 & 0.1501 & 33.14 \\
      & SD & 20.077 & 21.707 & 0.0062 & 0.0054 & 0.0110 & 8.06 \\
      & SE & 3.864 & 4.177 & 0.0012 & 0.0010 & 0.0021 & 1.55 \\
      & 95\% CI lower & 137.45 & 237.17 & 0.0110 & 0.0336 & 0.1460 & 30.00 \\
      & 95\% CI upper & 152.70 & 253.79 & 0.0164 & 0.0376 & 0.1543 & 36.28 \\
    \midrule
    $[t^{420},t^{700})$
      & Mean & 450.444 & 578.185 & 0.0250 & 0.0748 & 0.2083 & 59.80 \\
      & SD & 44.628 & 38.315 & 0.0055 & 0.0078 & 0.0380 & 5.76 \\
      & SE & 8.589 & 7.374 & 0.0011 & 0.0015 & 0.0731 & 1.11 \\
      & 95\% CI lower & 433.62 & 564.75 & 0.0228 & 0.0719 & 0.1934 & 57.58 \\
      & 95\% CI upper & 467.27 & 591.63 & 0.0271 & 0.0778 & 0.2234 & 62.02 \\
    \bottomrule
  \end{tabular}
  }
  \caption{[Interval recovery statistics for Table~\ref{tab:varmix3_support_freq}, 27 of three interval case] (a) Identified regions of one equation with 95\% confidence intervals of their boundaries. (b) For each region of one equation, we present the mean of left and right boundary location, coefficient error and dominant ratio.}
  \label{fig:noise15_ci_vs_gt_3seg}
  \label{tab:varmix3_stats_3seg}
\end{figure}

\begin{figure}
  \centering
  (a) \\
  \includegraphics[width=0.75\textwidth]{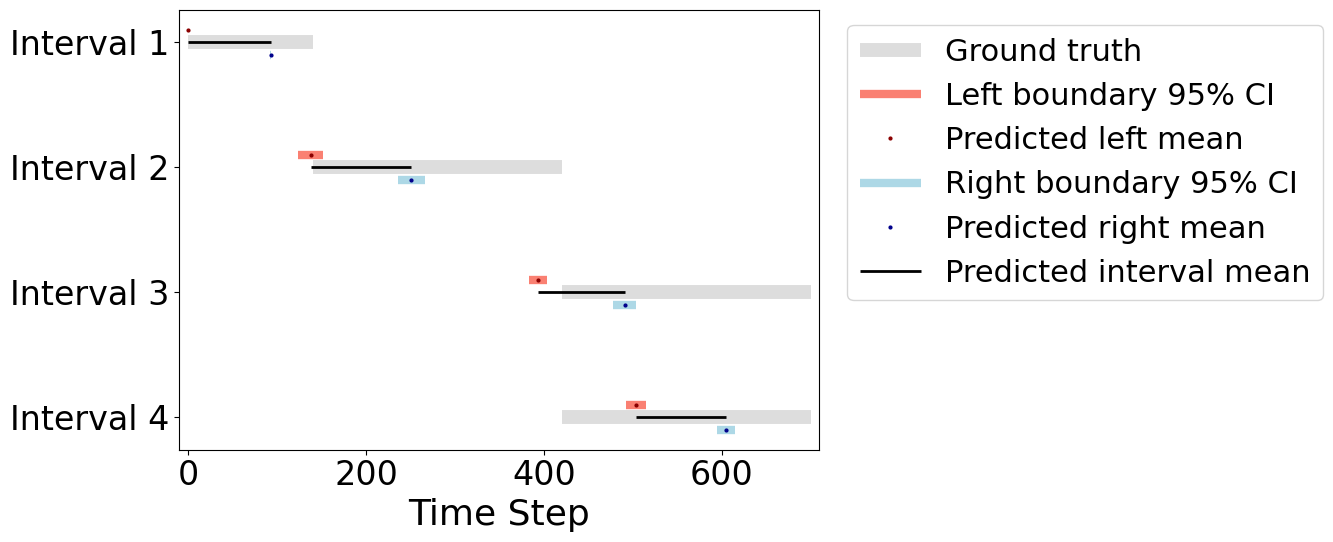}\\(b) \\
  {\small
  \begin{tabular}{llcccccc}
    \toprule
    Interval & Stat & Left ($t$-idx) & Right ($t$-idx) & $ E_2$ & $ E_\infty$ & $ E_{\mathrm{res}}$ & $\mathrm{R (\%)}$ \\
    \midrule
    $[t^{0},t^{140})$
      & Mean & 0.000 & 92.778 & 0.0115 & 0.0115 & 0.0260 & 84.51 \\
      & SD & 0.000 & 1.202 & 0.0021 & 0.0021 & 0.0005 & 2.84 \\
      & SE & 0.000 & 0.401 & 0.0007 & 0.0007 & 0.0002 & 0.95 \\
      & 95\% CI lower & 0.00 & 92.00 & 0.0101 & 0.0101 & 0.0257 & 83.56 \\
      & 95\% CI upper & 0.00 & 93.56 & 0.0130 & 0.0130 & 0.0263 & 85.46 \\
    \midrule
    $[t^{140},t^{420})$
      & Mean & 137.889 & 251.111 & 0.0133 & 0.0344 & 0.1487 & 33.23 \\
      & SD & 21.263 & 22.883 & 0.0037 & 0.0031 & 0.0110 & 7.99 \\
      & SE & 7.088 & 7.628 & 0.0012 & 0.0010 & 0.0037 & 2.66 \\
      & 95\% CI lower & 123.33 & 233.38 & 0.0108 & 0.0323 & 0.1415 & 30.57 \\
      & 95\% CI upper & 152.45 & 268.84 & 0.0157 & 0.0365 & 0.1559 & 35.89 \\
    \midrule
    $[t^{420},t^{700})$ (1st)
      & Mean & 393.889 & 490.889 & 0.0429 & 0.0876 & 0.2064 & 56.15 \\
      & SD & 15.235 & 19.896 & 0.0221 & 0.0202 & 0.0045 & 5.76 \\
      & SE & 5.078 & 6.632 & 0.0074 & 0.0067 & 0.0015 & 1.92 \\
      & 95\% CI lower & 383.07 & 476.34 & 0.0280 & 0.0745 & 0.2034 & 53.23 \\
      & 95\% CI upper & 404.71 & 505.44 & 0.0577 & 0.1008 & 0.2094 & 59.07 \\
    \midrule
    $[t^{420},t^{700})$ (2nd)
      & Mean & 503.556 & 605.222 & 0.0201 & 0.0793 & 0.2208 & 59.92 \\
      & SD & 17.111 & 15.147 & 0.0022 & 0.0078 & 0.0043 & 6.82 \\
      & SE & 5.704 & 5.049 & 0.0007 & 0.0026 & 0.0014 & 2.27 \\
      & 95\% CI lower & 491.65 & 595.32 & 0.0186 & 0.0738 & 0.2180 & 57.65 \\
      & 95\% CI upper & 515.46 & 615.12 & 0.0215 & 0.0849 & 0.2237 & 62.19 \\
    \bottomrule
  \end{tabular}
  }
  \caption{Interval recovery statistics for Table~\ref{tab:varmix3_support_freq}, 9 of four interval case] (a) Identified regions of one equation with 95\% confidence intervals of their boundaries. (b) For each region of one equation, we present the mean of left and right boundary location, coefficient error and dominant ratio.}
  \label{fig:noise15_ci_vs_gt_4seg}
  \label{tab:varmix3_stats_4seg}
\end{figure}

In Figure~\ref{fig:noise15_ci_vs_gt_3seg} (a), we present the 95\% confidence intervals of the identified boundaries, and in (b), we report the mean of left and right boundary, standard deviation (SD), standard error (SE), and 95\% CI lower and upper bounds, for the 27 trials that recovered the true supports \(\{u_x\}\),\(\{u_{xx}, (u^2)_x\}\), \(\{u_x, u_{xx}\}\) in three intervals.   We present coefficient errors and dominance ratios.  
In Figure~\ref{fig:noise15_ci_vs_gt_4seg}, we present the statistics of the 9 trials that recovered \(\{u_x\}\),\(\{u_{xx}, (u^2)_x\}\),\(\{u_x, u_{xx}\}\), \(\{u_x, u_{xx}\}\). 

\section{Additional experiment results}
\label{append:additional_experiments}

In the following, we present additional numerical experiments considered in Section~\ref{sec:numerical}.

\noindent\textbf{PDE changing between constant coefficient equations.}
In Figure~\ref{fig:transport_combined}, we consider a challenging scenario in which PDE changes between two equations across five intervals. The given data $\mathcal{D}$ is generated by solving a transport equation in which the diffusion term is intermittently activated:
\begin{equation}\label{add_transp}
  u_t = -u_x + c(t)\,u_{xx}
  \quad\text{on}\quad x\in[0,1],\;t\in[0,0.3),
\end{equation}
with $c(t)=0.05$ when $t \in [t^0,t^{100})$, activated to $0.10$ for $t \in [t^{200},t^{300})$ and $0.15$ for $t \in [t^{350},t^{500})$, otherwise, $c(t)=0$.
We use initial condition $u(x,0)=\exp\!\bigl(-(x+1)^2\bigr)$, and the numerical resolution is set to \(N_x \times N_t = 501 \times 501\).
\begin{figure}
\centering
{\small
\begin{tabular}{ll|ll}
\toprule
True Int. & SLW Int. & True Equation & SLW result \\ \midrule
$[t^{0},t^{100})$ & $[t^{0},t^{77})$ & $u_t = -u_x + 0.05\,u_{xx}$ & $u_t = -1.00000415\,u_x + 0.04999063\,u_{xx}$ \\
$[t^{100},t^{200})$ & $[t^{86},t^{174})$ & $u_t = -u_x$ & $u_t = -1.00074731\,u_x$ \\
$[t^{200},t^{300})$ & $[t^{192},t^{274})$ & $u_t = -u_x + 0.10\,u_{xx}$ & $u_t = -0.99999763\,u_x + 0.09999719\,u_{xx}$ \\
$[t^{300},t^{350})$ & $[t^{290},t^{324})$ & $u_t = -u_x$ & $u_t = -1.00011903\,u_x$ \\
$[t^{350},t^{500})$ & $[t^{341},t^{468})$ & $u_t = -u_x + 0.15\,u_{xx}$ & $u_t = -0.99999777\,u_x + 0.14999390\,u_{xx}$ \\ \bottomrule
\end{tabular}\\
\begin{tabular}{ccccccccc}
\toprule
Interval & Interval & Inclusion & Support & Support & $\mathrm{R (\%)}$ & $E_2$ & $E_\infty$ & $E_{\mathrm{res}}$ \\
TPR & PPV &  & TPR & PPV & & & & \\ \midrule
0.77 & 1.00 & 1.0 & 1.00 & 1.00 & 57.14 & $1.03\times10^{-5}$ & $1.87\times10^{-4}$ & $4.91\times10^{-6}$ \\
0.74 & 0.84 & 0.5 & 1.00 & 1.00 & 95.65 & $7.47\times10^{-4}$ & $7.47\times10^{-4}$ & $3.29\times10^{-3}$ \\
0.74 & 0.90 & 0.5 & 1.00 & 1.00 & 84.28 & $4.02\times10^{-6}$ & $2.86\times10^{-5}$ & $2.22\times10^{-6}$ \\
0.48 & 0.71 & 0.5 & 1.00 & 1.00 & 84.62 & $1.19\times10^{-4}$ & $1.19\times10^{-4}$ & $3.36\times10^{-4}$ \\
0.79 & 0.93 & 0.5 & 1.00 & 1.00 & 99.48 & $6.75\times10^{-6}$ & $4.15\times10^{-5}$ & $2.65\times10^{-6}$ \\ \bottomrule
\end{tabular}
}
\caption{PDE changing between two constant coefficient equations \eqref{add_transp}. The top table shows the SLW-Ident result, and the bottom table shows the accuracy measures for Interval 1 to 5 from top to bottom.}
\label{fig:transport_combined}
\end{figure}
The interval TPR in the first column is greater than 0.7 for all but Interval 4, and the interval PPV in the second column exceeds 0.8 for all but Interval 4. Interval 4 shows a slight shift in the identified region; even a small boundary offset has a noticeable impact on both TPR and PPV. Despite these boundary discrepancies, the support identification is accurate: all supports are recovered perfectly, with TPR = PPV = 1.00 across all detected intervals. The dominance ratios exceed 60\% in all but Interval 1, and the coefficient errors mostly stay at \(10^{-3}\) level or lower.

In Figure~\ref{fig:burgers_combined}, $\mathcal{D}$ is generated by solving a Burgers equation in which the advection term is intermittently activated:
\begin{equation}\label{add_burgers}
  u_t = 0.3\,u_{xx} - 0.5\,c(t)\,(u^2)_x
  \quad\text{on}\quad x\in[-5,5],\;t\in[0,10),
\end{equation}
with $c(t)=1.0$ when $t \in [t^0,t^{100})$, activated to $2.0$ for $t \in [t^{200},t^{300})$ and $1.5$ for $t \in [t^{350},t^{500})$, otherwise, $c(t)=0$.
We use initial condition $u(x,0)=\sin\!\bigl(\tfrac{\pi x}{4}\bigr)$, and the numerical resolution is set to \(N_x \times N_t = 501 \times 501\).
\begin{figure}
\centering
{\small
\begin{tabular}{ll|ll}
\toprule
True Int. & SLW Int. & True Equation & SLW result \\ \midrule
$[t^{0},t^{100})$ & $[t^{3},t^{79})$ & $u_t = 0.3\,u_{xx} - 0.5\,(u^2)_x$ & $u_t = 0.29998\,u_{xx} - 0.49942\,(u^2)_x$ \\
$[t^{100},t^{200})$ & $[t^{92},t^{174})$ & $u_t = 0.3\,u_{xx}$ & $u_t = 0.30000\,u_{xx}$ \\
$[t^{200},t^{300})$ & $[t^{191},t^{277})$ & $u_t = 0.3\,u_{xx} - 1.0\,(u^2)_x$ & $u_t = 0.30000\,u_{xx} - 0.99966\,(u^2)_x$ \\
$[t^{300},t^{350})$ & $[t^{290},t^{324})$ & $u_t = 0.3\,u_{xx}$ & $u_t = 0.30029\,u_{xx}$ \\
$[t^{350},t^{500})$ & $[t^{336},t^{468})$ & $u_t = 0.3\,u_{xx} - 0.75\,(u^2)_x$ & $u_t = 0.30000\,u_{xx} - 0.74762\,(u^2)_x$ \\ \bottomrule
\end{tabular}\\
\begin{tabular}{ccccccccc}
\toprule
Interval & Interval & Inclusion & Support & Support & $\mathrm{R (\%)}$ & $E_2$ & $E_\infty$ & $E_{\mathrm{res}}$ \\
TPR & PPV &  & TPR & PPV & & & & \\ \midrule
0.76 & 1.00 & 1.0 & 1.00 & 1.00 & 75.94 & $9.93\times10^{-4}$ & $1.16\times10^{-3}$ & $4.36\times10^{-6}$ \\
0.82 & 0.90 & 0.5 & 1.00 & 1.00 & 84.47 & $1.34\times10^{-5}$ & $1.34\times10^{-5}$ & $4.41\times10^{-5}$ \\
0.77 & 0.90 & 0.5 & 1.00 & 1.00 & 82.86 & $3.24\times10^{-4}$ & $3.38\times10^{-4}$ & $5.65\times10^{-6}$ \\
0.48 & 0.71 & 0.5 & 1.00 & 1.00 & 80.70 & $9.54\times10^{-4}$ & $9.54\times10^{-4}$ & $1.12\times10^{-3}$ \\
0.79 & 0.89 & 0.5 & 1.00 & 1.00 & 96.40 & $2.94\times10^{-3}$ & $3.17\times10^{-3}$ & $4.54\times10^{-6}$ \\ \bottomrule
\end{tabular}
}
\caption{PDE changing between two constant coefficient equations \eqref{add_burgers}. The top table shows the SLW-Ident result, and the bottom table shows the accuracy measures for Interval 1 to 5 from top to bottom.}
\label{fig:burgers_combined}
\end{figure}
The interval TPR is greater than 0.7 for all but Interval 4, and the interval PPV exceeds 0.8 for all but Interval 4. The method successfully identifies the governing PDE in each segment with perfect support recovery, Support TPR = Support PPV = 1.00 across all intervals. The dominance ratios $\mathrm{R}$ are generally high, confirming the stability of support selection within each detected interval. Coefficient estimation errors remain small, with average $E_2$ on the order of $10^{-3}$ or below.

In Figure~\ref{fig:allen_combined}, we present another example with a higher nonlinearity (containing a $u^3$ term) and two terms activated and deactivated concurrently. The given data $\mathcal{D}$ is generated by solving the following Allen--Cahn equation:
\begin{equation}\label{add_AC}
  u_t = c(t)\,u + 0.01\,u_{xx}  - c(t)\,u^3
  \quad\text{on}\quad x\in[-1,1],\;t\in[0,1],
\end{equation}
with $c(t)=0$ when $t \in [t^0,t^{140})$, activated to $0.01$ for $t \in [t^{140},t^{280})$ and $0.015$ for $t \in [t^{420},t^{490})$, otherwise, $c(t)=0$.
We use initial condition $u(x,0)=0.1\,\tanh(x)$, and the numerical resolution is set to \(N_x \times N_t = 701 \times 701\).
\begin{figure}
\centering
{\small
\begin{tabular}{ll|ll}
\toprule
True Int. & SLW Int. & True Equation & SLW result \\ \midrule
$[t^{0},t^{140})$ & $[t^{0},t^{118})$ & $u_t = 0.01\,u_{xx}$ & $u_t = 0.00999989\,u_{xx}$ \\
$[t^{140},t^{280})$ & $[t^{131},t^{256})$ & $u_t = 0.01\,u + 0.01\,u_{xx} - 0.01\,u^3$ & $u_t = 0.009972\,u + 0.01000035\,u_{xx}$ \\
$[t^{280},t^{420})$ & $[t^{267},t^{398})$ & $u_t = 0.01\,u_{xx}$ & $u_t = 0.00999918\,u_{xx}$ \\
$[t^{420},t^{490})$ & $[t^{412},t^{463})$ & $u_t = 0.015\,u + 0.01\,u_{xx} - 0.015\,u^3$ & $u_t = 0.01496417\,u + 0.01000045\,u_{xx}$ \\
$[t^{490},t^{700})$ & $[t^{477},t^{668})$ & $u_t = 0.01\,u_{xx}$ & $u_t = 0.0099994\,u_{xx}$ \\ \bottomrule
\end{tabular}\\
\begin{tabular}{ccccccccc}
\toprule
Interval & Interval & Inclusion & Support & Support & $\mathrm{R (\%)}$ & $E_2$ & $E_\infty$ & $E_{\mathrm{res}}$ \\
TPR & PPV &  & TPR & PPV & & & & \\ \midrule
0.84 & 1.00 & 1.0 & 1.00 & 1.00 & 86.58 & $1.10\times10^{-5}$ & $1.10\times10^{-5}$ & $1.92\times10^{-5}$ \\
0.89 & 1.00 & 0.5 & 0.67 & 1.00 & 76.98 & $5.77\times10^{-1}$ & $1.00$ & $5.41\times10^{-5}$ \\
0.84 & 0.90 & 0.5 & 1.00 & 1.00 & 93.62 & $8.42\times10^{-5}$ & $8.42\times10^{-5}$ & $1.15\times10^{-4}$ \\
0.61 & 0.84 & 0.5 & 0.67 & 1.00 & 55.00 & $6.40\times10^{-1}$ & $1.00$ & $6.87\times10^{-5}$ \\
0.85 & 0.93 & 0.5 & 1.00 & 1.00 & 98.58 & $6.17\times10^{-5}$ & $6.17\times10^{-5}$ & $7.24\times10^{-5}$ \\ \bottomrule
\end{tabular}
}
\caption{PDE changing between two constant coefficient equations \eqref{add_AC}. The top table shows the SLW-Ident result, and the bottom table shows the accuracy measures for Interval 1 to 5 from top to bottom.}
\label{fig:allen_combined}
\end{figure}
The interval TPR is greater than 0.8 for all but Interval 4, and the interval PPV exceeds 0.8 in all cases. In Intervals 2 and 4, the identified support omits the $u^3$ term, which leads to large coefficient error ($ E_2\approx 0.6$). To further assess recovery quality, we compute the dynamics error
\[
E_{\mathrm{dyn}} = \frac{1}{N_x\,|I_{\mathrm{gt}}|}
  \sum_{t^n\in I_{\mathrm{gt}}}\!
  \sum_{j=1}^{N_x}
    \bigl|\,U_{\mathrm{sim}}(x_j,t^n)-U(x_j,t^n)\bigr|,
\]
where $U_{\mathrm{sim}}$ is obtained by simulating the recovered PDE over the ground-truth interval $I_{\mathrm{gt}}$. The time-average dynamics error remains small ($ E_{\mathrm{dyn}}\approx 5.5\times10^{-5}$ and $7.0\times10^{-6}$), showing that the missing nonlinearity has little impact on the solution behavior.

In Figure~\ref{fig:fisher_combined}, the given data $\mathcal{D}$ is generated by solving a Fisher--KPP equation in which the reaction terms are intermittently activated:
\begin{equation} \label{add_fisher}
  u_t = c(t)\,u + u_{xx} - c(t)\,u^2
  \quad\text{on}\quad x\in[-10,10],\;t\in[0,10],
\end{equation}
with $c(t)=0.05$ when $t \in [t^0,t^{100})$, activated to $0.10$ for $t \in [t^{200},t^{300})$ and $0.15$ for $t \in [t^{350},t^{500})$, otherwise, $c(t)=0$.
We use initial condition
$
  u(x,0)
  = \exp\!\bigl(-(x+1)^2\bigr)
    +\sin\!\bigl(\tfrac{\pi x}{8}\bigr)
    +\cos\!\bigl(\tfrac{\pi x}{4}\bigr),$
and the numerical resolution is set to \(N_x \times N_t = 501 \times 501\).
\begin{figure}
\centering
{\small
\begin{tabular}{ll|ll}
\toprule
True Int. & SLW Int. & True Equation & SLW result \\ \midrule
$[t^{0},t^{100})$ & $[t^{0},t^{77})$ & $u_t = 0.05\,u + u_{xx} - 0.05\,u^2$ & $u_t = 0.04999\,u + 1.0000\,u_{xx} - 0.04999\,u^2$ \\
$[t^{100},t^{200})$ & $[t^{85},t^{177})$ & $u_t = u_{xx}$ & $u_t = 0.99639\,u_{xx}$ \\
$[t^{200},t^{300})$ & $[t^{191},t^{276})$ & $u_t = 0.10\,u + u_{xx} - 0.10\,u^2$ & $u_t = 0.09998\,u + 1.00000\,u_{xx} - 0.09998\,u^2$ \\
$[t^{300},t^{350})$ & $[t^{289},t^{328})$ & $u_t = u_{xx}$ & $u_t = 0.99807\,u_{xx}$ \\
$[t^{350},t^{500})$ & $[t^{341},t^{470})$ & $u_t = 0.15\,u + u_{xx} - 0.15\,u^2$ & $u_t = 0.14999\,u + 1.00000\,u_{xx} - 0.14999\,u^2$ \\ \bottomrule
\end{tabular}\\
\begin{tabular}{ccccccccc}
\toprule
Interval & Interval & Inclusion & Support & Support & $\mathrm{R (\%)}$ & $E_2$ & $E_\infty$ & $E_{\mathrm{res}}$ \\
TPR & PPV &  & TPR & PPV & & & & \\ \midrule
0.77 & 1.00 & 1.0 & 1.00 & 1.00 & 91.17 & $1.68\times10^{-6}$ & $2.42\times10^{-5}$ & $9.38\times10^{-8}$ \\
0.77 & 0.84 & 0.5 & 1.00 & 1.00 & 66.40 & $3.61\times10^{-3}$ & $3.61\times10^{-3}$ & $2.54\times10^{-3}$ \\
0.76 & 0.89 & 0.5 & 1.00 & 1.00 & 89.30 & $2.35\times10^{-5}$ & $1.67\times10^{-4}$ & $1.26\times10^{-6}$ \\
0.56 & 0.72 & 0.5 & 1.00 & 1.00 & 55.82 & $1.93\times10^{-3}$ & $1.93\times10^{-3}$ & $3.76\times10^{-4}$ \\
0.80 & 0.93 & 0.5 & 1.00 & 1.00 & 95.21 & $1.57\times10^{-5}$ & $7.22\times10^{-5}$ & $2.05\times10^{-6}$ \\ \bottomrule
\end{tabular}
}
\caption{PDE changing between two constant coefficient equations \eqref{add_fisher}. The top table shows the SLW-Ident result, and the bottom table shows the accuracy measures for Interval 1 to 5 from top to bottom.}
\label{fig:fisher_combined}
\end{figure}
The interval TPR in the first column is greater than 0.7 for all but Interval 4, and the interval PPV in the second column exceeds 0.8 for all but Interval 4. Despite these boundary discrepancies, the support identification is accurate: all supports are recovered perfectly, with TPR = PPV = 1.00 across all detected intervals. The dominance ratios exceed 60\% in all but Interval 4, and the coefficient errors mostly stay at \(10^{-3}\) level or lower.

After examining different settings in which changing supports are produced by activating and deactivating a single term, we consider the example in Figure~\ref{fig:mix2_combined}, in which the governing equation changes between completely different constant coefficient equations over five intervals:
\begin{align}
u_t &= 0.05\,u + u_{xx} - 0.05\,u^2, & t \in [t^{0},t^{120}), \nonumber\\
u_t &= u_x + 0.2\,u_{xx}, & t \in [t^{120},t^{210}), \nonumber\\
u_t &= -u_x, & t \in [t^{210},t^{300}), \label{eq:mix2}\\
u_t &= 0.1\,u + u_x + u_{xx} - 0.1\,u^2, & t \in [t^{300},t^{480}), \nonumber\\
u_t &= u_{xx}, & t \in [t^{480},t^{600}) \nonumber
\end{align}
on \(x\in[-8,8]\), \(t\in[0,10)\) with \(N_x=N_t=601\) and $
  u(x,0)
  = e^{-(x+1)^2}
    +\sin\!\bigl(\tfrac{\pi x}{8}\bigr)
    +\cos\!\bigl(\tfrac{\pi x}{4}\bigr).$
\begin{figure}
\centering
{\small
\begin{tabular}{ll|ll}
\toprule
True Int. & SLW Int. & True Equation & SLW result \\ \midrule
$[t^{0},t^{120})$ & $[t^{0},t^{97})$ & $u_t = 0.05\,u + u_{xx} - 0.05\,u^2$ & $u_t = 0.04996\,u + 1.00000\,u_{xx} - 0.04996\,u^2$ \\
$[t^{120},t^{210})$ & $[t^{107},t^{188})$ & $u_t = u_x + 0.2\,u_{xx}$ & $u_t = 0.99449\,u_x + 0.20207\,u_{xx}$ \\
$[t^{210},t^{300})$ & $[t^{197},t^{290})$ & $u_t = -u_x$ & $u_t = -1.00000\,u_x$ \\
$[t^{300},t^{480})$ & $[t^{344},t^{460})$ & $u_t = 0.1\,u + u_x + u_{xx} - 0.1\,u^2$ &
\begin{tabular}[t]{@{}l@{}}
$u_t = 0.09875\,u + 0.98747\,u_x$\\
${}+\,1.00017\,u_{xx} - 0.09870\,u^2$
\end{tabular} \\
$[t^{480},t^{600})$ & $[t^{468},t^{569})$ & $u_t = u_{xx}$ & $u_t = 1.00000\,u_{xx}$ \\ \bottomrule
\end{tabular}\\
\begin{tabular}{ccccccccc}
\toprule
Interval & Interval & Inclusion & Support & Support & $\mathrm{R (\%)}$ & $E_2$ & $E_\infty$ & $E_{\mathrm{res}}$ \\
TPR & PPV &  & TPR & PPV & & & & \\ \midrule
0.81 & 1.00 & 1.0 & 1.00 & 1.00 & 84.82 & $3.08\times10^{-4}$ & $8.61\times10^{-4}$ & $2.25\times10^{-4}$ \\
0.76 & 0.84 & 0.5 & 1.00 & 1.00 & 50.74 & $6.05\times10^{-3}$ & $1.32\times10^{-2}$ & $3.86\times10^{-4}$ \\
0.86 & 0.89 & 0.5 & 1.00 & 1.00 & 75.21 & $7.59\times10^{-2}$ & $7.59\times10^{-2}$ & $5.95\times10^{-2}$ \\
0.64 & 1.00 & 1.0 & 1.00 & 1.00 & 77.72 & $8.92\times10^{-3}$ & $1.31\times10^{-2}$ & $6.45\times10^{-5}$ \\
0.74 & 0.88 & 0.5 & 1.00 & 1.00 & 100.00 & $5.49\times10^{-4}$ & $5.49\times10^{-4}$ & $1.77\times10^{-3}$ \\ \bottomrule
\end{tabular}
}
\caption{PDE changing between five constant coefficient equations \eqref{eq:mix2}. The top table shows the SLW-Ident result, and the bottom table shows the accuracy measures for Interval 1 to 5 from top to bottom.}
\label{fig:mix2_combined}
\end{figure}
Interval TPR exceeds 0.7 in four of the five intervals and Interval PPV is above 0.7 in all cases.
The method successfully identifies the governing PDE in each segment with perfect support recovery, Support TPR = Support PPV = 1.00 across all intervals. The dominance ratios $\mathrm{R}$ are generally high, confirming the stability of support selection within each detected interval. Coefficient estimation errors remain small, with average $E_2$ on the order of $10^{-3}$ or below, indicating accurate recovery.

In Figure \ref{fig:transport_ramp_combined}, we consider a dataset $\mathcal{D}$ when a coefficient changes linearly during short transitions between constant values:
\begin{equation} \label{add_linear}
  u_t = -u_x + c(t)\,u_{xx}
  \quad\text{on}\quad x\in[0,1],\;t\in[0,0.3),
\end{equation}
with $c(t) = 0.05$ for $t^{0} \le t < t^{50}$, a short transition to zero by $c(t) = -5\,(t-t^{50})+0.05$ when $t^{50} \le t < t^{60}$, and $c(t)=0$ for $t^{60} \le t < t^{140}$. Another short transition $c(t) = 10\,(t-t^{140})$ when $t^{140} \le t < t^{150}$, and $c(t) = 0.10$ for $t^{150} \le t < t^{300}$. We use initial condition $u(x,0)=\exp\!\bigl(-(x+1)^2\bigr)$, and \(N_x \times N_t = 301 \times 301\).
Since the transition intervals are short, we take the ground-truth intervals of one equation to be
$[t^0,t^{50}):\{u_x,u_{xx}\}$,
$[t^{60},t^{140}):\{u_x\}$,
$[t^{150},t^{300}):\{u_x,u_{xx}\}$.
\begin{figure}
\centering
{\small
\begin{tabular}{ll|ll}
\toprule
True Int. & SLW Int. & True Equation & SLW result \\ \midrule
$[t^{0},t^{50})$ & $[t^{0},t^{33})$ & $u_t = -u_x + 0.05\,u_{xx}$ & $u_t = -1.00005\,u_x + 0.04987\,u_{xx}$ \\
$[t^{60},t^{140})$ & $[t^{45},t^{118})$ & $u_t = -u_x$ & $u_t = -1.00005\,u_x$ \\
$[t^{150},t^{300})$ & $[t^{144},t^{269})$ & $u_t = -u_x + 0.10\,u_{xx}$ & $u_t = -1.00000\,u_x + 0.10000\,u_{xx}$ \\ \bottomrule
\end{tabular}\\
\begin{tabular}{ccccccccc}
\toprule
Interval & Interval & Inclusion & Support & Support & $\mathrm{R (\%)}$ & $E_2$ & $E_\infty$ & $E_{\mathrm{res}}$ \\
TPR & PPV &  & TPR & PPV & & & & \\ \midrule
0.66 & 1.00 & 1.0 & 1.00 & 1.00 & 42.06 & $1.44\times10^{-4}$ & $2.69\times10^{-3}$ & $4.36\times10^{-5}$ \\
0.73 & 0.79 & 0.5 & 1.00 & 1.00 & 98.27 & $4.52\times10^{-5}$ & $4.53\times10^{-5}$ & $1.20\times10^{-3}$ \\
0.79 & 0.95 & 0.5 & 1.00 & 1.00 & 100.00 & $9.49\times10^{-8}$ & $1.97\times10^{-7}$ & $7.11\times10^{-7}$ \\ \bottomrule
\end{tabular}
}
\caption{PDE changing linearly during short transitions between constant coefficient equations \eqref{add_linear}. The top table shows the SLW-Ident result, and the bottom table shows the accuracy measures for Interval 1 to 3 from top to bottom.}
\label{fig:transport_ramp_combined}
\end{figure}
SLW-Ident localizes all three intervals of one equation with interval TPR greater than 0.6 and interval PPV greater than 0.7 even with the presence of short linear transition regions.
The reduction in interval TPR values is due to identified intervals being shorter and shifted. Since the coefficient varies continuously, transition point detection becomes more difficult. However, these boundary offsets do not affect support identification or coefficient recovery within the intervals of one equation. In all detected intervals, the true support is recovered exactly (TPR = PPV = 1.00), and the estimated coefficients remain highly accurate.

In Figure~\ref{fig:afkpp_ramp_combined}, we consider a dataset $\mathcal{D}$ when a coefficient changes linearly during short transitions between constant values:
\begin{equation}\label{add_short}
  u_t = c(t)\,u + u_x + u_{xx} - c(t)\,u^2
  \quad\text{on}\quad x\in[-8,8],\;t\in[0,10),
\end{equation}
with $c(t) = 0$ for $t^{0} \le t < t^{60}$, a short transition $c(t) = 2\,(t-t^{60})$ when $t^{60} \le t < t^{63}$, and $c(t)=0.2$ for $t^{63} \le t < t^{180}$. Another short transition $c(t) = -2\,(t-t^{180})+0.2$ when $t^{180} \le t < t^{183}$, and $c(t) = 0$ for $t^{183} \le t < t^{300}$. We use initial condition
$
  u(x,0)
  = \exp\!\bigl(-(x+1)^2\bigr)
    +\sin\!\bigl(\tfrac{\pi x}{8}\bigr)
    +\cos\!\bigl(\tfrac{\pi x}{4}\bigr),$
and \(N_x \times N_t = 301 \times 301\).
Since the transition intervals are short, we take the ground-truth intervals of one equation to be
$[t^0,t^{60}):\{u_x,u_{xx}\}$,
$[t^{63},t^{150}):\{u,u_x,u_{xx},u^2\}$,
$[t^{183},t^{300}):\{u_x,u_{xx}\}$.
\begin{figure}
\centering
{\small
\begin{tabular}{ll|ll}
\toprule
True Int. & SLW Int. & True Equation & SLW result \\ \midrule
$[t^{0},t^{60})$ & $[t^{0},t^{35})$ & $u_t = u_x + u_{xx}$ & $u_t = 0.99999\,u_x + 0.99996\,u_{xx}$ \\
$[t^{63},t^{150})$ & $[t^{62},t^{154})$ & $u_t = 0.2\,u + u_x + u_{xx} - 0.2\,u^2$ &
\begin{tabular}[t]{@{}l@{}}
$u_t = 0.20000\,u + 0.99999\,u_x$\\
${}+\,0.99999\,u_{xx} - 0.19999\,u^2$
\end{tabular} \\
$[t^{183},t^{300})$ & $[t^{170},t^{266})$ & $u_t = u_x + u_{xx}$ & $u_t = 1.00019\,u_x + 0.99877\,u_{xx}$ \\ \bottomrule
\end{tabular}\\
\begin{tabular}{ccccccccc}
\toprule
Interval & Interval & Inclusion & Support & Support & $\mathrm{R (\%)}$ & $E_2$ & $E_\infty$ & $E_{\mathrm{res}}$ \\
TPR & PPV &  & TPR & PPV & & & & \\ \midrule
0.58 & 1.00 & 1.0 & 1.00 & 1.00 & 85.00 & $2.67\times10^{-5}$ & $3.70\times10^{-5}$ & $1.67\times10^{-5}$ \\
0.78 & 1.00 & 0.5 & 1.00 & 1.00 & 78.32 & $2.77\times10^{-1}$ & $2.00\times10^{0}$ & $6.39\times10^{-7}$ \\
0.74 & 1.00 & 0.5 & 1.00 & 1.00 & 100.00 & $8.80\times10^{-4}$ & $1.23\times10^{-3}$ & $2.76\times10^{-4}$ \\ \bottomrule
\end{tabular}
}
\caption{PDE changing linearly during short transitions between constant coefficient equations \eqref{add_short}. Notice the gap between the true intervals. Recovery of equations is accurate.}
\label{fig:afkpp_ramp_combined}
\end{figure}
SLW-Ident localizes all three intervals of one equation with interval TPR greater than 0.5 and interval PPV greater than 0.8 even with the presence of short linear transition regions.
In all detected intervals, the true support is recovered exactly (TPR = PPV = 1.00), and the estimated coefficients remain highly accurate within the intervals of one equation.

\noindent\textbf{Noisy data of changing PDEs with time varying coefficients.}
Besides the $15\%$ noise example in Subsection~\ref{subsec:3pv_noise15}, we consider another case with 10\% Gaussian noise (NSR = 0.10).
We generate the given data $\mathcal{D}$ by PDE changing among three constant coefficient equations:
\begin{align}
u_t &= u_x, & t \in [t^{0},t^{233}), \nonumber\\
u_t &= -u_x + 0.05\,u_{xx}, & t \in [t^{234},t^{466}), \label{eq:noise10}\\
u_t &= 0.1\,u_{xx}, & t \in [t^{467},t^{700}) \nonumber
\end{align}
on \(x\in[0,1]\), \(t\in[0,0.3)\), with \(N_x=N_t=701\) and initial condition
\[
  u(x,0)
  = \sin\!\bigl(\tfrac{4\pi x}{0.7}\bigr)^3
    \cos\!\bigl(\tfrac{\pi x}{0.7}\bigr)
    +0.8\,e^{-\frac{(x-0.2)^2}{2(0.05)^2}}
    +0.6\,e^{-\frac{(x-0.5)^2}{2(0.08)^2}}
\quad(\text{for }x<0.7),
\]
and $0$ otherwise.
We present the result of a representative run in Figure~\ref{fig:transport_heat_noise10}.
\begin{figure}
\centering
{\small
\begin{tabular}{ll|ll}
\toprule
True Int. & SLW Int. & True Equation & SLW result \\ \midrule
$[t^{0},t^{233})$ & $[t^{0},t^{174})$ & $u_t = u_x$ & $u_t = 1.00085\,u_x$ \\
$[t^{234},t^{466})$ & $[t^{209},t^{411})$ & $u_t = -u_x + 0.05\,u_{xx}$ & $u_t = -0.99245\,u_x + 0.05006\,u_{xx}$ \\
$[t^{467},t^{700})$ & $[t^{448},t^{579})$ & $u_t = 0.1\,u_{xx}$ & $u_t = 0.10004\,u_{xx}$ \\ \bottomrule
\end{tabular}\\
\begin{tabular}{ccccccccc}
\toprule
Interval & Interval & Inclusion & Support & Support & $\mathrm{R (\%)}$ & $E_2$ & $E_\infty$ & $E_{\mathrm{res}}$ \\
TPR & PPV &  & TPR & PPV & & & & \\ \midrule
0.75 & 0.75 & 1.0 & 1.00 & 1.00 & 33.8 & $6.91\times10^{-3}$ & $6.91\times10^{-3}$ & $2.54\times10^{-2}$ \\
0.72 & 0.72 & 0.5 & 1.00 & 1.00 & 27.3 & $1.23\times10^{-2}$ & $2.65\times10^{-2}$ & $3.73\times10^{-2}$ \\
0.68 & 0.68 & 0.5 & 1.00 & 1.00 & 9.6 & $2.67\times10^{-2}$ & $2.67\times10^{-2}$ & $1.04\times10^{-1}$ \\ \bottomrule
\end{tabular}
}
\caption{PDE changing between three constant coefficient equations with 10\% Gaussian noise (NSR = 0.10). The top table shows the SLW-Ident result, and the bottom table shows the accuracy measures for Interval 1 to 3 from top to bottom.}
\label{fig:transport_heat_noise10}
\end{figure}
SLW-Ident recovers each interval's support perfectly; although identified intervals are shorter than the true ones, as shown by low Interval TPR values, the identified intervals lie in the true intervals, with large interval PPV.

\begin{figure}
  \centering
  (a) \\
    {\small
    \begin{tabular}{lc}
\toprule
    \textbf{Support Combination (Intervals 1--3)}
      & \textbf{Number of Identification} \\
    \midrule
    \(\{u_x\},\;\{u_x,u_{xx}\},\;\{u_{xx}\}\)
      & 32 \\
    \(\{u_x,(u^2)_x,(u^3)_x,(u^4)_x\},\;\{u_x\},\;\{u_{xx}\}\)
      & 16 \\
    \(\{u_x,(u^2)_x,(u^3)_x,(u^4)_x\},\;\{u_x\},\;\{u_x,\;u^2,\;(u^2)_x,\;u^3,\;(u^3)_x,\;(u^3)_{xx},\;(u^4)_x\}\)
      & 1 \\
    \(\{u_x,(u^2)_x,(u^3)_x,(u^4)_x\},\;\{u_x\},\;\{u_{xx},(u^2)_{xx},(u^3)_{xx},(u^4)_{xx}\}\)
      & 1 \\
    \bottomrule
  \end{tabular}
  }\\
  (b) \\
\includegraphics[width=0.75\textwidth]{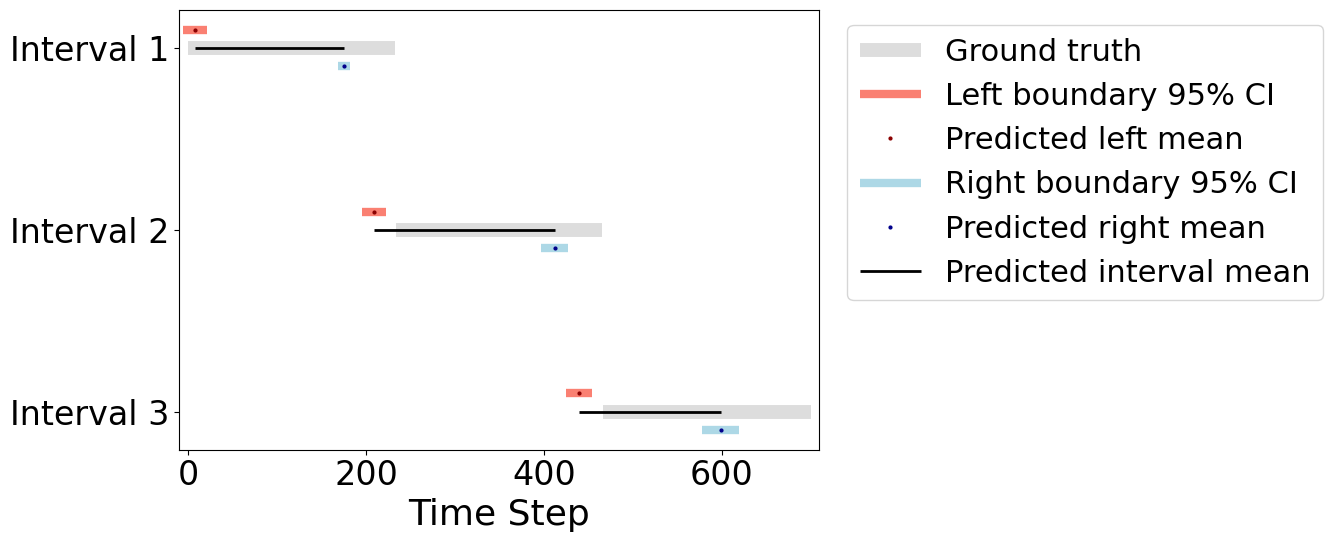}\\ (c) \\
{\small
  \begin{tabular}{llcccccc}
    \toprule
    Interval & Stat & Left ($t$-idx) & Right ($t$-idx) & $ E_2$ & $ E_\infty$ & $ E_{\mathrm{res}}$ & $\mathrm{R (\%)}$ \\
    \midrule
    $[t^{0},t^{233})$
      & Mean & 7.875 & 174.940 & 0.0073 & 0.0073 & 0.0252 & 32.45 \\
      & SD & 6.719 & 3.435 & 0.0011 & 0.0011 & 0.0006 & 3.03 \\
      & SE & 1.188 & 0.607 & 0.0002 & 0.0002 & 0.0001 & 0.54 \\
      & 95\% CI lower & 5.55 & 173.75 & 0.0069 & 0.0069 & 0.0249 & 31.39 \\
      & 95\% CI upper & 10.20 & 176.13 & 0.0077 & 0.0077 & 0.0254 & 33.51 \\
    \midrule
    $[t^{234},t^{466})$
      & Mean & 209.06 & 412.34 & 0.0130 & 0.0274 & 0.0374 & 26.07 \\
      & SD & 6.705 & 7.815 & 0.0027 & 0.0020 & 0.0005 & 3.39 \\
      & SE & 1.185 & 1.381 & 0.0005 & 0.0004 & 0.0001 & 0.60 \\
      & 95\% CI lower & 206.74 & 409.63 & 0.0120 & 0.0267 & 0.0372 & 24.89 \\
      & 95\% CI upper & 211.38 & 415.05 & 0.0139 & 0.0281 & 0.0376 & 27.25 \\
    \midrule
    $[t^{467},t^{700})$
      & Mean & 439.59 & 599.03 & 0.0279 & 0.0279 & 0.1077 & 13.01 \\
      & SD & 7.219 & 10.718 & 0.0033 & 0.0033 & 0.0035 & 3.14 \\
      & SE & 1.276 & 1.895 & 0.0006 & 0.0006 & 0.0006 & 0.56 \\
      & 95\% CI lower & 437.09 & 595.31 & 0.0267 & 0.0267 & 0.1065 & 11.91 \\
      & 95\% CI upper & 442.09 & 602.75 & 0.0290 & 0.0290 & 0.1089 & 14.11 \\
    \bottomrule
  \end{tabular}
}
\caption{PDE changing between three constant coefficient equations \eqref{eq:noise10} with 10\% Gaussian noise (NSR = 0.10), 50 independent experiments. In (a), the first column shows regions of one equation and the support identified in each interval. (b) Correctly identified 32 cases, regions of one equation with 95\% confidence intervals of their boundaries. (c) Recovery statistics for the 32 correct cases: mean, standard deviation (SD), standard error (SE),  95\% CI lower and upper bounds for boundary locations,  coefficient error and the dominant ratio. }
  \label{fig:noise10_ci_vs_gt}
\end{figure}

We repeated this experiment 50 times, each time with an independent noise realization. Figure~\ref{fig:noise10_ci_vs_gt} summarizes the identified support patterns, showing that in 64\% of the trials the proposed model finds the correct support
\(\{\!u_x\},\;\{u_x,u_{xx}\},\;\{u_{xx}\}\).
\begin{table}
  \centering

  \label{tab:noise10_support_freq}
\end{table}
Restricting attention to the 32 trials that recovered the true supports, Figure~\ref{fig:noise10_ci_vs_gt} overlays the 95\% confidence intervals of the identified boundaries and reports the corresponding mean, sample standard deviation (SD), standard error (SE), and 95\% CI lower and upper bounds for each quantity.
Overall, among the successful trials, the identified intervals remain well-localized and the coefficient errors stay small, confirming that SLW-Ident recovers not only the correct supports but also accurate interval boundaries and coefficients under 10\% noise.

\end{document}